\RequirePackage{ifpdf}
\ifpdf 
\documentclass[pdftex]{sigma}
\else
\documentclass{sigma}
\fi

\usepackage[all]{xy}
\usepackage[mathscr]{eucal}

\usepackage{stmaryrd}
\newcommand{\slar}{\shortleftarrow}

\newenvironment{enm}%
{\begin{enumerate}
\setlength{\itemsep}{0em}
}%
{\end{enumerate}}
\newenvironment{itm}%
{\begin{itemize}
\setlength{\itemsep}{0em}
}%
{\end{itemize}}
{\begin{description}
\setlength{\itemsep}{0em}
}%
{\end{description}}

\numberwithin{equation}{section}

\newtheorem{Theorem}{Theorem}[section]
\newtheorem{Proposition}[Theorem]{Proposition}
\newtheorem{Lemma}[Theorem]{Lemma}
\theoremstyle{definition}
\newtheorem{Definition}[Theorem]{Definition}
\newtheorem{Remark}[Theorem]{Remark}
\newtheorem{Example}[Theorem]{Example}

\newcommand{\thmref}[1]{Theo\-rem~\ref{#1}}
\newcommand{\secref}[1]{Section~\ref{#1}}
\newcommand{\lemref}[1]{Lemma~\ref{#1}}
\newcommand{\propref}[1]{Proposition~\ref{#1}}

\newcommand{\remref}[1]{Remark~\ref{#1}}

\newcommand{\Z}{\mathbb{Z}} 
\newcommand{\R}{\mathbb{R}} 
\newcommand{\C}{\mathbb{C}} 
\newcommand{\bbP}{\mathbb{P}} 

\newcommand{\GL}{\operatorname{GL}} 
\newcommand{\gl}{\operatorname{\mathfrak{gl}}}
\newcommand{\g}{\mathfrak{g}} 
\newcommand{\h}{\mathfrak{h}} 
\newcommand{\frb}{\mathfrak{b}}
\newcommand{\Lie}{\operatorname{Lie}} 
\newcommand{\ad}{\operatorname{ad}}

\newcommand{\End}{\operatorname{End}}
\newcommand{\Aut}{\operatorname{Aut}}
\newcommand{\Hom}{\operatorname{Hom}}
\newcommand{\Ker}{\operatorname{Ker}}
\newcommand{\range}{\operatorname{Im}}
\newcommand{\Coker}{\operatorname{Coker}}

\newcommand{\rank}{\operatorname{rank}}
\newcommand{\tr}{\operatorname{tr}}
\newcommand{\unit}{\mathrm{Id}}
\newcommand{\prj}{\operatorname{pr}}

\newcommand{\Rep}{\operatorname{Rep}}
\newcommand{\bM}{\operatorname{\mathbf{M}}}
\newcommand{\bA}{\mathbf{A}} 
\newcommand{\bD}{\mathbf{D}} 
\newcommand{\bC}{\mathbf{C}} 

\newcommand{\bV}{\mathbf{V}}

\newcommand{\bv}{\mathbf{v}}
\newcommand{\bw}{\mathbf{w}}
\newcommand{\bS}{\mathbf{S}}
\newcommand{\bT}{\mathbf{T}}

\newcommand{\lset}[2]%
{\left\{ \, \left. #1 \hspace{0.25em} \right| \, #2 \, \right\} }
\newcommand{\rset}[2]%
{\left\{ \, #1 \, \left| \hspace{0.25em} #2\right. \, \right\} }
\newcommand{\ov}{\overline}

\newcommand{\res}{\operatorname*{res}}
\newcommand{\mc}{\operatorname{\it mc}}

\newcommand{\ved}{\mathbf{d}}

\newcommand{\Mirr}{\operatorname{\mathcal{M}}^\mathrm{irr}}
\newcommand{\Mset}{\operatorname{\mathcal{M}}^\mathrm{set}}
\newcommand{\bO}{\mathcal{O}}
\newcommand{\irr}{\mathrm{irr}}
\newcommand{\st}{\mathrm{s}}
\newcommand{\sumfrac}[2]{\genfrac{}{}{0pt}{2}{#1}{#2}}
\newcommand{\hquad}{\hspace{.5em}}

\newcommand{\Qst}{\operatorname{\mathscr{N}}^\mathrm{s}}
\newcommand{\Qset}{\operatorname{\mathscr{N}}^\mathrm{set}}
\newcommand{\vin}{\mathsf{t}} 
\newcommand{\vout}{\mathsf{s}} 
\newcommand{\QQ}{\mathsf{Q}} 
\newcommand{\frF}{\mathscr{F}} 

\newcommand{\rd}{\mathop{\mathrm{d}\!}\mathstrut}

\begin{document}

\allowdisplaybreaks

\renewcommand{\PaperNumber}{087}

\FirstPageHeading

\ShortArticleName{Quiver Varieties with Multiplicities}

\ArticleName{Quiver Varieties with Multiplicities, Weyl Groups \\
of Non-Symmetric Kac--Moody Algebras, \\
and Painlev\'e Equations}

\Author{Daisuke YAMAKAWA~$^{\dag\ddag}$}

\AuthorNameForHeading{D. Yamakawa}

\Address{$^\dag$~Centre de math\'ematiques Laurent Schwartz,
\'Ecole polytechnique,\\
\hphantom{$^\dag$}~CNRS UMR 7640, ANR S\'EDIGA,
91128 Palaiseau Cedex, France}

\Address{$^\ddag$~Department of Mathematics, Graduate School of Science,
Kobe University,\\
\hphantom{$^\ddag$}~Rokko, Kobe 657-8501, Japan}
\EmailD{\href{mailto:yamakawa@math.kobe-u.ac.jp}{yamakawa@math.kobe-u.ac.jp}}

\ArticleDates{Received March 19, 2010, in f\/inal form October 18, 2010;  Published online October 26, 2010}

\Abstract{To a f\/inite quiver equipped with a positive integer on each of its vertices,
we associate a holomorphic symplectic manifold having some parameters.
This coincides with Nakajima's quiver variety
with no stability parameter/framing
if the integers attached on the vertices are all equal to one.
The construction of ref\/lection functors for quiver varieties
are generalized to our case,
in which these relate to simple ref\/lections in the Weyl group of
some symmetrizable, possibly non-symmetric Kac--Moody algebra.
The moduli spaces of
meromorphic connections on the rank 2 trivial bundle over the Riemann sphere
are described as our manifolds.
In our picture, the list of Dynkin diagrams for Painlev\'e equations
is slightly dif\/ferent from (but equivalent to) Okamoto's.
}


\Classification{53D30; 16G20; 20F55; 34M55}

\section{Introduction}\label{sec:intro}

First, we brief\/ly explain our main objects in this article.
Let
\begin{itemize}\itemsep=0pt
\item $\QQ$ be a quiver, i.e., a directed graph,
with the set of vertices $I$
(our quivers are always assumed to be f\/inite
and have no arrows joining a vertex with itself);
\item $\ved =(d_i)_{i \in I} \in \Z_{>0}^I$ be
a collection of positive integers indexed by the vertices.
\end{itemize}
We think of each number $d_i$ as the `multiplicity' of the vertex $i \in I$,
so the pair $(\QQ,\ved)$ as a~`quiver with multiplicities'.
In this article,
we associate to such $(\QQ,\ved)$ a holomorphic symplectic manifold
$\Qst_{\QQ,\ved}(\lambda,\bv)$ having parameters
\begin{itemize}\itemsep=0pt
\item $\lambda=(\lambda_i(z))_{i \in I}$, where
$\lambda_i(z)=\lambda_{i,1}z^{-1} + \lambda_{i,2}z^{-2}
+ \cdots + \lambda_{i,d_i}z^{-d_i} \in z^{-d_i}\C[z]/\C[z]$;
\item $\bv=(v_i)_{i \in I} \in \Z_{\geq 0}^I$,
\end{itemize}
and call it the {\em quiver variety with multiplicities},
because if $d_i=1$ for all $i \in I$,
it then coincides with (the stable locus of)
Nakajima's quiver variety
$\mathfrak{M}^\mathrm{reg}_\zeta(\bv,\bw)$ \cite{MR1302318} with
\[
\bw=0 \in \Z_{\geq 0}^I, \qquad
\zeta=(\zeta_\R,\zeta_\C)=\bigl( 0,(\lambda_{i,1})_{i \in I} \bigr)
\in \sqrt{-1}\R^I \times \C^I.
\]
As in the case of quiver variety,
$\Qst_{\QQ,\ved}(\lambda,\bv)$ is def\/ined as
a holomorphic symplectic quotient
with respect to some algebraic group action (see \secref{sec:q-multi}).
However, the group used here is {\em non-reductive}
unless $d_i=1$ or $v_i=0$ for all $i \in I$.
Therefore a number of basic facts
in the theory of holomorphic symplectic quotients
(e.g.\ the hyper-K\"ahler quotient description) cannot be applied to our
$\Qst_{\QQ,\ved}(\lambda,\bv)$,
and for the same reason,
they seem to provide new geometric objects relating to quivers.

The def\/inition of $\Qst_{\QQ,\ved}(\lambda,\bv)$
is motivated by the theory of Painlev\'e equations.
It is known due to Okamoto's work
\cite{MR854008,MR916698,MR914314,MR927186} that
all Painlev\'e equations except the f\/irst one
have (extended) af\/f\/ine Weyl group symmetries;
see the table below, where $P_J$ denotes the Painlev\'e equation
of type $J$ ($J=\text{II, III, \dots ,VI}$).
\[
\hfill
{\renewcommand{\arraystretch}{1.2}
\begin{array}{|c||c|c|c|c|c|}
\hline
\text{Equations} & P_\mathrm{VI} & P_\mathrm{V} & P_\mathrm{IV} & P_\mathrm{III} & P_\mathrm{II}
\\ \hline
\text{Symmetries} & D_4^{(1)} & A_3^{(1)} & A_2^{(1)} & C_2^{(1)} & A_1^{(1)}
\\ \hline
\end{array}}
\hfill
\]
On the other hand, each of them is known to
govern an isomonodromic deformation
of rank two meromorphic connections on $\bbP^1$~\cite{MR625446};
the number of poles and the pole orders of connections remain
unchanged during the deformation, and are determined from
(if we assume that the connections have only `unramif\/ied' singularities)
the type of the Painlev\'e equation (see e.g.\ \cite{0902.1702}).
See the table below, where $k_1 + k_2 + \cdots + k_n$
means that the connections in the deformation
have $n$ poles of order $k_i$, $i=1,2, \dots ,n$ and no other poles.
\[
\hfill
{\renewcommand{\arraystretch}{1.2}
\begin{array}{|c||c|c|c|c|c|}
\hline
\text{Equations} & P_\mathrm{VI} & P_\mathrm{V} & P_\mathrm{IV} & P_\mathrm{III} & P_\mathrm{II}
\\ \hline
\text{Connections} & 1+1+1+1 & 2+1+1 & 3+1 & 2+2 & 4
\\ \hline
\end{array}}
\hfill
\]
Roughly speaking, we thus have a non-trivial correspondence
between some Dynkin diagrams and
rank two meromorphic connections.

In fact, such a relationship can be understood in terms of quiver varieties
except in the case of $P_\mathrm{III}$.
Crawley-Boevey \cite{MR1980997} described
the moduli spaces of Fuchsian systems
(i.e., meromorphic connections on the trivial bundle over $\bbP^1$
having only simple poles)
as quiver varieties associated with `star-shaped' quivers.
In particular, the moduli space of rank two Fuchsian systems
having exactly four poles are
described as a quiver variety of type $D_4^{(1)}$,
which is consistent with the above correspondence for $P_\mathrm{VI}$.
The quiver description in the cases of $P_\mathrm{II}$, $P_\mathrm{IV}$
and $P_\mathrm{V}$ was obtained by Boalch\footnote{His description in the case of $P_\mathrm{V}$
is based on the work of Harnad~\cite{MR1309553}.} \cite{0806.1050};
more generally, he proved that the moduli spaces of meromorphic connections
on the trivial bundle over $\bbP^1$
having one higher order pole (and possibly simple poles)
are quiver varieties.

A remarkable point is that
their quiver description provides
Weyl group symmetries of the moduli spaces\footnote{Actually in each Painlev\'e case, this action reduces to
an action of the corresponding f\/inite Weyl group,
which together with `Schlesinger transformations'
give the full symmetry; see \cite[Section~6]{MR2500553}.}
at the same time,
because for any quiver,
the associated quiver varieties are known to have such symmetry.
This is generated by the so-called {\em reflection functors}
(see \thmref{thm:intro-ref} below),
whose existence was f\/irst announced by Nakajima
(see \cite[Section~9]{MR1302318}, where
he also gave its geometric proof in some important cases),
and then shown by several researchers including himself
\cite{MR1620538,MR1775358,MR2023313,MR1428466}.

The purpose of quiver varieties with multiplicities
is to generalize their description to the case of $P_\mathrm{III}$;
the starting point is the following observation
(see \propref{prop:moduli} for a further generalized, precise statement;
see also Remarks \ref{rem:stable-bundle} and \ref{rem:moduli}):
\begin{Proposition}\label{prop:intro-moduli}
Consider a `star-shaped quiver of length one'
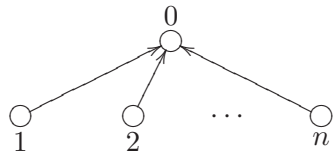
\[
\hfill
\begin{xy}
\ar@{} (25,0);(30,0)_*=0{\cdots}
\ar@{<-} (20,10)*++!D{0}*\cir<4pt>{}="A";(0,0)*++!U{1}*\cir<4pt>{}
\ar@{<-} "A";(15,0)*++!U{2}*\cir<4pt>{}
\ar@{<-} "A";(40,0)*++!U{n}*\cir<4pt>{}
\end{xy}
\hfill
\]
Here the set of vertices is $I=\{\, 0,1, \dots ,n \,\}$.
Take multiplicities $\ved \in \Z_{>0}^I$ with $d_0=1$ and
set $\bv \in \Z_{\geq 0}^I$ by $v_0=2$, $v_i=1$ $(i=1,\dots ,n)$.
Then $\Qst_{\QQ,\ved}(\lambda,\bv)$ is
isomorphic to the moduli space of stable meromorphic connections
on the rank two trivial bundle over $\bbP^1$ having $n$ poles of
order~$d_i$, $i=1, \dots ,n$ of prescribed formal type.
\end{Proposition}

On the other hand, to any quiver with multiplicities,
we associate a generalized Cartan mat\-rix~$\bC$
as follows:
\[
\bC = 2 \unit - \bA \bD,
\]
where $\bA$ is the adjacency matrix of the underlying graph,
namely, the matrix whose $(i,j)$ entry is
the number of edges joining $i$ and $j$,
and $\bD$ is the diagonal matrix with entries
given by the multiplicities $\ved$.
It is symmetrizable as $\bD \bC$ is symmetric, but may be not symmetric.

Now as stated below, our
quiver varieties with multiplicities admit ref\/lection functors;
this is the main result of this article.
\begin{Theorem}[see \secref{sec:reflection}]\label{thm:intro-ref}
For any quiver with multiplicities $(\QQ,\ved)$,
there exist linear maps
\[
s_i \colon \ \Z^I \to \Z^I, \qquad
r_i \colon \ \bigoplus_{i \in I}\bigl( z^{-d_i}\C[z]/\C[z] \bigr) \to
\bigoplus_{i \in I}\bigl( z^{-d_i}\C[z]/\C[z] \bigr)
\qquad (i \in I)
\]
generating actions of the Weyl group of the associated Kac--Moody algebra,
such that for any $(\lambda,\bv)$ and $i \in I$ with $\lambda_{i,d_i} \neq 0$,
one has a natural symplectomorphism
\[
\frF_i \colon \ \Qst_{\QQ,\ved}(\lambda,\bv) \xrightarrow{\simeq}
\Qst_{\QQ,\ved}(r_i(\lambda),s_i(\bv)).
\]
If $d_i=1$ for all $i \in I$, then
the maps $\frF_i$ coincide with the reflection functors.
\end{Theorem}

In the case of star-shaped quivers,
the original ref\/lection functor at the central vertex
can be interpreted in terms of
Katz's {\em middle convolution}~\cite{MR1366651}
for Fuchsian systems (see \cite[Appendix~A]{MR2107041}).
A similar assertion also holds
in the situation of \propref{prop:intro-moduli};
the map $\frF_0$ at the central vertex~$0$
can be interpreted in terms of
the `generalized middle convolution' \cite{0808.0699,0911.3863}
(see \secref{subsec:middle}).

For instance, consider the star-shaped quivers with multiplicities
given below
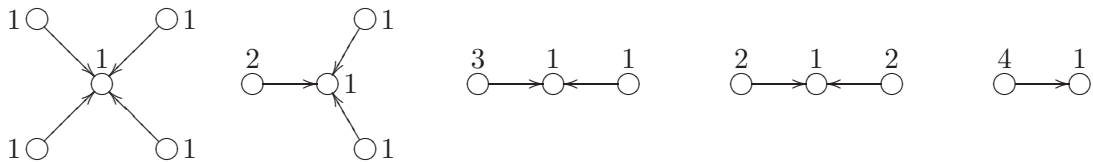
\begin{gather*}
\begin{xy}
  \ar@{<-} (0,0) *++!D{1} *\cir<4pt>{}="A";
    (8.61,-8.61)   *++!L{1} *\cir<4pt>{}
  \ar@{<-} "A";(-8.61,-8.61) *++!R{1} *\cir<4pt>{}
  \ar@{<-} "A";(-8.61,8.61) *++!R{1} *\cir<4pt>{}
  \ar@{<-} "A";(8.61,8.61) *++!L{1} *\cir<4pt>{}
  \ar@{->} (20,0) *++!D{2} *\cir<4pt>{};
    (30,0)   *++!L{1} *\cir<4pt>{}="B"
  \ar@{<-} "B";(35,8.61) *++!L{1} *\cir<4pt>{}
  \ar@{<-} "B";(35,-8.61) *++!L{1} *\cir<4pt>{}
  \ar@{->} (50,0) *++!D{3} *\cir<4pt>{};
    (60,0)   *++!D{1} *\cir<4pt>{}="C"
  \ar@{<-} "C";(70,0) *++!D{1} *\cir<4pt>{}
  \ar@{->} (85,0) *++!D{2} *\cir<4pt>{};
    (95,0)   *++!D{1} *\cir<4pt>{}="D"
  \ar@{<-} "D";(105,0) *++!D{2} *\cir<4pt>{}
  \ar@{->} (120,0) *++!D{4} *\cir<4pt>{};
    (130,0)   *++!D{1} *\cir<4pt>{}
\end{xy}
\end{gather*}
\propref{prop:intro-moduli} says that
the associated $\Qst_{\QQ,\ved}(\lambda,\bv)$
with a particular choice of $\bv$ give
the moduli spaces for
$P_\mathrm{VI}$, $P_\mathrm{V}$,
$P_\mathrm{IV}$, $P_\mathrm{III}$ and $P_\mathrm{II}$, respectively.
On the other hand, the associated Kac--Moody algebras are respectively given by\footnote{We follow Kac \cite{MR1104219} for
the notation of (twisted) af\/f\/ine Lie algebras.}
\[
D_4^{(1)}, \quad A_5^{(2)}, \quad D_4^{(3)}, \quad C_2^{(1)}, \quad A_2^{(2)}.
\]
Interestingly, this list of Kac--Moody algebras is dif\/ferent from
the table given before;
however we can clarify the relationship between
our description and Boalch's by using
a sort of `shifting trick' established by him
(see \secref{subsec:shifting}).
This trick, which may be viewed as a geometric phenomenon arising from
the `normalization of the leading coef\/f\/icient
in the principal part of the connection at an irregular singular point',
connects two quiver varieties with multiplicities associated to
{\em different} quivers with multiplicities;
more specif\/ically, we prove the following:
\begin{Theorem}[see \secref{sec:normalization}]\label{thm:intro-normalization}
Suppose that a quiver with multiplicities $(\QQ,\ved)$
has a pair of vertices $(i,j)$ such that
\[
d_i >1, \qquad d_j =1, \qquad
a_{ik}=a_{ki}=\delta_{jk} \quad \text{for any}\hquad k \in I,
\]
where $\bA=(a_{ij})$ is the adjacency matrix of the underlying graph.
Then it determines another quiver with multiplicities
$(\check{\QQ},\check{\ved})$
and a map $(\lambda,\bv) \mapsto (\check{\lambda},\check{\bv})$
between parameters such that
the following holds:
if $\lambda_{i,d_i} \neq 0$, then $\Qst_{\QQ,\ved}(\lambda,\bv)$
and $\Qst_{\check{\QQ},\check{\ved}}(\check{\lambda},\check{\bv})$
are symplectomorphic to each other.
\end{Theorem}

We call the transformation
$(\QQ,\ved) \mapsto (\check{\QQ},\check{\ved})$,
whose precise def\/inition is given in \secref{subsec:normalization},
the {\em normalization}.
Using this theorem,
we can translate the above list of Dynkin diagrams into the original one
(see \secref{subsec:example}).

There is a close relationship between
two Kac--Moody root systems connected via the normalization
(see \secref{subsec:weyl}).
In particular, we have the following relation between
the Weyl groups $W$, $\check{W}$ associated to
$(\QQ,\ved)$, $(\check{\QQ},\check{\ved})$:
\[
W \simeq \check{W} \rtimes \Z/2\Z,
\]
where the semidirect product is taken with respect to
some Dynkin automorphism of order~2
(such a Dynkin automorphism canonically
exists by the def\/inition of normalization).
For instance,
in the cases of $P_\mathrm{V}$, $P_\mathrm{IV}$ and $P_\mathrm{II}$,
we have
\begin{gather*}
W\big(A_5^{(2)}\big) \simeq W\big(A_3^{(1)}\big) \rtimes \Z/2\Z, \\
W\big(D_4^{(3)}\big) \simeq W\big(A_2^{(1)}\big) \rtimes \Z/2\Z, \\
W\big(A_2^{(2)}\big) \simeq W\big(A_1^{(1)}\big) \rtimes \Z/2\Z,
\end{gather*}
which mean that our list of Dynkin diagrams for Painlev\'e equations
is a variant of Okamoto's obtained by (partially) extending the Weyl groups.

\section{Preliminaries}\label{sec:pre}

In this section
we brief\/ly recall the def\/inition of Nakajima's quiver variety~\cite{MR1302318}.

\subsection{Quiver}\label{subsec:quiver}

Recall that a (f\/inite) {\em quiver} is a
quadruple $\QQ=(I,\Omega,\vout,\vin)$ consisting of
two f\/inite sets~$I$,~$\Omega$
(the set of {\em vertices}, resp.\ {\em arrows})
and two maps
$\vout, \vin \colon \Omega \to I$
(assigning to each arrow
its {\em source}, resp.\ {\em target}).
Throughout this article, for simplicity,
we assume that {\em our quivers $\QQ$
have no arrow $h \in \Omega$ with $\vout(h)=\vin(h)$}.

For given $\QQ$, we denote by $\ov{\QQ}=(I,\ov{\Omega},\vout,\vin)$
the quiver obtained from $\QQ$ by reversing the orientation of each arrow;
the set $\ov{\Omega}=\{\, \ov h \mid h \in \Omega \,\}$
is just a copy of $\Omega$,
and $\vout(\ov h):=\vin(h)$, $\vin(\ov h):=\vout(h)$
for $h \in \Omega$.
We set $H := \Omega \sqcup \ov{\Omega}$,
and extend the map $\Omega \to \ov{\Omega}$, $h \mapsto \ov h$ to
an involution of $H$ in the obvious way.
The resulting quiver $\QQ + \ov{\QQ} =(I,H,\vout,\vin)$
is called the {\em double} of $\QQ$.

The {\em underlying graph} of $\QQ$,
which is obtained by forgetting the orientation of each arrow,
determines a symmetric matrix $\bA=(a_{ij})_{i,j \in I}$,
called the {\em adjacency matrix}, as follows:
\[
a_{ij} := \sharp\{\,\text{edges joining $i$ and $j$}\,\}
= \sharp\{\, h \in H \mid \vout(h)=i,\ \vin(h)=j \,\}.
\]

Let $\bV=\bigoplus_{i \in I} V_i$ be a nonzero f\/inite-dimensional
$I$-graded $\C$-vector space.
A representation of $\QQ$ over $\bV$ is an element of the vector space
\[
\Rep_\QQ(\bV):= \bigoplus_{h \in \Omega}
\Hom_\C(V_{\vout(h)},V_{\vin(h)}),
\]
and its {\em dimension vector} is given by
$\bv := \dim \bV \equiv (\dim V_i)_{i \in I}$.
Isomorphism classes of representations of $\QQ$ with dimension vector $\bv$
just correspond to orbits in $\Rep_\QQ(\bV)$
with respect to the action of the group
$\GL(\bV):=\prod\limits_{i \in I} \GL_\C(V_i)$ given by
\[
g=(g_i) \colon (B_h)_{h \in \Omega} \longmapsto
\bigl( g_{\vin(h)}B_h g_{\vout(h)}^{-1} \bigr)_{h \in \Omega},
\qquad g \in \GL(\bV).
\]
We denote the Lie algebra of $\GL(\bV)$ by $\gl(\bV)$; explicitly,
$\gl(\bV):=\bigoplus_{i \in I} \gl_\C(V_i)$.
For $\zeta =(\zeta_i)_{i \in I} \in \C^I$,
we denote its image under the natural map $\C^I \to \gl(\bV)$
by $\zeta\,\unit_\bV$,
and also use the same letter $\zeta\,\unit_\bV$
for $\zeta \in \C$ via the diagonal embedding $\C \hookrightarrow \C^I$.
Note that the central subgroup
$\C^\times \simeq
\{\,\zeta\,\unit_\bV \mid \zeta \in \C^\times\,\} \subset \GL(\bV)$
acts trivially on $\Rep_\QQ(\bV)$,
so we have the induced action of the quotient group
$\GL(\bV)/\C^\times$.

Let $B=(B_h)_{h \in \Omega} \in \Rep_\QQ(\bV)$.
An $I$-graded subspace $\bS=\bigoplus_{i \in I}S_i$ of $\bV$ is said to be
{\em $B$-invariant} if
$B_h(S_{\vout(h)}) \subset S_{\vin(h)}$ for all $h \in \Omega$.
If $\bV$ has no $B$-invariant subspace except $\bS=0, \bV$,
then $B$ is said to be {\em irreducible}.
Schur's lemma\footnote{One can apply Schur's lemma
thanks to the following well-known fact:
the category of representations of~$\QQ$
is equivalent to that of an algebra~$\C \QQ$,
the so-called {\em path algebra}; see e.g.~\cite{Etingof}.}
implies that the stabilizer of each irreducible~$B$
is just the central subgroup $\C^\times \subset \GL(\bV)$,
and a standard fact in
Mumford's geometric invariant theory
\cite[Corollary~2.5]{MR1304906} (see also \cite{MR1315461})
implies that the action of $\GL(\bV)/\C^\times$ on the subset
$\Rep_\QQ^\irr(\bV)$ consisting of all irreducible representations
over $\bV$ is proper.

\subsection{Quiver variety}\label{subsec:q-variety}

Suppose that a quiver $\QQ$
and a nonzero f\/inite-dimensional $I$-graded $\C$-vector space
$\bV = \bigoplus_{i \in I}V_i$ are given.
We set
\[
\bM_\QQ(\bV) := \Rep_{\QQ+\ov{\QQ}}(\bV)
= \Rep_\QQ(\bV) \oplus \Rep_{\ov{\QQ}}(\bV),
\]
and regard it as the cotangent bundle of $\Rep_\QQ(\bV)$
by using the trace pairing.
Introducing the function
\[
\epsilon \colon H \to \{ \pm 1 \}, \qquad \epsilon(h):=
\begin{cases}
\hphantom{-}1 & \text{for}\hquad h \in \Omega,\\
-1& \text{for}\hquad h \in \ov{\Omega},
\end{cases}
\]
we can write the canonical symplectic form on $\bM_\QQ(\bV)$ as
\[
\omega := \sum_{h \in \Omega} \tr \rd B_h \wedge \rd B_{\ov{h}}
= \frac12 \sum_{h \in H} \epsilon(h) \tr \rd B_h \wedge \rd B_{\ov{h}},
\qquad (B_h)_{h \in H} \in \bM_\QQ(\bV).
\]
The natural $\GL(\bV)$-action on $\bM_\QQ(\bV)$ is Hamiltonian
with respect to $\omega$ with the moment map
\begin{gather}\label{eq:moment}
\mu=(\mu_i)_{i \in I} \colon \bM_\QQ(\bV) \to \gl(\bV), \qquad
\mu_i(B)=\sum_{\sumfrac{h \in H:}{\vin(h)=i}} \epsilon (h) B_h B_{\ov h}
\end{gather}
vanishing at the origin,
where we identify $\gl(\bV)$
with its dual using the trace pairing.

\begin{Definition}\label{dfn:stable}
A point $B \in \bM_\QQ(\bV)$ is said to be {\em stable}
if it is irreducible as a representation of $\QQ + \ov{\QQ}$.
\end{Definition}

For a $\GL(\bV)$-invariant Zariski closed subset $Z$ of $\bM_\QQ(\bV)$,
let $Z^\st$ be the subset of all stable points in $Z$.
It is a $\GL(\bV)$-invariant Zariski open subset of $Z$,
on which the group $\GL(\bV)/\C^\times$ acts
freely and properly.

\begin{Definition}\label{dfn:q-variety}
For $\zeta \in \C^I$ and $\bv \in \Z_{\geq 0}^I \setminus \{ 0 \}$,
taking an $I$-graded $\C$-vector space $\bV$ with $\dim \bV = \bv$
we def\/ine
\[
\Qst_\QQ(\zeta,\bv) := \mu^{-1}(-\zeta\,\unit_\bV)^\st / \GL(\bV),
\]
which we call the {\em quiver variety}.
\end{Definition}

\begin{Remark}\label{rem:nakajima}
In Nakajima's notation (see \cite{MR1302318} or \cite{MR2023313}),
$\Qst_\QQ(\zeta,\bv)$ is denoted by
$\mathfrak{M}^\mathrm{reg}_{(0,\zeta)}(\bv,0)$.
\end{Remark}

\section{Quiver variety with multiplicities}\label{sec:q-multi}

\subsection{Def\/inition}\label{subsec:dfn}

For a positive integer $d$, we set
\[
R_d := \C[[z]]/z^d\C[[z]], \qquad R^d := z^{-d}\C[[z]]/\C[[z]].
\]
The $\C$-algebra $R_d$ has a typical basis $\{\, z^{d-1}, \dots ,z,1 \,\}$,
with respect to which
the multiplication by $z$ in $R_d$ is represented by
the nilpotent single Jordan block
\[
J_d := \begin{pmatrix}
                        0     & 1 &        & 0 \\
                              & 0 & \ddots &   \\
                              &   & \ddots & 1   \\
                        0     &   &        & 0 \end{pmatrix}
\in \End(\C^d)=\End_\C(R_d).
\]
The vector space $R^d$ may be identif\/ied with
the $\C$-dual $R_d^* =\Hom_\C (R_d,\C)$ of $R_d$ via the pairing
\[
R_d \otimes_\C R^d \to \C, \qquad (f,g) \mapsto
\res_{z=0} \bigl( f(z)g(z) \bigr).
\]

For a f\/inite-dimensional $\C$-vector space $V$, we set
\[
\g_d(V) :=
\gl(V) \otimes_\C R_d = \gl(V)[[z]]/z^d \gl(V)[[z]].
\]
Note that $\g_d(V)$ is naturally isomorphic to
$\End_{R_d}(V \otimes_\C R_d)$ as an $R_d$-module;
hence it is the Lie algebra of the complex algebraic group
\[
G_d(V) := \Aut_{R_d}(V \otimes_\C R_d) \simeq
\lset{g(z) = \sum_{k=0}^{d-1}g_k z^k \in \g_d(V)}{\det g_0 \neq 0}.
\]
The inverse element of $g(z) \in G_d(V)$ is given by
taking modulo $z^d \gl(V)[[z]]$ of
the formal inverse $g(z)^{-1} \in \gl(V)[[z]]$.
The adjoint action of $g(z)$ is described as
\[
(g \cdot \xi)(z)= g(z)\xi(z)g(z)^{-1} \mod z^d \gl(V)[[z]],
\qquad \xi(z) \in \g_d(V).
\]
Using the above $R_d^* \simeq R^d$ and the trace pairing,
we always identify the $\C$-dual $\g_d^*(V)$ of $\g_d(V)$ with
$\gl(V) \otimes_\C R^d = z^{-d}\gl(V)[[z]]/\gl(V)[[z]]$.
Then the coadjoint action of $g(z) \in G_d(V)$ is also described as
\[
(g \cdot \eta)(z)= g(z)\eta(z)g(z)^{-1} \mod \gl(V)[[z]],
\qquad \eta(z)=\sum_{k=1}^d \eta_k z^{-k} \in \g^*_d(V).
\]

The natural inclusion
$\g_d(V) \hookrightarrow \End_\C(V \otimes_\C R_d)
=\End_\C(V) \otimes_\C \End_\C(R_d)$ is represented by
\[
\xi(z)=\sum_{k=0}^{d-1} \xi_k z^k \longmapsto
\sum_{k=0}^{d-1} \xi_k \otimes J_d^k,
\]
whose image is just the centralizer of $\unit_V \otimes J_d$.
Accordingly, its transpose
can be written as
\[
\gl_\C(V \otimes_\C R_d) \simeq \gl_\C (V \otimes_\C R_d)^*
\to \g_d^*(V),\qquad
X \mapsto \sum_{k=1}^d \tr_{R_d}
\bigl[ X \big(\unit_V \otimes J_d^{k-1}\big) \bigr]  z^{-k},
\]
where $\tr_{R_d} \colon \End_\C(V \otimes_\C R_d)
= \End_\C(V)\otimes_\C \End_\C(R_d) \to \End_\C(V)$
denotes the trace of the $\End_\C(R_d)$-part.

Now suppose that a quiver $\QQ$
and a collection of positive integers $\ved =(d_i)_{i \in I} \in \Z_{>0}^I$
are given.
We call the pair $(\QQ,\ved)$ as a {\em quiver with multiplicities}
and $d_i$ as the {\em multiplicity} of the vertex $i$.
Set
\[
R_\ved := \bigoplus_{i \in I} R_{d_i}, \qquad
R^\ved := \bigoplus_{i \in I} R^{d_i},
\]
and for a nonzero f\/inite-dimensional $I$-graded $\C$-vector space
$\bV=\bigoplus_{i \in I} V_i$,
set
\begin{gather*}
\bV_\ved \equiv \bV \otimes_\C R_\ved
:= \bigoplus_{i \in I} V_i \otimes_\C R_{d_i}, \\
\bM_{\QQ,\ved}(\bV) := \bM_\QQ(\bV_\ved)
= \bigoplus_{h\in H}
\Hom_\C \big(V_{\vout(h)} \otimes_\C R_{d_{\vout(h)}},
V_{\vin(h)} \otimes_\C R_{d_{\vin(h)}}\big), \\
G_\ved(\bV) := \prod_{i \in I} G_{d_i}(V_i), \qquad
\g_\ved(\bV):= \bigoplus_{i \in I} \g_{d_i}(V_i).
\end{gather*}
The group $G_\ved(\bV)$ naturally acts on $\bM_{\QQ,\ved}(\bV)$ as a
subgroup of $\GL(\bV_\ved)$.
Note that the subgroup
$\C^\times \subset \GL(\bV_\ved)$
is contained in $G_\ved(\bV)$ and acts trivially on $\bM_{\QQ,\ved}(\bV)$.
As in the case of $\gl(\bV)$,
for $\lambda =(\lambda_i(z))_{i \in I} \in R^\ved$
we denote its image under the natural map
$R^\ved=\g^*_\ved(\C^I) \to \g^*_\ved(\bV)$
by~$\lambda\,\unit_\bV$.

Let $\omega$ be the canonical symplectic form on $\bM_{\QQ,\ved}(\bV)$;
\[
\omega = \frac12 \sum_{h \in H} \epsilon(h) \tr \rd B_h \wedge \rd B_{\ov{h}},
\qquad (B_h)_{h \in H} \in \bM_{\QQ,\ved}(\bV).
\]
Then the $G_\ved(\bV)$-action is Hamiltonian
whose moment map $\mu_\ved$ is given by the composite of
the $\GL(\bV_\ved)$-moment map
$\mu =(\mu_i) \colon \bM_{\QQ,\ved}(\bV) \to \gl(\bV_\ved)$
(see \eqref{eq:moment} for the def\/inition)
and the transpose
$\prj=(\prj_i)$ of the inclusion
$\g_\ved(\bV) \hookrightarrow \gl(\bV_\ved)$;
\begin{gather*}
\mu_\ved  =(\mu_{\ved,i})_{i \in I}
\colon \bM_{\QQ,\ved}(\bV) \to \g_\ved^*(\bV), \\
\mu_{\ved,i}(B)  := \prj_i \circ\, \mu_i (B)
= \sum_{k=1}^{d_i} \sum_{\sumfrac{h \in H:}{\vin(h)=i}}
\epsilon(h) \tr_{R_{d_i}}
\bigl[ B_h B_{\ov{h}} N_i^{k-1} \bigr]  z^{-k},
\end{gather*}
where $N_i := \unit_{V_i} \otimes J_{d_i}$.

\begin{Definition}\label{dfn:conf-stable}
A point $B \in \bM_{\QQ,\ved}(\bV)$
is said to be {\em stable} if
$\bV_\ved$ has no nonzero proper $B$-invariant subspace
$\bS=\bigoplus_{i \in I} S_i$ such that
$S_i \subset V_i \otimes_\C R_{d_i}$
is an $R_{d_i}$-submodule for each $i \in I$.
\end{Definition}

The above stability can be interpreted in terms of
the irreducibility of representations of a~quiver.
Letting $\widetilde{\Omega} := \Omega \sqcup \{ \ell_i \}_{i \in I}$
and extending the maps $\vout$, $\vin$ to $\widetilde{\Omega}$
by $\vout(\ell_i)=\vin(\ell_i)=i$,
we obtain a new quiver
$\widetilde{\QQ}=(I,\widetilde{\Omega},\vout,\vin)$.
Consider the vector space
\[
\Rep_{\widetilde{\QQ}+\ov{\QQ}}(\bV_\ved)
\simeq \bM_{\QQ,\ved}(\bV) \oplus \gl(\bV_\ved)
\]
associated with the quiver
$\widetilde{\QQ}+\ov{\QQ}=(I,\widetilde{\Omega}\sqcup \ov{\Omega},\vout,\vin)$.
Note that in the above def\/inition,
a vector subspace $S_i \subset V_i \otimes R_{d_i}$
is an $R_{d_i}$-submodule
if and only if it is invariant under the action of
$N_i=\unit_{V_i} \otimes J_{d_i}$,
which corresponds to the multiplication by $z$.
Thus letting
\begin{gather}\label{eq:embed}
\iota \colon \ \bM_{\QQ,\ved}(\bV) \hookrightarrow
\Rep_{\widetilde{\QQ}+\ov{\QQ}}(\bV_\ved), \qquad
B \mapsto (B, (N_i)_{i \in I}),
\end{gather}
we see that a point $B \in \bM_{\QQ,\ved}(\bV)$ is
stable if and only if its image $\iota(B)$
is irreducible as a~representation of $\widetilde{\QQ}+\ov{\QQ}$.

For a $G_\ved(\bV)$-invariant Zariski closed subset $Z$ of
$\bM_{\QQ,\ved}(\bV)$,
let $Z^\st$ be the subset of all stable points in $Z$.

\begin{Proposition}
The group $G_\ved(\bV)/\C^\times$ acts
freely and properly on $Z^\st$.
\end{Proposition}

\begin{proof}
Note that
the closed embedding $\iota$ def\/ined in \eqref{eq:embed}
is equivariant under the action of
$G_\ved(\bV) \subset \GL(\bV_\ved)$.
Hence the freeness of
the $G_\ved(\bV)/\C^\times$-action on $Z^\st$
follows from that of the $\GL(\bV_\ved)/\C^\times$-action on
$\Rep_{\widetilde{\QQ}+\ov{\QQ}}(\bV_\ved)^\irr$ and
\[
\iota(\bM_{\QQ,\ved}(\bV)^\st) =
\iota(\bM_{\QQ,\ved}(\bV)) \cap \Rep_{\widetilde{\QQ}+\ov{\QQ}}(\bV_\ved)^\irr,
\]
which we have already checked.
Furthermore, the above implies that the embedding
$Z^\st \hookrightarrow \Rep_{\widetilde{\QQ}+\ov{\QQ}}(\bV_\ved)^\irr$
induced from $\iota$ is closed.
Consider the following commutative diagram:
\[
\xymatrix{
G_\ved(\bV)/\C^\times \times Z^\st
\ar[r] \ar[d] \ar@{}[rd]|{\circlearrowright} & Z^\st \ar[d] \\
\GL(\bV_\ved)/\C^\times \times
\Rep_{\widetilde{\QQ}+\ov{\QQ}}(\bV_\ved)^\irr \ar[r] &
\Rep_{\widetilde{\QQ}+\ov{\QQ}}(\bV_\ved)^\irr,
}
\]
where the vertical arrows are the maps induced from $\iota$,
and the horizonal arrows are the action maps
$(g,x) \mapsto g \cdot x$.
Since the bottle horizontal arrow is proper
and both vertical arrows are closed,
the properness of the top horizontal arrow follows from
well-known basic properties of proper maps
(see e.g.\ \cite[Corollary~4.8]{Hartshorne}).
\end{proof}

\begin{Definition}\label{dfn:q-multi}
For $\lambda \in R^\ved$ and $\bv \in \Z_{\geq 0}^I \setminus \{ 0 \}$,
taking an $I$-graded $\C$-vector space $\bV$ with $\dim \bV = \bv$
we def\/ine
\[
\Qst_{\QQ,\ved}(\lambda,\bv) :=
\mu_\ved^{-1}(-\lambda\,\unit_\bV)^\st / G_\ved(\bV),
\]
which we call the {\em quiver variety with multiplicities}.
\end{Definition}

We also use the following set-theoretical quotient:
\[
\Qset_{\QQ,\ved}(\lambda,\bv) :=
\mu_\ved^{-1}(-\lambda\,\unit_\bV) / G_\ved(\bV).
\]

It is clear from the def\/inition that if $(\QQ,\ved)$ is multiplicity-free,
i.e., $d_i=1$ for all $i \in I$,
then $\Qst_{\QQ,\ved}(\lambda,\bv)$ coincides with the ordinal quiver variety
$\Qst_\QQ(\zeta,\bv)$ with $\zeta_i = \res\limits_{z=0}\lambda_i(z)$.
Even when $\ved$ is non-trivial, for simplicity,
we often refer to $\Qst_{\QQ,\ved}(\lambda,\bv)$ just as the `quiver variety'.

\subsection{Properties}\label{subsec:property}

Here we introduce some basic properties of quiver varieties
with multiplicities.

First, we associate a symmetrizable Kac--Moody algebra
to a quiver with multiplicities $(\QQ,\ved)$.
Let $\bA=(a_{ij})_{i,j \in I}$ be the adjacency matrix of
the underlying graph of $\QQ$ and
set $\bD := (d_i\delta_{ij})_{i,j \in I}$.
Consider the generalized Cartan matrix
\[
\bC=(c_{ij})_{i,j \in I} := 2 \unit - \bA \bD.
\]
Note that it is symmetrizable as $\bD \bC = 2 \bD - \bD \bA \bD$
is symmetric.
Let
\[
\bigl( \g(\bC), \h, \{ \alpha_i \}_{i \in I}, \{ \alpha_i^\vee \}_{i \in I}
\bigr)
\]
be the corresponding Kac--Moody algebra
with its Cartan subalgebra, simple roots and simple coroots.
As usual we set
\[
Q := \sum_{i \in I} \Z \alpha_i, \qquad
Q_+ := \sum_{i \in I} \Z_{\geq 0} \alpha_i.
\]
The diagonal matrix $\bD$
induces a non-degenerate invariant symmetric bilinear form
$(\; , \,)$
on $\h^*$ satisfying
\[
( \alpha_i, \alpha_j) = d_i c_{ij}
= 2 d_i \delta_{ij} - d_i a_{ij} d_j,
\qquad  i,j \in I.
\]

From now on, we regard a dimension vector $\bv \in \Z_{\geq 0}^I$
of the quiver variety as an element of~$Q_+$~by
\[
\Z_{\geq 0}^I \xrightarrow{\simeq} Q_+,
\qquad \bv=(v_i)_{i \in I} \mapsto \sum_{i \in I} v_i \alpha_i.
\]
Let $\res \colon R^\ved \to \C^I$ be the map def\/ined by
\[
\res \colon \ \lambda=(\lambda_i(z)) \longmapsto
\Bigl( \res_{z=0} \lambda_i(z) \Bigr),
\]
and for $(\bv, \zeta) \in Q \times \C^I$,
let $\bv \cdot \zeta := \sum\limits_{i \in I} v_i \zeta_i$
be the scalar product.

\begin{Proposition}\label{prop:dim}\qquad
\begin{enumerate}\itemsep=0pt
\item[{\rm (i)}] The quiver variety $\Qst_{\QQ,\ved}(\lambda,\bv)$
is a holomorphic symplectic manifold of dimension
$2 -(\bv, \bv)$ if it is nonempty.

\item[{\rm (ii)}] If $\bv \cdot \res \lambda \neq 0$,
then $\Qset_{\QQ,\ved}(\lambda,\bv) = \varnothing$.

\item[{\rm (iii)}] If two quivers $\QQ_1$, $\QQ_2$ have the same underlying graph,
then the associated quiver varieties
$\Qst_{\QQ_1,\ved}(\lambda,\bv)$, $\Qst_{\QQ_2,\ved}(\lambda,\bv)$
are symplectomorphic to each other.
\end{enumerate}
\end{Proposition}

\begin{proof}
(i) Assume that $\Qst_{\QQ,\ved}(\lambda,\bv)$ is nonempty.
Since the action of $G_\ved(\bV)/\C^\times$
on the level set $\mu_\ved^{-1}(-\lambda\,\unit_\bV)^\st$ is free and proper,
the Marsden--Weinstein reduction theorem implies that
$\Qst_{\QQ,\ved}(\lambda,\bv)$ is a holomorphic symplectic manifold and
\begin{gather*}
\dim \Qst_{\QQ,\ved}(\lambda,\bv)  =
\dim \bM_{\QQ,\ved}(\bV) - 2 \dim G_\ved(\bV)/\C^\times
 = \sum_{i,j\in I} a_{ij} d_i v_i d_j v_j -2 \sum_{i \in I} d_i v_i^2 + 2 \\
\phantom{\dim \Qst_{\QQ,\ved}(\lambda,\bv)} = {}^t \bv \bD \bA \bD \bv -2 {}^t \bv \bD \bv +2 = 2-(\bv, \bv).
\end{gather*}

(ii) Assume $\Qset_{\QQ,\ved}(\lambda,\bv) \neq \varnothing$
and take a point $[B] \in \Qset_{\QQ,\ved}(\lambda,\bv)$.
Then we have
\[
\sum_{k=1}^{d_i} \sum_{\sumfrac{h \in H:}{\vin(h)=i}}
\epsilon(h) \tr_{R_{d_i}}
\bigl[ B_h B_{\ov{h}} N_i^{k-1} \bigr]\, z^{-k} = -\lambda_i(z)\,\unit_\bV
\]
for any $i \in I$.
Taking $\res\limits_{z=0} \circ \tr$ of both sides and sum over all $i$,
we obtain
\[
\sum_{h \in H} \epsilon(h) \tr ( B_h B_{\ov{h}} )
= -\sum_{i \in I} v_i \res_{z=0} \lambda_i(z) = -\bv \cdot \res \lambda.
\]
Here the left hand side is zero because
\[
\sum_{h \in H} \epsilon(h) \tr ( B_h B_{\ov{h}} ) =
\sum_{h \in H} \epsilon(\ov h) \tr ( B_{\ov{h}}B_h ) =
-\sum_{h \in H} \epsilon(h) \tr ( B_h B_{\ov{h}} ).
\]
Hence $\bv \cdot \res \lambda =0$.

(iii) By the assumption, we can identify the double quivers
$\QQ_1 + \ov{\QQ}_1$ and $\QQ_2 + \ov{\QQ}_2$.
Let $H$ be the set of arrows for them.
Then both the sets of arrows $\Omega_1$, $\Omega_2$ for $\QQ_1$, $\QQ_2$
are subsets of $H$.
Now the linear map
$\bM_{\QQ_1,\ved}(\bV) \to \bM_{\QQ_2,\ved}(\bV)=\bM_{\QQ_1,\ved}(\bV)$
def\/ined by
\[
B \mapsto B', \qquad B'_h:=
\begin{cases}
\hphantom{-}B_h &
\text{if}\hquad h \in \Omega_1 \cap \Omega_2
\hquad\text{or}\hquad h \in \ov{\Omega}_1 \cap \ov{\Omega}_2,
\\
-B_h & \text{otherwise},
\end{cases}
\]
induces a desired symplectomorphism
$\Qst_{\QQ_1,\ved}(\lambda,\bv) \xrightarrow{\simeq}
\Qst_{\QQ_2,\ved}(\lambda,\bv)$.
\end{proof}

Now f\/ix $i \in I$ and set
\[
\widehat{V}_i := \bigoplus_{\vin(h)=i}
V_{\vout(h)} \otimes_\C R_{d_{\vout(h)}}.
\]
Then using it we can decompose the vector space $\bM_{\QQ,\ved}(\bV)$ as
\begin{gather}\label{eq:decomposition}
\bM_{\QQ,\ved}(\bV) =\Hom(\widehat{V}_i,V_i \otimes_\C R_{d_i})
\oplus \Hom(V_i \otimes_\C R_{d_i}, \widehat{V}_i)
\oplus \bM^{(i)}_{\QQ,\ved}(\bV),
\end{gather}
where
\[
\bM^{(i)}_{\QQ,\ved}(\bV) := \bigoplus_{\vin(h), \vout(h) \neq i}
\Hom(V_{\vout(h)} \otimes_\C R_{d_{\vout(h)}},
V_{\vin(h)} \otimes_\C R_{d_{\vin(h)}}).
\]
According to this decomposition,
for a point $B \in \bM_{\QQ,\ved}(\bV)$ we put
\begin{gather*}
B_{i \slar}  := \bigl( \epsilon(h) B_h \bigr)_{\vin(h)=i}
\in \Hom(\widehat{V}_i,V_i \otimes_\C R_{d_i}), \\
B_{\slar i}  := \bigl( B_{\ov{h}} \bigr)_{\vin(h)=i}
\in \Hom(V_i \otimes_\C R_{d_i},\widehat{V}_i), \\
B_{\neq i}  := \bigl( B_h \bigr)_{\vin(h), \vout(h) \neq i}
\in \bM^{(i)}_{\QQ,\ved}(\bV).
\end{gather*}
We regard these as coordinates for $B$ and
write $B=(B_{i \slar},B_{\slar i},B_{\neq i})$.
Note that the symplectic form can be written as
\begin{gather}\label{eq:decomposeform}
\omega = \tr \rd B_{i \slar} \wedge \rd B_{\slar i}
+ \frac12 \sum_{\vout(h),\vin(h) \neq i}
\epsilon(h) \tr \rd B_h \wedge \rd B_{\ov{h}},
\end{gather}
and also the $i$-th component of the moment map can be written as
\[
\mu_{\ved,i}(B) = \prj_i (B_{i \slar} B_{\slar i})
= \sum_{k=1}^{d_i}
\tr_{R_{d_i}} \bigl[ B_{i \slar}  B_{\slar i} N_i^{k-1} \bigr]  z^{-k}.
\]

\begin{Lemma}\label{lem:stability}
Fix $i \in I$ and
suppose that $B$ satisfies at least one of the following two conditions:
\begin{enm}
\item $B$ is stable and $\bv \neq \alpha_i$;
\item the top coefficient
$\tr_{R_{d_i}}(B_{i \slar}  B_{\slar i} N_i^{d_i-1})$
of $\prj_i(B_{i \slar}B_{\slar i})$ is invertible.
\end{enm}
Then $(B_{i \slar},B_{\slar i})$ satisfies
\begin{gather}\label{eq:stability}
\Ker B_{\slar i} \cap \Ker N_i =0, \qquad
\range B_{i \slar} + \range N_i = V_i \otimes_\C R_{d_i}.
\end{gather}
\end{Lemma}

\begin{proof}
First, assume (i) and set
\begin{gather*}
\bS =\bigoplus_{j\in I} S_j, \qquad
S_j :=
\begin{cases}
\Ker B_{\slar i} \cap \Ker N_i & \text{if}\hquad j=i, \\
0 & \text{if}\hquad j \neq i,
\end{cases} \\
\bT =\bigoplus_{j\in I} T_j, \qquad
T_j :=
\begin{cases}
\range B_{i \slar} + \range N_i & \text{if}\hquad j=i, \\
V_j \otimes_\C R_{d_j} & \text{if}\hquad j \neq i.
\end{cases}
\end{gather*}
Then both $\bS$ and $\bT$ are $B$-invariant and
$N_j(S_j) \subset S_j,\, N_j(T_j) \subset T_j$ for all $j \in I$.
Since $B$ is stable,
we thus have $\bS =0$ or $\bS=\bV_\ved$,
and $\bT=0$ or $\bT=\bV_\ved$.
By the assumption $\bv \neq \alpha_i$ and the def\/initions of $\bS$ and $\bT$,
only the case $(\bS,\bT)=(0,\bV_\ved)$ occurs.
Hence $(B_{i \slar},B_{\slar i})$ satisf\/ies~\eqref{eq:stability}.

Next assume (ii).
Set $A(z)=\sum A_k z^{-k}:=\prj_i(B_{i \slar}B_{\slar i})$ and
\[
\widetilde{A} :=
\begin{pmatrix}
0 & 0 & \cdots & 0 \\
\vdots & \vdots & \ddots & \vdots \\
0 & 0 & \cdots & 0 \\
A_{d_i} & A_{d_i-1} & \cdots & A_1
\end{pmatrix}
\in \End_\C(V_i \otimes_\C \C^{d_i})=\End_\C(V_i \otimes_\C R_{d_i}).
\]
Then we have
\[
\tr_{R_{d_i}}(\widetilde{A}N^{k-1}) = A_k
= \tr_{R_{d_i}}(B_{i \slar}  B_{\slar i}N^{k-1}),
\qquad k=1,2,\dots ,d_i,
\]
i.e., $\widetilde{A}-B_{i \slar}B_{\slar i} \in \Ker \prj_i$.
Here, since the group $G_{d_i}(V_i)$ coincides with the centralizer of $N_i$
in $\GL_\C(V_i \otimes_\C R_{d_i})$, we have $\Ker \prj_i = \range \ad_{N_i}$.
Hence there is $C \in \End_\C(V_i \otimes_\C R_{d_i})$ such that
\[
B_{i \slar}  B_{\slar i} = \widetilde{A} + [N_i,C].
\]
Now note that both $\Ker N_i$ and $\Coker N_i$ are naturally isomorphic to $V_i$,
and the natural injection $\iota \colon \Ker N_i \to V_i \otimes_\C R_{d_i}$ and
projection $\pi \colon V_i \otimes_\C R_{d_i} \to \Coker N_i$
can be respectively identif\/ied with
the following block matrices:
\[
\begin{pmatrix} \unit_V \\ 0 \\ \vdots \\ 0
\end{pmatrix} \colon V_i \to V_i \otimes_\C R_{d_i}, \qquad
\begin{pmatrix}
0 & 0 & \cdots & 0 & \unit_V
\end{pmatrix} \colon V_i \otimes_\C R_{d_i} \to V_i.
\]
Thus we have
\begin{gather}\label{eq:momentrel}
\pi B_{i \slar}  B_{\slar i} \iota =
\pi (\widetilde{A} + [N_i,C]) \iota =
\pi \widetilde{A} \iota = A_{d_i}.
\end{gather}
By the assumption $A_{d_i}$ is invertible.
Hence $\pi B_{i \slar}$ is surjective and $B_{\slar i} \iota$ is injective.
\end{proof}

The following lemma is a consequence of results obtained in
\cite{0911.3863}:
\begin{Lemma}\label{lem:orbit1}
Suppose that the set
\begin{gather*}
Z_i :=
\{\, (B_{i \slar},B_{\slar i}) \in \Hom(\widehat{V}_i,V_i\otimes_\C R_{d_i})
\oplus \Hom(V_i\otimes_\C R_{d_i},\widehat{V}_i) \mid \\
\phantom{Z_i := \{\, }{} \prj_i(B_{i \slar}B_{\slar i})=-\lambda_i(z)\,\unit_{V_i},\quad
\text{\rm $(B_{i \slar}, B_{\slar i})$
satisfies \eqref{eq:stability}}
\,\}
\end{gather*}
is nonempty.
Then the quotient of it modulo the action of $G_{d_i}(V_i)$
is a smooth complex manifold having a
symplectic structure induced from
$\tr \rd B_{i \slar} \wedge \rd B_{\slar i}$,
and is symplectomorphic to
a $G_{d_i}(\widehat{V}_i)$-coadjoint orbit via the map given by
\[
\Phi_i \colon \ (B_{i \slar},B_{\slar i}) \longmapsto
- B_{\slar i}(z-N_i)^{-1}B_{i \slar} \in \g^*_{d_i}(\widehat{V}_i).
\]
\end{Lemma}

\begin{proof}
Take any point $(B_{i \slar},B_{\slar i})$ in the above set
and let $\bO$ be the $G_{d_i}(\widehat{V}_i)$-coadjoint orbit
through $\Phi_i(B_{i \slar},B_{\slar i})$.
By Proposition~4,~(a), Theorem~6 and Lemma~3 in \cite{0911.3863},
there exist
\begin{itm}
\item a f\/inite-dimensional $\C$-vector space $W$;
\item a nilpotent endomorphism $N \in \End(W)$ with $N^{d_i}=0$;
\item a coadjoint orbit $\bO_N \subset (\Lie G_N)^*$
of the centralizer $G_N \subset \GL(W)$ of $N$,
\end{itm}
such that the quotient modulo the natural $G_N$-action of the set
\[
\rset{(Y,X) \in \Hom(\widehat{V}_i,W)
\oplus \Hom(W,\widehat{V}_i)}{%
\begin{aligned}
&\prj_N(YX) \in \bO_N,\\
&\Ker X \cap \Ker N =0, \quad
\range Y + \range N =\widehat{V}_i
\end{aligned}
},
\]
where $\prj_N$ is the transpose of the inclusion
$\Lie G_N \hookrightarrow \gl(W)$,
is a smooth manifold having a symplectic
structure induced from $\tr \rd X \wedge \rd Y$,
and is symplectomorphic to $\bO$ via the map
$(Y,X) \mapsto X(z\,\unit_{W} -N)^{-1}Y$.
Note that if $W=V_i\otimes_\C R_{d_i}$, $N=N_i$ and
$\bO_N$ is a single element $\lambda_i(z)\,\unit_{V_i}$,
then we obtain the result by the coordinate change
$(Y,X)=(B_{i \slar},-B_{\slar i})$.
Indeed, this is the case thanks to Proposition~4,~(c) and
Theorem~6 (the uniqueness assertion) in \cite{0911.3863}.
\end{proof}

Note that in the above lemma, the assumption $Z_i \neq \varnothing$
implies $\dim V_i \leq \dim \widehat{V}_i$;
indeed, if $(B_{i \slar},B_{\slar i}) \in Z_i$,
then $B_{\slar i}|_{\Ker N_i}$ is injective by condition \eqref{eq:stability},
and hence
\[
\dim V_i = \dim \Ker N_i = \rank \left( B_{\slar i}|_{\Ker N_i} \right)
\leq \dim \widehat{V}_i.
\]
The following lemma tells us that if the top coef\/f\/icient
of $\lambda_i(z)$ is nonzero,
then the converse is true and
the corresponding coadjoint orbit can be explicitly described:
\begin{Lemma}\label{lem:orbit2}
Suppose $\dim V_i \leq \dim \widehat{V}_i$ and that the top coefficient
$\lambda_{i,d_i}$ of $\lambda_i(z)$ is nonzero.
Then the set $Z_i$ in Lemma~{\rm \ref{lem:orbit1}} is nonempty
and the coadjoint orbit
contains an element of the form
\[
\Lambda_i(z)= \begin{pmatrix}
\lambda_i(z)\,\unit_{V_i} & 0 \\
0 & 0\,\unit_{V'_i} \end{pmatrix},
\]
where $V_i$ is regarded as a subspace of $\widehat{V}_i$ and
$V'_i \subset \widehat{V}_i$ is a complement of it.
\end{Lemma}

\begin{proof}
Suppose $\dim V_i \leq \dim \widehat{V}_i$
and that the top coef\/f\/icient of $\lambda_i(z)$ is nonzero. We set
\begin{gather*}
B_{i \slar}  := \begin{pmatrix}
0 & 0\, \\ \vdots & \vdots\, \\ 0 & 0\, \\ \unit_{V_i} & 0\,
\end{pmatrix}
\colon \widehat{V}_i = V_i \oplus V'_i \to V_i \otimes_\C R_{d_i}, \\
B_{\slar i}  := - \begin{pmatrix}
\lambda_{i,d_i}\,\unit_{V_i} & \cdots & \lambda_{i,1}\,\unit_{V_i} \\
0 & \cdots & 0
\end{pmatrix}
\colon V_i \otimes_\C R_{d_i} \to \widehat{V}_i,
\end{gather*}
where $\lambda_{i,k}$ denotes the coef\/f\/icient in $\lambda_i(z)$ of $z^{-k}$.
Then we have
\[
\tr_{R_{d_i}}(B_{i \slar}  B_{\slar i} N_i^{k-1})
= -\lambda_{i,k}\,\unit_{V_i}, \qquad k=1,2, \dots ,d_i,
\]
i.e., $\prj_i(B_{i \slar}  B_{\slar i})=-\lambda_i(z)\,\unit_{V_i}$.
The assumption $\lambda_{i,d_i} \neq 0$ and
\lemref{lem:stability} imply that
$(B_{i \slar},B_{\slar i})$ satisf\/ies \eqref{eq:stability}.
Hence $(B_{i \slar},B_{\slar i}) \in Z_i$.
Moreover we have
\begin{gather*}
\Phi_i(B_{i \slar},B_{\slar i})=
-B_{\slar i}(z-N_i)^{-1}B_{i \slar}  =
-\sum_{k=1}B_{\slar i}N_i^{k-1}B_{i \slar}\,z^{-k} \\
\hphantom{\Phi_i(B_{i \slar},B_{\slar i})=
-B_{\slar i}(z-N_i)^{-1}B_{i \slar}}{}
=\sum_{k=1}
\begin{pmatrix}
\lambda_{i,k}\,\unit_{V_i} & 0 \\
0 & 0\,\unit_{V'_i} \end{pmatrix} z^{-k} = \Lambda_i(z).\tag*{\qed}
\end{gather*}
\renewcommand{\qed}{}
\end{proof}

\section{Ref\/lection functor}\label{sec:reflection}

In this section we construct {\em reflection functors}
for quiver varieties with multiplicities.

\subsection{Main theorem}\label{subsec:main}

Recall that the Weyl group $W(\bC)$ of the Kac--Moody algebra $\g(\bC)$
is the subgroup of $\GL(\h^*)$ generated by the simple ref\/lections
\[
s_i(\beta) := \beta -
\langle \beta, \alpha^\vee_i \rangle \alpha_i
= \beta - \frac{2(\beta,\alpha_i)}{(\alpha_i,\alpha_i)}\alpha_i,
\qquad i \in I,\quad \beta \in \h^*.
\]
The fundamental relations for the generators
$s_i,\, i \in I$ are
\begin{gather}\label{eq:coxeter}
s_i^2 =\unit, \qquad
(s_is_j)^{m_{ij}} =\unit, \qquad i,j \in I,\quad i \neq j,
\end{gather}
where the numbers $m_{ij}$ are determined from $c_{ij}c_{ji}$
as the table below
(we use the convention $r^\infty=\unit$ for any $r$)
$$
{\renewcommand{\arraystretch}{1.2}
\begin{array}{|c|c|c|c|c|c|}
\hline
c_{ij}c_{ji} & 0 & 1 & 2 & 3 & \geq 4 \\ \hline
m_{ij} & 2 & 3 & 4 & 6 & \infty \\ \hline
\end{array}}
$$

We will def\/ine a $W(\bC)$-action on the parameter space
$R^\ved \times Q$ for the quiver variety.
The action on the second component $Q$ is given by just
the restriction of the standard action on $\h^*$, namely,
\[
s_i \colon \ \bv=\sum_{i \in I} v_i \alpha_i \longmapsto
\bv - \langle \bv, \alpha^\vee_i \rangle \alpha_i
=\bv - \sum_{j \in I} c_{ij}v_j \alpha_i.
\]
The action on the f\/irst component $R^\ved$ is unusual.
We def\/ine $r_i \in \GL(R^\ved)$ by
\[
r_i(\lambda) = \lambda' \equiv (\lambda'_j(z)), \qquad
\lambda'_j(z) :=
\begin{cases}
-\lambda_i(z) & \text{if}\hquad j =i, \\
\lambda_j(z) - z^{-1}c_{ij} \res\limits_{z=0} \lambda_i(z) & \text{if}\hquad j \neq i.
\end{cases}
\]

\begin{Lemma}\label{lem:rel}
The above $r_i$, $i \in I$ satisfy relations \eqref{eq:coxeter}.
\end{Lemma}

\begin{proof}
The relations $r_i^2 = \unit$, $i \in I$ are obvious.
To check the relation $(r_ir_j)^{m_{ij}}=\unit$ for $i\neq j$,
f\/irst note that the transpose of $s_i \colon Q \to Q$
relative to the scalar product is given by
\[
{}^t s_i \colon \C^I \to \C^I, \qquad
{}^t s_i (\zeta) = \zeta - \zeta_i\sum_{j \in I} c_{ij} \alpha_j.
\]
Now let $\lambda \in R^\ved$. We decompose it as
\[
\lambda = \lambda^0 + \res (\lambda)\,z^{-1},
\qquad \res \lambda^0 =0.
\]
Then we easily see that
\[
r_i(\res(\lambda)\,z^{-1})={}^t s_i (\res(\lambda))\,z^{-1},
\]
and hence that
\[
(r_ir_j)^{m_{ij}} (\lambda) = (r_ir_j)^{m_{ij}} (\lambda^0)
+ \res (\lambda)\,z^{-1}.
\]
Therefore we may assume that $\res \lambda =0$.
Set $\lambda' \equiv (\lambda'_k(z)) := (r_ir_j)^{m_{ij}}(\lambda)$.
Then we have
\[
\lambda'_k(z) =
\begin{cases}
(-1)^{m_{ij}} \lambda_k(z) & \text{if}\hquad k=i,j, \\
\lambda_k(z) & \text{if}\hquad k \neq i,j.
\end{cases}
\]
If $m_{ij}$ is odd, by the def\/inition we have
$c_{ij}c_{ji}=1$.
In particular, $i \neq j$ and
\[
a_{ij} d_j a_{ji} d_i = c_{ij} c_{ji} =1.
\]
This implies $d_i=d_j=1$ and hence that
$\lambda_i(z) = \lambda_j(z)=0$.
\end{proof}

The main result of this section is as follows:
\begin{Theorem}\label{thm:ref}
Let $\lambda=(\lambda_i(z)) \in R^\ved$ and suppose that
the top coefficient $\lambda_{i,d_i}$ of $\lambda_i(z)$
for fixed $i \in I$ is nonzero.
Then there exists a bijection
\[
\frF_i \colon \ \Qset_{\QQ,\ved}(\lambda,\bv) \xrightarrow{\simeq}
\Qset_{\QQ,\ved}(r_i(\lambda),s_i(\bv))
\]
such that $\frF_i^2 = \unit$
and the restriction gives a symplectomorphism
\[
\frF_i \colon \ \Qst_{\QQ,\ved}(\lambda,\bv) \xrightarrow{\simeq}
\Qst_{\QQ,\ved}(r_i(\lambda),s_i(\bv)).
\]
\end{Theorem}

We call the above map $\frF_i$
the {\em $i$-th reflection functor}.

\subsection{Proof of the main theorem}\label{subsec:proof}

Fix $i \in I$ and
suppose that the top coef\/f\/icient
$\lambda_{i,d_i}$ of $\lambda_i(z)$ is nonzero.
Recall the decomposition~\eqref{eq:decomposition} of $\bM_{\QQ,\ved}(\bV)$:
\[
\bM_{\QQ,\ved}(\bV) =\Hom(\widehat{V}_i,V_i \otimes_\C R_{d_i})
\oplus \Hom(V_i \otimes_\C R_{d_i}, \widehat{V}_i)
\oplus \bM^{(i)}_{\QQ,\ved}(\bV),
\]
and the set $Z_i$ given in \lemref{lem:orbit1}.
\lemref{lem:stability} and the assumption $\lambda_{i,d_i} \neq 0$
imply that any
$B=(B_{i \slar},B_{\slar i},B_{\neq i})
\in \mu_{\ved,i}^{-1}(-\lambda_i(z)\,\unit_{V_i})$
satisf\/ies condition \eqref{eq:stability}.
Thus we have
\[
\mu_{\ved,i}^{-1}(-\lambda_i(z)\,\unit_{V_i})
=Z_i \times \bM^{(i)}_{\QQ,\ved}(\bV).
\]
By \lemref{lem:orbit2},
it is nonempty if and only if
\[
v_i \leq \dim \widehat{V}_i
= \sum_j a_{ij} d_j v_j = 2v_i -\sum_j c_{ij}v_j,
\]
i.e., the $i$-th component of $s_i(\bv)$ is non-negative.
We assume this condition,
because otherwise both
$\Qset_{\QQ,\ved}(\lambda,\bv)$ and
$\Qset_{\QQ,\ved}(r_i(\lambda),s_i(\bv))$ are empty
(since $s_i(\bv) \notin \Z_{\geq 0}^I$).
Fix a $\C$-vector space~$V'_i$ of dimension
$\dim \widehat{V}_i - \dim V_i$ and an identif\/ication
$\widehat{V}_i = V_i \oplus V'_i$.
As the group $G_{d_i}(V_i)$ acts trivially on $\bM^{(i)}_{\QQ,\ved}(\bV)$,
Lemmas~\ref{lem:orbit1} and \ref{lem:orbit2} imply that
\[
\mu_{\ved,i}^{-1}(-\lambda_i(z)\,\unit_{V_i})/G_{d_i}(V_i) =
Z_i/G_{d_i}(V_i) \times \bM^{(i)}_{\QQ,\ved}(\bV)
\simeq \bO \times \bM^{(i)}_{\QQ,\ved}(\bV),
\]
where $\bO$ is the $G_{d_i}(\widehat{V}_i)$-coadjoint orbit
through
\[
\Lambda(z)= \begin{pmatrix}
\lambda_i(z)\,\unit_{V_i} & 0 \\
0 & 0\,\unit_{V'_i} \end{pmatrix}.
\]
Now let us def\/ine an $I$-graded $\C$-vector space $\bV'$
with $\dim \bV'=s_i(\bv)$ by
\[
\bV'=\bigoplus_{j \in I} V'_j, \qquad
V'_j := \begin{cases}
V'_i & \text{if}\hquad j=i, \\
V_j & \text{if}\hquad j\neq i,
\end{cases}
\]
and consider the associated symplectic vector space $\bM_{\QQ,\ved}(\bV')$.
Note that $\widehat{V}'_i =\widehat{V}_i$.
Thus by interchanging the roles of $\bV$ and $\bV'$,
$\lambda_i$ and $-\lambda_i$ in Lemmas~\ref{lem:orbit1} and \ref{lem:orbit2},
we obtain an isomorphism
\[
\mu_{\ved,i}^{-1}(\lambda_i(z)\,\unit_{V'_i})/G_{d_i}(V'_i) \simeq
\bO' \times \bM^{(i)}_{\QQ,\ved}(\bV')=\bO' \times \bM^{(i)}_{\QQ,\ved}(\bV),
\]
where $\bO'$ is the $G_{d_i}(\widehat{V}_i)$-coadjoint orbit through
\[
\begin{pmatrix} 0\,\unit_{V_i} & 0 \\
0 & -\lambda_i(z)\,\unit_{V'_i} \end{pmatrix}
=\Lambda(z)-\lambda_i(z)\,\unit_{\widehat{V}_i},
\]
i.e., $\bO' = \bO - \lambda_i(z)\,\unit_{\widehat{V}_i}$.
Hence the scalar shift
$\bO \xrightarrow{\simeq} \bO - \lambda_i(z)\,\unit_{\widehat{V}_i}$
induces an isomorphism
\[
\widetilde{\frF}_i \colon \
\mu_{\ved,i}^{-1}(-\lambda_i(z)\,\unit_{V_i})/G_{d_i}(V_i) \xrightarrow{\simeq}
\mu_{\ved,i}^{-1}(\lambda_i(z)\,\unit_{V'_i})/G_{d_i}(V'_i),
\]
which is characterized as follows: if
\[
\widetilde{\frF}_i[B]=[B'], \qquad
B=(B_{i \slar},B_{\slar i},B_{\neq i}),\qquad
B' =(B'_{i \slar},B'_{\slar i},B_{\neq i}'),
\]
one has
\begin{gather}
B_{\neq i}  = B_{\neq i}', \label{eq:refrelother}\\
- B'_{\slar i}(z-N'_i)^{-1}B'_{i \slar}
 =- B_{\slar i}(z-N_i)^{-1}B_{i \slar} - \lambda_i(z)\,\unit_{\widehat{V}_i},
\label{eq:shift}
\end{gather}
where $N'_i := \unit_{V'_i} \otimes_\C J_{d_i}
\in \End_\C(V'_i \otimes_\C R_{d_i})$.
Note that
\begin{gather}\label{eq:refstability}
\Ker B'_{\slar i} \cap \Ker N'_i =0, \qquad
\range B'_{i \slar} + \range N'_i = V'_i \otimes_\C R_{d_i}
\end{gather}
by \lemref{lem:stability}.

\begin{Lemma}\label{lem:refmoment}
If $\mu_\ved(B)=-\lambda\,\unit_\bV$,
then $\mu_{\ved}(B')=-r_i(\lambda)\,\unit_{\bV'}$.
\end{Lemma}

\begin{proof}
Let $\lambda'=(\lambda'_j(z)) := r_i(\lambda)$.
The identity $\mu_{\ved,i}(B')=\lambda_i(z)\,\unit_{V'_i}$
is clear from the construction.
We check $\mu_{\ved,j}(B')=-\lambda'_j(z)\,\unit_{V'_j}$ for $j \neq i$.
Taking the residue of both sides of~\eqref{eq:shift}, we have
\[
 B'_{\slar i} B'_{i \slar} =  B_{\slar i} B_{i \slar}
+ \lambda_{i,1}\,\unit_{\widehat{V}_i},
\]
which implies that
\[
\epsilon(h) B'_{\ov{h}} B'_h = \epsilon(h) B_{\ov{h}} B_h
+ \lambda_{i,1}\, \unit_{V_{\vout(h)}} \qquad
\text{if} \quad \vin(h)=i.
\]
On the other hand, \eqref{eq:refrelother} means that
$B'_h = B_h$ whenever $\vin(h), \vout(h) \neq i$.
Thus for $j \neq i$, we obtain
\begin{gather}
\sum_{\vin(h)=j} \epsilon(h) B'_h B'_{\ov{h}}
 = \sum_{h \colon i \to j} \epsilon(h) B'_h B'_{\ov{h}} +
\sum_{\vin(h)=j, \vout(h) \neq i} \epsilon(h) B'_h B'_{\ov{h}} \notag\\
\hphantom{\sum_{\vin(h)=j} \epsilon(h) B'_h B'_{\ov{h}}}{}
= \sum_{h \colon i \to j}
\bigl( \epsilon(h) B_h B_{\ov{h}} - \lambda_{i,1}\,\unit_{V_j} \bigr)
+\sum_{\vin(h)=j, \vout(h) \neq i} \epsilon(h) B_h B_{\ov{h}} \notag\\
\hphantom{\sum_{\vin(h)=j} \epsilon(h) B'_h B'_{\ov{h}}}{}
=\sum_{\vin(h)=j} \epsilon(h) B_h B_{\ov{h}}
- a_{ij}\lambda_{i,1}\,\unit_{V_j}.\label{eq:refmoment2}
\end{gather}
Note that
\[
\prj_j(\unit_{V_j}) = \sum_{k=1}^{d_j} \tr_{R_{d_j}} (N_j^{k-1}) z^{-k}
= d_j \unit_{V_j} z^{-1}.
\]
Therefore the image under $\prj_j$
of both sides of~\eqref{eq:refmoment2} gives
\[
\mu_{\ved,j}(B') = \mu_{\ved,j}(B)
- a_{ij}\lambda_{i,1} \prj_i(\unit_{V_j})
=\mu_{\ved,j}(B)
+ c_{ij}\lambda_{i,1}\,\unit_{V_j}z^{-1}.
\]
The result follows.
\end{proof}

\begin{Lemma}\label{lem:refirred}
If $B$ is stable, then so is~$B'$.
\end{Lemma}

\begin{proof}
Suppose that there exists a $B'$-invariant subspace
$\bS'=\bigoplus_j S'_j {\subset} \bV'_\ved$ such that
$N'_j(S'_j) {\subset} S'_j$.
We def\/ine an $I$-graded subspace $\bS=\bigoplus_j S_j$ of
$\bV_\ved$ by
\[
S_j := \begin{cases}
\displaystyle \sum\limits_{k=1}^{d_i} N_i^{k-1} B_{i \slar}\bigl( \widehat{S}'_i \bigr)
& \text{if}\hquad j=i,\\
S'_j & \text{if}\hquad j \neq i,
\end{cases}
\]
where $\widehat{S}'_i := \bigoplus_{\vin(h)=i} S'_{\vout(h)} = \widehat{S}_i$.
Then $B_{i \slar}\bigl( \widehat{S}_i \bigr) \subset S_i$ and
\begin{gather*}
B_{\slar i} (S_i)  = \sum_{k=1}^{d_i}
B_{\slar i} N_i^{k-1} B_{i \slar}\bigl( \widehat{S}'_i \bigr)
 =\sum_{k=1}^{d_i}
\bigl( B'_{\slar i} (N'_i)^{k-1} B'_{i \slar} -\lambda_{i,k} \bigr)
\bigl( \widehat{S}'_i \bigr)
 \subset \widehat{S}'_i = \widehat{S}_i.
\end{gather*}
Hence $\bS$ is $B$-invariant.
Clearly $N_j(S_j) \subset S_j$ for all $j \in I$.
Therefore the stability condition for~$B$ implies that
$\bS=0$ or $\bS=\bV_\ved$.
First, assume $\bS=0$.
Then $S'_j=S_j=0$ for $j \neq i$, and hence
$B'_{\slar i}(S'_i) \subset \widehat{S}'_i=0$, i.e.,
$S'_i \subset \Ker B'_{\slar i}$.
If $S'_i$ is nonzero, then the kernel of the restriction
$N'_i |_{S'_i}$ is nonzero because it is nilpotent.
However it implies $\Ker B'_{\slar i} \cap \Ker N'_i \neq 0$,
which contradicts to \eqref{eq:refstability}.
Hence $S'_i=0$.
Next assume $\bS=\bV_\ved$.
Then $S'_j=S_j=V_j \otimes_\C R_{d_j}$ for $j \neq i$, and hence
$S'_i \supset B'_{i \slar}\bigl( \widehat{S}'_i \bigr) = \range B'_{i \slar}$.
If $V'_i/S'_i$ is nonzero,
then the endomorphism of $V'_i/S'_i$ induced from $N'_i$
has a nonzero cokernel because it is nilpotent.
However it implies
$\range B'_{i \slar} + \range N'_i \neq V'_i \otimes_\C R_{d_i}$,
which contradicts to \eqref{eq:refstability}.
Hence $S'_i=V'_i \otimes_\C R_{d_i}$.
\end{proof}

\begin{proof}[Proof of \thmref{thm:ref}]
As the map $\widetilde{\frF}_i$ is clearly
$\prod\limits_{j \neq i}G_{d_j}(V_j)$-equivariant,
\lemref{lem:refmoment} implies
that it induces a bijection
\[
\frF_i \colon \  \Qset_{\QQ,\ved}(\lambda,\bv) \to \Qset_{\QQ,\ved}(r_i(\lambda),s_i(\bv)),
\qquad [B] \mapsto [B'],
\]
which preserves the stability by \lemref{lem:refirred}.
We easily obtain the relation $\frF_i^2 = \unit$
by noting that $\frF_i$ is induced from the scalar shift
$\bO \to \bO'=\bO -\lambda_i(z)\,\unit_{\widehat{V}_i}$
and the $i$-th component of~$r_i(\lambda)$ is~$-\lambda_i$.
Consider the restriction
\[
\frF_i \colon \ \Qst_{\QQ,\ved}(\lambda,\bv) \to \Qst_{\QQ,\ved}(r_i(\lambda),s_i(\bv)),
\qquad [B] \mapsto [B'].
\]
By \lemref{lem:orbit1} and \eqref{eq:shift},
we have
\[
\tr \rd B_{i \slar} \wedge \rd B_{\slar i} =
\tr \rd B'_{i \slar} \wedge \rd B'_{\slar i},
\]
because the scalar shift
$\bO \to \bO'$
is a symplectomorphism.
Substituting it and \eqref{eq:refrelother}
into \eqref{eq:decomposeform},
we see that the above map $\frF_i$
is a symplectomorphism.
\end{proof}

\begin{Remark}\label{rem:original-ref}
It is clear from \eqref{eq:refrelother}, \eqref{eq:shift}
and \eqref{eq:refstability} that if $d_j=1$ for all $j \in I$,
then $\frF_i$ coincides with the original $i$-th ref\/lection functor
for quiver varieties
(see conditions (a), (b1) and (c) in \cite[Section 3]{MR1623674}).
\end{Remark}

\begin{Remark}\label{rem:ref-relation}
It is known (see e.g.\ \cite{MR1775358}) that if $d_i=1$ for all $i \in I$,
then the ref\/lection functors~$\frF_i$ satisfy relations~\eqref{eq:coxeter}.
We expect that this fact is true for any~$(\QQ,\ved)$.
\end{Remark}

\subsection{Application}\label{subsec:app}

In this subsection we introduce a basic application of ref\/lection functors.

\begin{Lemma}
Let $(\lambda,\bv) \in R^\ved \times Q_+$, $i \in I$.
Suppose that the top coefficient of $\lambda_i(z)$ is zero
and $\bv \neq \alpha_i$.
Then $\Qst_{\QQ,\ved}(\lambda,\bv) \neq \varnothing$ implies $(\bv,\alpha_i) \leq 0$.
\end{Lemma}

\begin{proof}
Take any point $[B] \in \Qst_{\QQ,\ved}(\lambda,\bv)$.
Let $\iota \colon \Ker N_i \to V_i \otimes_\C R_{d_i}$ be the inclusion
and $\pi \colon V_i \otimes_\C R_{d_i} \to \Coker N_i$ be the projection.
Then \lemref{lem:stability} together with the assumption
$\bv \neq \alpha_i$ implies that
$B_{\slar i} \iota$ is injective and $\pi B_{i \slar}$ is surjective.
On the other hand, \eqref{eq:momentrel}
and the assumption for $\lambda_i(z)$ imply that
\[
\xymatrix{
V_i \simeq \Ker N_i \ar[r]^-{B_{\slar i} \iota} &
\widehat{V}_i \ar[r]^-{\pi B_{i \slar}} &
\Coker N_i \simeq V_i
}
\]
is a complex. Thus we have
\[
0 \leq \dim \widehat{V}_i - 2 \dim V_i = \sum_j a_{ij} d_j v_j -2v_i,
\]
which is equivalent to $(\bv,\alpha_i) \leq 0$.
\end{proof}

Now applying Crawley-Boevey's argument
in \cite[Lemma~7.3]{MR1834739} to our quiver varieties with multiplicities,
we obtain the following:
\begin{Proposition}\label{prop:root}
If $\Qst_{\QQ,\ved}(\lambda,\bv) \neq \varnothing$,
then $\bv$ is a positive root of $\g(\bC)$.
\end{Proposition}

\begin{proof}
Assume $\Qst_{\QQ,\ved}(\lambda,\bv) \neq \varnothing$ and that
$\bv$ is not a real root.
We show that $\bv$ is an imaginary root using \cite[Theorem~5.4]{MR1104219};
namely, show that there exists $w \in W(\bC)$ such that
$w(\bv)$ has a~connected support
and $(w(\bv),\alpha_i) \leq 0$ for any $i \in I$.

Assume that there is $i \in I$ such that $(\bv,\alpha_i) >0$.
The above lemma implies that
the top coef\/f\/i\-cient of $\lambda_i(z)$ is nonzero,
which together with \thmref{thm:ref}
implies that \mbox{$\Qst_{\QQ,\ved}(r_i(\lambda),s_i(\bv)) \neq \varnothing$}.
In particular we have $s_i(\bv) \in Q_+$,
and further $\bv -s_i(\bv) \in \Z_{>0} \alpha_i$
by the assumption \mbox{$(\bv,\alpha_i)>0$}.
We then replace $(\lambda,\bv)$ with
$(r_i(\lambda),s_i(\bv))$, and repeat this argument.
As the components of $\bv$ decrease,
it eventually stops after f\/inite number of steps,
and we f\/inally obtain a pair $(\lambda,\bv) \in R^\ved \times Q_+$
such that $(\bv,\alpha_i) \leq 0$ for all $i \in I$.
Additionally, the property $\Qst_{\QQ,\ved}(\lambda,\bv) \neq \varnothing$
clearly implies that the support of $\bv$ is connected.
The result follows.
\end{proof}

\section{Normalization}\label{sec:normalization}

In this section we give an application of Boalch's `shifting trick'
to quiver varieties with multiplicities.

\subsection{Shifting trick}\label{subsec:shifting}

\begin{Definition}\label{dfn:pole}
Let $(\QQ,\ved)$ be a quiver with multiplicities.
A vertex $i \in I$ is called a {\em pole vertex}
if there exists a unique vertex $j \in I$ such that
\[
d_j =1, \qquad a_{ik}=a_{ki}=\delta_{jk} \quad \text{for any}\hquad k \in I.
\]
The vertex $j$ is called the {\em base vertex} for the pole $i$.
If furthermore $d_i>1$, the pole $i \in I$ is said to be {\em irregular}.
\end{Definition}

Let $i \in I$ be a pole vertex with the base $j \in I$.
Then $\widehat{V}_i = V_j \otimes_\C R_{d_j} =V_j$.
In what follows we assume that the top coef\/f\/icient of $\lambda_i(z)$
is nonzero.
As the set
$\Qset_{\QQ,\ved}(\lambda,\bv)$ is empty unless
$\dim V_i \leq \dim \widehat{V}_i =\dim V_j$,
we also assume $V_i \subset V_j$ and f\/ix an identif\/ication
$V_j \simeq V_i \oplus V_j/V_i$.
Recall the isomorphism given in the previous section:
\[
\mu_{\ved,i}^{-1}(-\lambda_i(z)\,\unit_{V_i})/G_{d_i}(V_i) \simeq
\bO \times \bM^{(i)}_{\QQ,\ved}(\bV),
\]
where $\bO$ is the $G_{d_i}(V_j)$-coadjoint orbit
through the element of the form
\[
\Lambda(z)= \begin{pmatrix}
\lambda_i(z)\,\unit_{V_i} & 0 \\
0 & 0\,\unit_{V_j/V_i} \end{pmatrix}.
\]
Let us decompose $\Lambda(z)$ as
\[
\Lambda(z)=\Lambda^0(z) + z^{-1}\res_{z=0} \Lambda(z)
\]
according to the decomposition
\[
\g^*_{d_i}(V_j) = \frb^*_{d_i}(V_j) \oplus z^{-1} \gl(V_j),
\]
where
\[
\frb^*_{d_i}(V_j) :=
\Ker \Bigl[ \,\res_{z=0} \colon \g^*_{d_i}(V_j) \to \gl(V_j) \Bigr]
\simeq z^{-d_i}\gl(V_j)[[z]]/z^{-1}\gl(V_j)[[z]].
\]
The above is naturally
dual to the Lie algebra $\frb_{d_i}(V_j)$ of
the unipotent subgroup
\[
B_{d_i}(V_j) := \big\{\, g(z) \in G_{d_i}(V_j) \mid g(0) = \unit_{V_j}\,\big\}.
\]
The coadjoint action of $g(z) \in B_{d_i}(V_j)$ is given by
\[
(g \cdot \eta)(z) = g(z)\eta(z)g(z)^{-1} \mod z^{-1}\gl(V_j)[[z]],
\qquad \eta(z) = \sum_{k=2}^{d_j} \eta_k z^{-k} \in \frb^*_{d_i}(V_j).
\]
Now consider the $B_{d_i}(V_j)$-coadjoint orbit $\check{\bO}$
through $\Lambda^0(z)$. Let
\[
K:=\GL(V_i) \times \GL(V_j/V_i) \subset \GL(V_j)
\]
be the Levi subgroup associated to the decomposition
$V_j = V_i \oplus V_j/V_i$.
The results in this section is based on the following two facts:

\begin{Lemma}\label{lem:boalch1}
The orbit $\check{\bO}$
is invariant under the conjugation action by $K$,
and there exists a~$K$-equivariant
algebraic symplectomorphism
\[
\check{\bO} \simeq
\Hom(V_j/V_i,V_i)^{\oplus (d_i-2)} \oplus \Hom(V_i,V_j/V_i)^{\oplus (d_i-2)}
\]
sending $\Lambda^0(z) \in \check{\bO}$ to the origin.
\end{Lemma}

\begin{Lemma}\label{lem:boalch2}
Let $M$ be a holomorphic symplectic manifold
with a Hamiltonian action of $\GL(V_j)$ and
a moment map $\mu_M \colon M \to \gl(V_j)$.
Then for any $\zeta \in \C$, the map
\[
\check{\bO} \times M \to \g^*_{d_i}(V_j) \times M, \qquad
(B(z), x) \mapsto \big(B(z)-z^{-1}\mu_M(x)-z^{-1}\zeta\,\unit_{V_j},x\big)
\]
induces a bijection between
\begin{enm}
\item the $($set-theoretical$)$ symplectic quotient of $\check{\bO} \times M$ by
the diagonal $K$-action at the level
$-\res\limits_{z=0}\Lambda(z)-\zeta\,\unit_{V_j}$; and
\item that of $\bO \times M$ by the diagonal $\GL(V_j)$-action
at the level $-\zeta\,\unit_{V_j}$.
\end{enm}
Furthermore, under this bijection
a point in the space {\rm (i)} represents a free $K$-orbit
if and only if the corresponding point in the space {\rm (ii)}
represents a free $\GL(V_j)$-orbit,
at which the two symplectic forms are intertwined.
\end{Lemma}

\lemref{lem:boalch2} is what we call `Boalch's shifting trick'.
We directly check the above two facts in Appendix~\ref{app:misc}.

\begin{Remark}\label{rem:HTL}
Let $\Lambda^1, \Lambda^2, \dots , \Lambda^k \in \End(V)$ be
mutually commuting endomorphisms of
a $\C$-vector space $V$,
and suppose that $\Lambda^2, \dots ,\Lambda^k$ are semisimple.
To such endomorphisms we associate
\[
\Lambda(z) := \sum_{j=1}^k \Lambda^j z^{-j} \in \g^*_k(V),
\]
which is called a {\em normal form}.
Let $\Sigma \subset \g^*_k(\C)$ be the subset consisting of all
residue-free elements $\lambda(z)=\sum\limits_{j=2}^k \lambda^j z^{-j}$
with $(\lambda^2, \dots ,\lambda^k)$
being a simultaneous eigenvalue of $(\Lambda^2,\dots ,\Lambda^k)$,
and let $V=\bigoplus_{\lambda \in \Sigma} V_\lambda$
be the eigenspace decomposition.
Then we can express $\Lambda(z)$ as
\[
\Lambda(z)=\bigoplus_{\lambda \in \Sigma}
\left( \lambda(z)\,\unit_{V_\lambda} + \frac{\Gamma_\lambda}{z} \right),
\qquad \Gamma_\lambda=\Lambda^1 |_{V_\lambda} \in \End(V_\lambda).
\]
It is known that any $A(z) \in \g^*_k(V)$ whose leading term
is regular semisimple is equivalent to
some normal form under the coadjoint action.

Note that $\Lambda(z)$ treated in
Lemmas~\ref{lem:boalch1} and \ref{lem:boalch2}
is a normal form.
A generalization of \lemref{lem:boalch1}
for an arbitrary normal form
has been announced in \cite[Appendix~C]{0806.1050}.
\lemref{lem:boalch2} is known
in the case where $\Lambda(z)$ is a normal form
whose leading term is regular semisimple~\cite{MR1864833};
however, as mentioned in \cite{0806.1050},
the arguments in \cite[Section~2]{MR1864833}
needed to prove this fact can be generalized
to the case where $\Lambda(z)$ is an arbitrary normal form.
\end{Remark}

We apply \lemref{lem:boalch2} to the case where
$M=\bM^{(i)}_{\QQ,\ved}(\bV)$, $\zeta=\res\limits_{z=0}\lambda_j(z)$.
In this case, the symp\-lectic quotient of the space (ii) by
the action of $\prod\limits_{k \neq i,j} G_{d_k}(V_k)$
turns out to be
$\mu_\ved^{-1}(-\lambda\,\unit_\bV)/G_\ved(\bV)$ $=\Qset_{\QQ,\ved}(\lambda,\bv)$.
On the other hand, by \lemref{lem:boalch1},
the symplectic quotient of the space (i) by
the action of $\prod\limits_{k \neq i,j} G_{d_k}(V_k)$
coincides with the symplectic quotient of
\begin{gather}\label{eq:normalization1}
\Hom(V_j/V_i,V_i)^{\oplus (d_i-2)} \oplus \Hom(V_i,V_j/V_i)^{\oplus (d_i-2)}
\oplus \bM^{(i)}_{\QQ,\ved}(\bV)
\end{gather}
by the action of
\begin{gather}\label{eq:normalization2}
\GL(V_i) \times \GL(V_j/V_i) \times \prod_{k \neq i,j} G_{d_k}(V_k),
\end{gather}
at the level given by
\begin{gather}\label{eq:normalization3}
-\Bigl( \res_{z=0} \bigl( \lambda_i(z) + \lambda_j(z) \bigr),\,
\res_{z=0} \lambda_j(z),\,
(\lambda_k(z))_{k \neq i,j} \Bigr).
\end{gather}

\subsection{Normalization}\label{subsec:normalization}

The observation in the previous subsection leads us to def\/ine the following:
\begin{Definition}\label{dfn:normalization}
Let $i \in I$ be an irregular pole vertex
of a quiver with multiplicities $(\QQ,\ved)$
and $j \in I$ be the base vertex for $i$.
Then def\/ine $\check{\ved}=(\check{d}_k) \in \Z_{>0}^I$ by
\[
\check{d}_i := 1, \qquad \check{d}_k := d_k \quad \text{for}\hquad k \neq i,
\]
and let $\check{\QQ}=(I,\check{\Omega},\vout,\vin)$
be the quiver obtained from $(\QQ,\ved)$ as the following:
\begin{enm}
\item f\/irst, delete a unique arrow joining $i$ and $j$; then
\item for each arrow $h$ with $\vin(h)=j$,
draw an arrow from $\vout(h)$ to $i$;
\item for each arrow $h$ with $\vout(h)=j$,
draw an arrow from $i$ to $\vin(h)$;
\item f\/inally, draw $d_i-2$ arrows from $j$ to $i$.
\end{enm}
The transformation $(\QQ,\ved) \mapsto (\check{\QQ},\check{\ved})$
is called the {\em normalization} at $i$.
\end{Definition}

The adjacency matrix $\check{\bA}=(\check{a}_{kl})$
of the underlying graph of $\check{\QQ}$ satisf\/ies
\[
\check{a}_{kl} = \check{a}_{lk} =
\begin{cases}
d_i -2 & \text{if}\hquad (k,l)=(i,j), \\
a_{jl} & \text{if}\hquad k=i,\,l \neq j, \\
a_{kl} & \text{if}\hquad k,l \neq i.
\end{cases}
\]

\begin{Example}\label{ex:normalization}
(i) Suppose that $(\QQ,\ved)$
has the graph with multiplicities given below
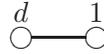
\[
\hfill
\begin{xy}
  \ar@{-} (0,0) *++!D{d} *\cir<4pt>{};
    (10,0)   *++!D{1} *\cir<4pt>{}
\end{xy}
\hfill
\]
Here we assume $d>1$.
The left vertex is an irregular pole, at which we can perform
the normalization and the resulting $(\check{\QQ},\check{\ved})$
has the underlying graph with multiplicities drawn below
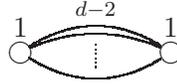
\[
\hfill
\begin{xy}
  \ar@/^10pt/ @{-} (0,0) *++!D{1} *\cir<4pt>{}="A";
    (20,0) *++!D{1} *\cir<4pt>{}="B"
  \ar@/^10pt/ @{-} "A";"B"^{d-2}
  \ar@/^7pt/ @{-} "A";"B"
  \ar@/_10pt/ @{-} "A";"B"
  \ar@{.} (10,1);(10,-2)
\end{xy}
\hfill
\]
The number of edges joining the two vertices are $d-2$.
If $d=3$,
the Kac--Moody algebra associated to $(\QQ,\ved)$ is
of type $G_2$,
while the one associated to $(\check{\QQ},\check{\ved})$
is of type $A_2$.
If $d=4$,
the Kac--Moody algebra associated to $(\QQ,\ved)$ is
of type $A_2^{(2)}$,
while the one associated to $(\check{\QQ},\check{\ved})$
is of type $A_1^{(1)}$.

(ii) Suppose that $(\QQ,\ved)$
has the graph with multiplicities given below
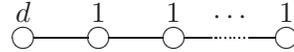
\[
\hfill
\begin{xy}
  \ar@{-} (0,0) *++!D{d} *\cir<4pt>{};
    (10,0)   *++!D{1} *\cir<4pt>{}="C"
  \ar@{-} "C";(20,0) *++!D{1} *\cir<4pt>{}="D"
  \ar@{-} "D";(25,0) \ar@{.} (25,0);(30,0)^*!U{\cdots}
  \ar@{-} (30,0);(35,0) *++!D{1} *\cir<4pt>{}
\end{xy}
\hfill
\]
Here we assume $d>1$ and the number of vertices is $n \geq 3$.
The vertex on the far left is an irregular pole,
at which we can perform the normalization and
the resulting $(\check{\QQ},\check{\ved})$
has the underlying graph with multiplicities drawn below
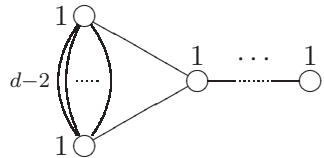
\[
\hfill
\begin{xy}
  \ar@{-} (15,0) *++!D{1} *\cir<4pt>{}="A";
    (0,8.61)   *++!R{1} *\cir<4pt>{}="B"
  \ar@{-} "A";(0,-8.61) *++!R{1} *\cir<4pt>{}="C"
  \ar@/^10pt/ @{-} "B";"C"
  \ar@/_7pt/ @{-} "B";"C"
  \ar@/_10pt/ @{-} "B";"C"_{d-2}
  \ar@{.} (-1,0);(1.7,0)
  \ar@{-} "A";(20,0) \ar@{.} (20,0);(25,0)^*!U{\cdots}
  \ar@{-} (25,0);(30,0) *++!D{1} *\cir<4pt>{}
\end{xy}
\hfill
\]
If $d=2$, then
the Kac--Moody algebra associated to $(\QQ,\ved)$ is
of type $C_n$, while
the one associated to $(\check{\QQ},\check{\ved})$
is of type $A_3$ if $n=3$ and
of type $D_n$ if $n >3$.
If $(d,n)=(3,3)$,
the Kac--Moody algebra associated to $(\QQ,\ved)$ is
of type $D_4^{(3)}$,
while the one associated to $(\check{\QQ},\check{\ved})$
is of type $A_2^{(1)}$.

(iii) Suppose that $(\QQ,\ved)$
has the graph with multiplicities given below
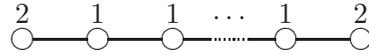
\[
\hfill
\begin{xy}
  \ar@{-} (0,0) *++!D{2} *\cir<4pt>{};
    (10,0)   *++!D{1} *\cir<4pt>{}="A"
  \ar@{-} "A";(20,0) *++!D{1} *\cir<4pt>{}="B"
  \ar@{-} "B";(25,0) \ar@{.} (25,0);(30,0)^*!U{\cdots}
  \ar@{-} (30,0);(35,0) *++!D{1} *\cir<4pt>{}="C"
  \ar@{-} "C";(45,0) *++!D{2} *\cir<4pt>{}
\end{xy}
\hfill
\]
Here the number of vertices is $n \geq 3$.
The associated Kac--Moody algebra is of type $C_{n-1}^{(1)}$.
It has two irregular poles.
Let us perform the normalization at the vertex on the far right.
If $n=3$, the resulting $(\check{\QQ},\check{\ved})$
has the underlying graph with multiplicities drawn below
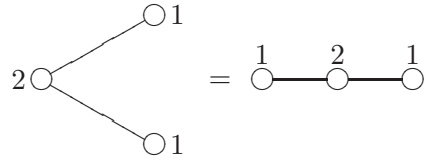
\[
\hfill
\begin{xy}
  \ar@{-} (-15,0) *++!R{2} *\cir<4pt>{}="A";
    (0,8.61)   *++!L{1} *\cir<4pt>{}
  \ar@{-} "A";(0,-8.61) *++!L{1} *\cir<4pt>{}
\end{xy}=
\begin{xy}
  \ar@{-} (0,0) *++!D{1} *\cir<4pt>{};
    (10,0)   *++!D{2} *\cir<4pt>{}="B"
  \ar@{-} "B";(20,0) *++!D{1} *\cir<4pt>{}
\end{xy}
\hfill
\]
The associated Kac--Moody algebra is of type $D_3^{(2)}$.
If $n \geq 4$, the resulting $(\check{\QQ},\check{\ved})$
has the underlying graph with multiplicities drawn below
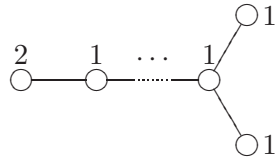
\[
\hfill
\begin{xy}
  \ar@{-} (0,0) *++!D{2} *\cir<4pt>{};
    (10,0)   *++!D{1} *\cir<4pt>{}="B"
    \ar@{-} "B";(15,0) \ar@{.} (15,0);(20,0)^*!U{\cdots}
  \ar@{-} (20,0);(25,0) *++!D{1} *\cir<4pt>{}="C"
  \ar@{-} "C";(30,8.61) *++!L{1} *\cir<4pt>{}
  \ar@{-} "C";(30,-8.61) *++!L{1} *\cir<4pt>{}
\end{xy}
\hfill
\]
The associated Kac--Moody algebra is of type $A_{2n-3}^{(2)}$.
The vertex on the far left is still an irregular pole,
at which we can perform the normalization again.
If $n=4$, the resulting $(\check{\QQ},\check{\ved})$
has the underlying graph with multiplicities drawn below
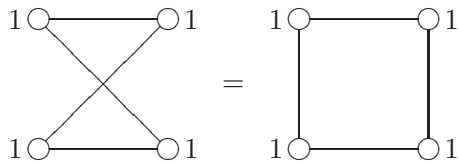
\[
\hfill
\begin{xy}
  \ar@{-} (8.61,-8.61) *++!L{1} *\cir<4pt>{}="A";
          (-8.61,-8.61) *++!R{1} *\cir<4pt>{}="B"
  \ar@{-} (8.61,8.61) *++!L{1} *\cir<4pt>{}="C";
          (-8.61,8.61) *++!R{1} *\cir<4pt>{}="D"
  \ar@{-} "A";"D"
  \ar@{-} "B";"C"
\end{xy}=
\begin{xy}
  \ar@{-} (8.61,-8.61) *++!L{1} *\cir<4pt>{}="A";
          (-8.61,-8.61) *++!R{1} *\cir<4pt>{}="B"
  \ar@{-} (8.61,8.61) *++!L{1} *\cir<4pt>{}="C";
          (-8.61,8.61) *++!R{1} *\cir<4pt>{}="D"
  \ar@{-} "A";"C"
  \ar@{-} "B";"D"
\end{xy}
\hfill
\]
The associated Kac--Moody algebra is of type $A_3^{(1)}$.
If $n > 4$, the resulting $(\check{\QQ},\check{\ved})$
has the underlying graph with multiplicities drawn below
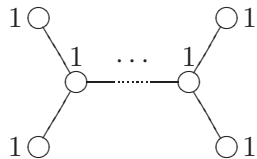
\[
\hfill
\begin{xy}
  \ar@{-} (5,8.61) *++!R{1} *\cir<4pt>{};
    (10,0)   *++!D{1} *\cir<4pt>{}="B"
  \ar@{-} (5,-8.61) *++!R{1} *\cir<4pt>{};"B"
    \ar@{-} "B";(15,0) \ar@{.} (15,0);(20,0)^*!U{\cdots}
  \ar@{-} (20,0);(25,0) *++!D{1} *\cir<4pt>{}="C"
  \ar@{-} "C";(30,8.61) *++!L{1} *\cir<4pt>{}
  \ar@{-} "C";(30,-8.61) *++!L{1} *\cir<4pt>{}
\end{xy}
\hfill
\]
The associated Kac--Moody algebra is of type $D_{n-1}^{(1)}$.
\end{Example}

In the situation discussed in the previous subsection, let
$\check{\bV}=\bigoplus_k \check{V}_k$
be the $I$-graded vector space def\/ined by
\[
\check{V}_j := V_j/V_i, \qquad \check{V}_k := V_k
\quad \text{for}\hquad k \neq j.
\]
Then we see that
the group in \eqref{eq:normalization2}
coincides with $G_{\check{\ved}}(\check{\bV})$.
Furthermore, the following holds:

\begin{Lemma}
The symplectic vector space in \eqref{eq:normalization1}
coincides with $\bM_{\check{\QQ},\check{\ved}}(\check{\bV})$.
\end{Lemma}

\begin{proof}
The def\/initions of $\check{\QQ},\check{\ved}, \check{\bV}$ imply
\begin{gather*}
\Rep_{\check{\QQ}}(\check{\bV}_{\check{\ved}}) =
\bigoplus_{\sumfrac{h \in \check{\Omega}}{h\colon j \to i}}
\Hom(V_j/V_i,V_i)
\\
\phantom{\Rep_{\check{\QQ}}(\check{\bV}_{\check{\ved}}) =}{}   \oplus \bigoplus_{k\neq i,j}\left(
\bigoplus_{\sumfrac{h \in \check{\Omega}}{h\colon k\to i}}
\Hom(V_k \otimes_\C R_{d_k},V_i)
\oplus
\bigoplus_{\sumfrac{h \in \check{\Omega}}{h\colon i\to k}}
\Hom(V_i,V_k \otimes_\C R_{d_k}) \right)
\\
\phantom{\Rep_{\check{\QQ}}(\check{\bV}_{\check{\ved}}) =}{} \oplus \bigoplus_{k\neq i,j}\left(
\bigoplus_{\sumfrac{h \in \check{\Omega}}{h\colon k\to j}}
\Hom(V_k \otimes_\C R_{d_k},V_j/V_i)
\oplus
\bigoplus_{\sumfrac{h \in \check{\Omega}}{h\colon j\to k}}
\Hom(V_j/V_i,V_k \otimes_\C R_{d_k}) \right)
\\
\phantom{\Rep_{\check{\QQ}}(\check{\bV}_{\check{\ved}}) =}{} \oplus \bigoplus_{k,l\neq i,j}
\bigoplus_{\sumfrac{h \in \check{\Omega}}{h\colon k\to l}}
\Hom(V_k \otimes_\C R_{d_k},V_l \otimes_\C R_{d_l})
\\
\phantom{\Rep_{\check{\QQ}}(\check{\bV}_{\check{\ved}})}{} = \Hom(V_j/V_i,V_i)^{\oplus (d_i-2)}
\\
\phantom{\Rep_{\check{\QQ}}(\check{\bV}_{\check{\ved}}) =}{} \oplus \bigoplus_{k\neq i,j}\left(
\bigoplus_{\sumfrac{h \in \Omega}{h\colon k\to j}}
\Hom(V_k \otimes_\C R_{d_k},V_j)
\oplus
\bigoplus_{\sumfrac{h \in \Omega}{h\colon j\to k}}
\Hom(V_j,V_k \otimes_\C R_{d_k}) \right)
\\
\phantom{\Rep_{\check{\QQ}}(\check{\bV}_{\check{\ved}}) =}{} \oplus \bigoplus_{k,l\neq i,j}
\bigoplus_{\sumfrac{h \in \Omega}{h\colon k\to l}}
\Hom(V_k \otimes_\C R_{d_k},V_l \otimes_\C R_{d_l})
\\
\phantom{\Rep_{\check{\QQ}}(\check{\bV}_{\check{\ved}})}{} = \Hom(V_j/V_i,V_i)^{\oplus (d_i-2)} \oplus
\bigoplus_{\sumfrac{h \in \Omega}{\vout(h),\vin(h)\neq i}}
\Hom(V_{\vout(h)} \otimes_\C R_{d_{\vout(h)}},
V_{\vin(h)} \otimes_\C R_{d_{\vin(h)}}).
\end{gather*}
Taking the cotangent bundle,
we thus see that
$\bM_{\check{\QQ},\check{\ved}}(\check{\bV})$
coincides with \eqref{eq:normalization1}.
\end{proof}

Set $\check{\bv} := \dim \check{\bV}$ and
\[
\check{\lambda}=(\check{\lambda}_k(z)) \in R^{\check{\ved}},\qquad
\check{\lambda}_k(z) :=
\begin{cases}
z^{-1}\res\limits_{z=0}\bigl( \lambda_i(z)+\lambda_j(z) \bigr) & \text{if}\hquad k=i, \\
z^{-1}\res\limits_{z=0} \lambda_j(z)  & \text{if}\hquad k=j, \\
\lambda_k(z) & \text{if}\hquad k \neq i,j.
\end{cases}
\]
Then the value given in \eqref{eq:normalization3}
coincides with $-\check{\lambda}$.
Note that
\begin{gather}
\check{\bv} \cdot \res \check{\lambda}
 = v_i \res_{z=0} \bigl( \lambda_i(z) +\lambda_j(z) \bigr)
+ (v_j-v_i) \res_{z=0} \lambda_j(z)
+ \sum_{k \neq i,j} v_k \res_{z=0} \lambda_k(z)
 = \bv \cdot \res \lambda.\label{eq:product}
\end{gather}
Now we state the main result of this section.

\begin{Theorem}\label{thm:normalization}
Let $i \in I$ be an irregular pole vertex
of a quiver with multiplicities $(\QQ,\ved)$
and $j \in I$ be the base vertex for $i$.
Let $(\check{\QQ},\check{\ved})$ be the quiver with multiplicities
obtained by the normalization of $(\QQ,\ved)$ at $i$.
Take $(\lambda, \bv) \in R^\ved \times Q_+$
such that
the top coefficient $\lambda_{i,d_i}$
of $\lambda_i(z)$ is nonzero.
Then the quiver varieties $\Qst_{\QQ,\ved}(\lambda,\bv)$
and $\Qst_{\check{\QQ},\check{\ved}}(\check{\lambda},\check{\bv})$
are symplectomorphic to each other.
\end{Theorem}

\begin{proof}
We have already constructed a bijection between
$\Qset_{\QQ,\ved}(\lambda,\bv)$ and
$\Qset_{\check{\QQ},\check{\ved}}(\check{\lambda},\check{\bv})$.
Thanks to \lemref{lem:boalch2},
in order to prove the assertion
it is suf\/f\/icient to
check that the bijection maps $\Qst_{\QQ,\ved}(\lambda,\bv)$
onto $\Qst_{\check{\QQ},\check{\ved}}(\check{\lambda},\check{\bv})$.
It immediately follows from the three lemmas below.
\end{proof}

\begin{Lemma}\label{lem:normalirred1}
A point $B \in \mu_\ved^{-1}(-\lambda\,\unit_\bv)$
is stable if and only if the corresponding
$(A(z),B_{\neq i})=\Big(\sum\limits_{l=1}^{d_i}A_lz^{-l}, B_{\neq i}\Big)
\in \bO \times \bM^{(i)}_{\QQ,\ved}(\bV)$
satisfies the following condition:
if a collection of subspaces
$S_k \subset V_k \otimes_\C R_{d_k},\, k \neq i$
satisfies
\begin{alignat}{3}
& N_k(S_k) \subset S_k \quad &&\text{for}\hquad k \neq i,j;& \notag\\
& B_h(S_{\vout(h)}) \subset S_{\vin(h)}
\quad & &\text{for}\hquad h \in H
\hquad\text{with}\hquad \vin(h), \vout(h) \neq i;&
\label{eq:normalirred1}\\
& A_l(S_j) \subset S_j \quad && \text{for}\hquad l=1, \dots ,d_i, & \notag
\end{alignat}
then $S_k=0\ (k \neq i)$ or
$S_k =V_k \otimes_\C R_{d_k}\ (k \neq i)$.
\end{Lemma}

\begin{proof}
This is similar to \lemref{lem:refirred}.
First, assume that $B$ is stable and
that a collection of subspaces
$S_k \subset V_k \otimes_\C R_{d_k},\, k \neq i$
satisf\/ies \eqref{eq:normalirred1}.
We def\/ine
\[
S_i := \sum_{l=1}^{d_i} N_i^{l-1} B_{i \slar}(S_j),
\]
and set $\bS := \bigoplus_{k \in I} S_k \subset \bV_\ved$.
Then $N_i(S_i) \subset S_i$, $B_{i \slar}(S_j) \subset S_i$ and
\[
B_{\slar i}(S_i) = \sum_l B_{\slar i}N_i^{l-1}B_{i \slar}(S_j)
= \sum_l A_l(S_j) \subset S_j
\]
imply that $\bS$ is $B$-invariant.
Since $B$ is stable,
we thus have $\bS=0$ or $\bS=\bV_\ved$.

Next assume that the pair $(A(z),B_{\neq i})$ satisf\/ies the condition
in the statement.
Let $\bS=\bigoplus_k S_k$ be a
$B$-invariant subspace of $\bV_\ved$
satisfying $N_k(S_k) \subset S_k$ for all $k \in I$.
Then clearly the collection~$S_k$, $k \neq i$ satisf\/ies
\eqref{eq:normalirred1}, and hence $S_k=0$ $(k\neq i)$ or
$S_k=V_k \otimes_\C R_{d_k}\ (k \neq i)$.
If $S_k=0$, $k\neq i$, we have $B_{\slar i}(S_i)=0$, which implies $S_i=0$
since $\Ker B_{\slar i} \cap \Ker N_i =0$ by \lemref{lem:stability}
and $N_i |_{S_i}$ is nilpotent.
Dualizing the argument, we easily see that
$S_i = V_i \otimes_\C R_{d_i}$ if
$S_k=V_k \otimes_\C R_{d_k}$, $k \neq i$.
\end{proof}

\begin{Lemma}\label{lem:normalirred2}
A point
$B' \in \mu_{\check{\ved}}^{-1}(-\check{\lambda}\,\unit_{\check{\bv}})$
is stable if and only if the corresponding
$(A^0(z),B_{\neq i})=\Big(\sum\limits_{l=2}^{d_i}A^0_lz^{-l},B_{\neq i}\Big)
\in \check{\bO} \times \bM^{(i)}_{\QQ,\ved}(\bV)$
satisfies the following condition:
if an $I$-graded subspace $\bS=\bigoplus_k S_k$ of
$\check{\bV}_{\check{\ved}}=\check{\bV} \otimes_\C R_{\check{\ved}}$ satisfies
\begin{alignat}{3}
& N_k(S_k) \subset S_k \quad &&\text{for}\hquad k \neq i,j; & \notag\\
& B_h(S_{\vout(h)}) \subset S_{\vin(h)}
\quad &&\text{for}\hquad h \in H
\hquad\text{with}\hquad (\vin(h), \vout(h)) \neq (i,j),\,(j,i);&
\label{eq:normalirred2}\\
& A^0_l(S_i \oplus S_j) \subset
S_i \oplus S_j \quad &&\text{for}\hquad l=2, \dots ,d_i, & \notag
\end{alignat}
then $\bS=0$ or $\bS =\check{\bV}_{\check{\ved}}$.
\end{Lemma}

\begin{proof}
In Appendix~\ref{app:misc},
we show that all the block components of $A^0_l$
relative to the decomposition $V_j = \check{V}_i \oplus \check{V}_j$
are described as a (non-commutative) polynomial in $B'_h$
over $h \in H$ with
$(\vin(h),\vout(h))=(i,j)$ or $(j,i)$, and vice versa
(see \remref{rem:depend-B}, where
$A^0$ is denoted by $B$ and $B'_h$
for such $h$ are denoted by $a'_k$, $b'_k$).
Hence an $I$-graded subspace $\bS$ of
$\check{\bV}_{\check{\ved}}$ satisf\/ies~\eqref{eq:normalirred2}
if and only if it is $B'$-invariant and $N_k(S_k) \subset S_k$
for $k \neq i,j$.
\end{proof}

\begin{Lemma}\label{lem:normalirred3}
Let $(A^0(z),B_{\neq i}) \in \check{\bO} \times \bM^{(i)}_{\QQ,\ved}(\bV)$
and let $(A(z),B_{\neq i}) \in \bO \times \bM^{(i)}_{\QQ,\ved}(\bV)$ be
the corresponding pair under the map given in Lemma~{\rm \ref{lem:boalch2}}.
Then $(A^0(z),B_{\neq i})$ satisfies the condition in Lemma~{\rm \ref{lem:normalirred2}}
if and only if $(A(z),B_{\neq i})$ satisfies the one in Lemma~{\rm \ref{lem:normalirred1}}.
\end{Lemma}

\begin{proof}
By def\/inition we have
\[
A(z)=A^0(z)
-z^{-1}\sum_{\sumfrac{\vin(h)=j,}{\vout(h)\neq i}}
\epsilon(h) B_h B_{\ov{h}}
-\lambda_j(z)\,\unit_{V_j},
\]
so the `if' part is clear.
To prove the `only if' part, note that if a collection of subspaces
$S_k \subset V_k \otimes_\C R_{d_k},\, k \neq i$
satisf\/ies \eqref{eq:normalirred1}, then in particular $S_j$
is preserved by the action of
\[
A_{d_i} = A^0_{d_i} = \lambda_{i,d_i}\,\unit_{\check{V}_i}
\oplus 0\, \unit_{\check{V}_j},
\]
and hence is homogeneous relative to the decomposition
$V_j = \check{V}_i \oplus \check{V}_j$;
\[
S_j=(S_j \cap \check{V}_i) \oplus (S_j \cap \check{V}_j).
\]
Now the result immediately follows.
\end{proof}

\subsection{Weyl groups}\label{subsec:weyl}

Let $(\QQ,\ved)$ be a quiver with multiplicities having
an irregular pole vertex $i \in I$ with base $j \in I$,
and let $(\check{\QQ},\check{\ved})$ be the one
obtained by the normalization of $(\QQ,\ved)$ at $i$.
In this subsection we discuss on the relation between the two Weyl groups
associated to $(\QQ,\ved)$ and $(\check{\QQ},\check{\ved})$.

Recall our notation for objects relating to the Kac--Moody algebra;
$\bC = 2 \unit - \bA \bD$ is the generalized Cartan matrix associated to
$(\QQ,\ved)$,
and $\h$, $Q$, $\alpha_k$, $s_k$, the Cartan subalgebra, the root lattice,
the simple roots, and the simple ref\/lections,
of the corresponding Kac--Moody algebra~$\g(\bC)$.
In what follows we denote by
$\check{\bC}$, $\check{\bD}$, $\check{\h}$, $\check{Q}$, $\check{\alpha}_k$,
$\check{s}_k$, the similar objects associated to $(\check{\QQ},\check{\ved})$.

Let $\varphi \colon Q \to \check{Q}$ be the linear map
def\/ined by $\bv \mapsto \check{\bv}=\bv-v_i \check{\alpha}_j$.
The same letter is also used
on the matrix representing $\varphi$ with respect to the simple roots.

\begin{Lemma}
The identity
${}^t \varphi \check{\bD} \check{\bC} \varphi = \bD \bC$ holds.
\end{Lemma}

\begin{proof}
To prove it, we express the matrices in block form with respect to
the decomposition of the index set
$I=\{ i \} \sqcup \{ j \} \sqcup (I \setminus \{ i,j \})$.
First, $\varphi$ is expressed as
\[
\varphi =
\begin{pmatrix}
1  & 0 & 0 \\
-1 & 1 & 0 \\
0  & 0 & \unit
\end{pmatrix}.
\]
By the properties of $i$ and $j$, the matrices
$\bD$ and $\bA$ are respectively expressed as
\[
\bD =
\begin{pmatrix}
d_i & 0 & 0 \\
0 & 1 & 0 \\
0 & 0 & \bD'
\end{pmatrix}, \qquad
\bA =
\begin{pmatrix}
0 & 1 & 0 \\
1 & 0 & {}^t \mathbf{a} \\
0 & \mathbf{a} & \bA'
\end{pmatrix},
\]
where $\bD'$ (resp.\ $\bA'$) is the sub-matrix of $\bD$ (resp.\ $\bA$)
obtained by restricting the index set to $I \setminus \{ i,j \}$,
and $\mathbf{a}=(a_{kj})_{k\neq i,j}$.
By the def\/inition of the normalization,
the matrices $\check{\bD}$ and $\check{\bA}$ are then respectively
expressed as
\[
\check{\bD} =
\begin{pmatrix}
1 & 0 & 0 \\
0 & 1 & 0 \\
0 & 0 & \bD'
\end{pmatrix}, \qquad
\check{\bA} =
\begin{pmatrix}
0 & d_i-2 & {}^t \mathbf{a} \\
d_i-2 & 0 & {}^t \mathbf{a} \\
\mathbf{a} & \mathbf{a} & \bA'
\end{pmatrix}.
\]
Now we check the identity. We have
\begin{gather*}
\bD \bC = 2\bD - \bD \bA \bD
 = \begin{pmatrix}
2d_i & 0 & 0 \\
0 & 2 & 0 \\
0 & 0 & 2\bD'
\end{pmatrix} -
\begin{pmatrix}
0 & d_i & 0 \\
d_i & 0 & {}^t \mathbf{a}\bD' \\
0 & \bD'\mathbf{a} & \bD'\bA'\bD'
\end{pmatrix} \\
\hphantom{\bD \bC = 2\bD - \bD \bA \bD}{}
=
\begin{pmatrix}
2d_i & -d_i & 0 \\
-d_i & 2 & -{}^t \mathbf{a}\bD' \\
0 & -\bD'\mathbf{a} & 2\unit -\bD'\bA'\bD'
\end{pmatrix}.
\end{gather*}
On the other hand,
\begin{gather*}
\check{\bD} \check{\bC} = 2\check{\bD} - \check{\bD} \check{\bA} \check{\bD}
= \begin{pmatrix}
2 & 0 & 0 \\
0 & 2 & 0 \\
0 & 0 & 2\bD'
\end{pmatrix} -
\begin{pmatrix}
0 & d_i-2 & {}^t \mathbf{a}\bD' \\
d_i-2 & 0 & {}^t \mathbf{a}\bD' \\
\bD'\mathbf{a} & \bD'\mathbf{a} & \bD'\bA'\bD'
\end{pmatrix} \\
\hphantom{\check{\bD} \check{\bC} = 2\check{\bD} - \check{\bD} \check{\bA} \check{\bD}}{}
=
\begin{pmatrix}
2 & 2-d_i & -{}^t \mathbf{a}\bD' \\
2-d_i & 2 & -{}^t \mathbf{a}\bD' \\
-\bD'\mathbf{a} & -\bD'\mathbf{a} & 2\unit -\bD'\bA'\bD'
\end{pmatrix}.
\end{gather*}
Hence
\begin{gather*}
{}^t \varphi \check{\bD} \check{\bC} \varphi
=
\begin{pmatrix}
1 & -1 & 0 \\
0 & 1  & 0 \\
0 & 0  & \unit
\end{pmatrix}
\begin{pmatrix}
2 & 2-d_i & -{}^t \mathbf{a}\bD' \\
2-d_i & 2 & -{}^t \mathbf{a}\bD' \\
-\bD'\mathbf{a} & -\bD'\mathbf{a} & 2\unit -\bD'\bA'\bD'
\end{pmatrix}
\begin{pmatrix}
1  & 0 & 0 \\
-1 & 1 & 0 \\
0  & 0 & \unit
\end{pmatrix} \\
\hphantom{{}^t \varphi \check{\bD} \check{\bC} \varphi}{}
=
\begin{pmatrix}
d_i & -d_i & 0 \\
2-d_i & 2 & -{}^t \mathbf{a}\bD' \\
-\bD'\mathbf{a} & -\bD'\mathbf{a} & 2\unit -\bD'\bA'\bD'
\end{pmatrix}
\begin{pmatrix}
1  & 0 & 0 \\
-1 & 1 & 0 \\
0  & 0 & \unit
\end{pmatrix} \\
\hphantom{{}^t \varphi \check{\bD} \check{\bC} \varphi}{}
=
\begin{pmatrix}
2d_i & -d_i & 0 \\
-d_i & 2 & -{}^t \mathbf{a}\bD' \\
0 & -\bD'\mathbf{a} & 2\unit -\bD'\bA'\bD'
\end{pmatrix}
= \bD \bC.\tag*{\qed}
\end{gather*}
\renewcommand{\qed}{}
\end{proof}

The above lemma implies that the map $\varphi$ preserves the inner product.
Furthermore it also implies $\rank \check{\bC} = \rank \bC$,
which means $\dim \h = \dim \check{\h}$.
Thus we can extend $\varphi$ to an isomorphism
$\tilde{\varphi} \colon \h^* \to \check{\h}^*$ preserving the inner product.

Note that by the def\/inition of normalization,
the permutation of the indices $i$ and $j$,
which we denote by $\sigma$,
has no ef\/fect on the matrix $\check{\bC}$.
Hence it def\/ines an involution of $W(\check{\bC})$,
or equivalently, a homomorphism $\Z/2\Z \to \Aut(W(\check{\bC}))$.

\begin{Proposition}\label{prop:weyl}
Under the isomorphism $\tilde{\varphi}$,
the Weyl group $W(\bC)$ associated to $\bC$ is
isomorphic to the semidirect product $W(\check{\bC}) \rtimes \Z/2\Z$
of the one associated to $\check{\bC}$ and $\Z/2\Z$
by the permutation $\sigma$.
\end{Proposition}

\begin{proof}
By the construction of $\tilde{\varphi}$ we have
\[
\tilde{\varphi}(\alpha_k)=
\begin{cases}
\check{\alpha}_i-\check{\alpha}_j & \text{if}\hquad k=i, \\
\check{\alpha}_k & \text{if}\hquad k\neq i.
\end{cases}
\]
As $\tilde{\varphi}$ preserves the inner product,
the above implies that for $k \neq i$,
the map $\tilde{\varphi}s_k\tilde{\varphi}^{-1}$ coincides with
the ref\/lection $\check{s}_k$ relative to
$\tilde{\varphi}(\alpha_k)=\check{\alpha}_k$,
and the map $\tilde{\varphi}s_i\tilde{\varphi}^{-1}$ coincides with
the ref\/lection relative to
$\tilde{\varphi}(\alpha_i)=\check{\alpha}_i -\check{\alpha}_j$.
Note that since the matrix $\check{\bD} \check{\bC}$ is invariant
under the permutation~$\sigma$,
we have
\[
(\check{\alpha}_i+\check{\alpha}_j,\check{\alpha}_i-\check{\alpha}_j)=0,\qquad
(\check{\alpha}_k,\check{\alpha}_i-\check{\alpha}_j)=0,
\qquad  k \neq i,j,
\]
which imply that
$\tilde{\varphi}s_i\tilde{\varphi}^{-1}(\check{\alpha}_k)
=\check{\alpha}_{\sigma(k)}$
for any $k \in I$.
Hence the map
$(\tilde{\varphi}s_i\tilde{\varphi}^{-1})\check{s}_k (\tilde{\varphi}s_i\tilde{\varphi}^{-1})^{-1}$,
which is the ref\/lection relative to
$\tilde{\varphi}s_i\tilde{\varphi}^{-1}(\check{\alpha}_k)$,
coincides with $\check{s}_{\sigma(k)}$ for each~$k$.
Now the result immediately follows.
\end{proof}

We can easily check that
\[
\res (\check{\lambda}) = {}^t\varphi^{-1} \bigl( \res (\lambda) \bigr)
\qquad \text{for}\hquad \lambda \in R^\ved.
\]
Note that
the action of $W(\check{\bC})$ on $R^{\check{\ved}}$
naturally extends to an action of
$W(\check{\bC}) \rtimes \Z/2\Z$.
We see from the above relation
that the map $R^\ved \to R^{\check{\ved}}$, $\lambda \mapsto \check{\lambda}$
is equivariant, and hence so is the map
$R^\ved \times Q \to R^{\check{\ved}} \times \check{Q},\;
(\lambda,\bv) \mapsto (\check{\lambda},\check{\bv})$,
with respect to the isomorphism
$W(\bC) \simeq W(\check{\bC}) \rtimes \Z/2\Z$ given in \propref{prop:weyl}.

\section[Naive moduli of meromorphic connections on $\bbP^1$]{Naive moduli of meromorphic connections on $\boldsymbol{\bbP^1}$}\label{sec:moduli}

This f\/inal section is devoted to study moduli spaces of
meromorphic connections on the trivial bundle over~$\bbP^1$
with some particular type of singularities.

\subsection{Naive moduli}\label{subsec:moduli}

When constructing the moduli spaces of meromorphic connections,
one usually f\/ix the `formal type' of singularities.
However, we f\/ix here the `truncated formal type',
and consider the corresponding `naive' moduli space.
Actually in generic case,
such a naive moduli space gives the moduli space in the usual sense,
which will be explained in~\remref{rem:moduli}.

Fix $n \in \Z_{>0}$ and
\begin{itm}
\item a nonzero f\/inite-dimensional $\C$-vector space $V$;
\item positive integers $k_1, k_2, \dots , k_n$;
\item mutually distinct points $t_1,t_2, \dots ,t_n$ in $\C$.
\end{itm}
Then consider a system
\[
\frac{\rd u}{\rd z} = A(z)\,u(z), \qquad
A(z)=\sum_{i=1}^n \sum_{j=1}^{k_i} \frac{A_{i,j}}{(z-t_i)^j},\qquad
A_{i,j} \in \End(V)
\]
of linear ordinary dif\/ferential equations with rational coef\/f\/icients.
It has a pole at $t_i$ of order at most $k_i$ for each $i$,
and (possibly) a simple pole at $\infty$ with residue $-\sum_i A_{i,1}$.
We identify such a system
with its coef\/f\/icient matrix $A(z)$,
which may be regarded as an element of $\bigoplus_i \g^*_{k_i}(V)$
via $A(z) \mapsto (A_i)$, $A_i(z) :=\sum_k A_{i,j}z^{-j}$.

After Boalch~\cite{MR1864833}, we introduce the following
(the terminologies used here are dif\/ferent from his):
\begin{Definition}
For a system $A(z)=(A_i) \in \bigoplus_i \g^*_{k_i}(V)$
and each $i=1, \dots ,n$,
the $G_{k_i}(V)$-coadjoint orbit through $A_i$
is called the {\em truncated formal type} of~$A(z)$ at~$t_i$.

For given coadjoint orbits
$\bO_i \subset \g^*_{k_i}(V)$, $i=1, \dots ,n$, the set
\[
\Mset(\bO_1, \dots ,\bO_n) :=
\rset{A(z) \in \prod_{i=1}^n \bO_i}%
{\sum_{i=1}^n \res_{z=t_i}A(z)=0}/\GL(V)
\]
is called the {\em naive moduli space of systems
having a pole
of truncated formal type $\bO_i$ at each $t_i$, $i=1,\dots ,n$}.
\end{Definition}

Note that $\prod_i \bO_i$ is a holomorphic symplectic manifold,
and the map
\[
\prod_{i=1}^n \bO_i \to \bigoplus_{i=1}^n \g^*_{k_i}(V), \qquad
A(z) \mapsto \sum_{i=1}^n \res_{z=t_i}A(z)
\]
is a moment map with respect to
the simultaneous $\GL(V)$-conjugation action.
Hence the set $\Mset(\bO_1, \dots , \bO_n)$
is a set-theoretical symplectic quotient.

It is also useful to introduce the following
{\em $\zeta$-twisted} naive moduli space:
\[
\Mset_\zeta(\bO_1, \dots ,\bO_n) :=
\rset{A(z) \in \prod_{i=1}^n \bO_i}%
{\sum_{i=1}^n \res_{z=t_i}A(z)=-\zeta\,\unit_V}/\GL(V)
\quad (\zeta \in \C).
\]

\begin{Definition}\label{dfn:irredsystem}
A system $A(z) \in \bigoplus_i \g^*_{k_i}(V)$
is said to be {\em irreducible} if
there is no nonzero proper subspace $S \subset V$
preserved by all the coef\/f\/icient matrices $A_{i,j}$.
\end{Definition}

If $A(z) \in \bigoplus_i \g^*_{k_i}(V)$ is irreducible,
Schur's lemma shows that
the stabilizer of $A(z)$ with respect to the $\GL(V)$-action
is equal to $\C^\times$,
and furthermore one can show that
the action on the set of irreducible systems in $\prod_i \bO_i$ is proper.

\begin{Definition}
For $\zeta \in \C$,
the holomorphic symplectic manifold
\begin{gather*}
\Mirr_\zeta(\bO_1, \dots, \bO_n) :=
\rset{A(z) \in \prod_{i=1}^n \bO_i}%
{\begin{aligned}
&\text{\rm $A(z)$ is irreducible,} \\
&\sum_i \res_{z=t_i}A(z)=-\zeta\,\unit_V
\end{aligned}
}/\GL(V)
\end{gather*}
is called the {\em $\zeta$-twisted naive moduli space of irreducible systems}
having a pole
of truncated formal type $\bO_i$ at each $t_i$, $i=1,\dots ,n$.
In the $0$-twisted (untwisted) case,
we simply write
$\Mirr_0(\bO_1, \dots ,\bO_n) \equiv \Mirr(\bO_1,\dots ,\bO_n)$.
\end{Definition}

If we have a specif\/ic element $\Lambda_i(z) \in \bO_i$ for each $i$,
the following notation is also useful:
\[
\Mset_\zeta(\Lambda_1, \dots, \Lambda_n) \equiv
\Mset_\zeta(\bO_1, \dots, \bO_n), \qquad
\Mirr_\zeta(\Lambda_1, \dots, \Lambda_n) \equiv
\Mirr_\zeta(\bO_1, \dots, \bO_n).
\]

\begin{Remark}\label{rem:stable-bundle}
Recall that a holomorphic vector bundle with meromorphic connection
$(E,\nabla)$ over a compact Riemann surface
is {\em stable} if for any nonzero proper subbundle $F \subset E$
preserved by $\nabla$,
the inequality $\deg F/\rank F < \deg E/\rank E$ holds.
It is easy to see that if the base space is $\bbP^1$ and
$E$ is trivial, then $(E,\nabla)$ is stable if and only if
it has no nonzero proper {\em trivial} subbundle $F \subset E$
preserved by $\nabla$.
This implies that a system $A(z) \in \bigoplus_i \g^*_{k_i}(V)$
is irreducible if and only if
the associated vector bundle with meromorphic connection
$(\bbP^1 \times V, \rd - A(z)\rd z)$ is stable.
\end{Remark}

\begin{Remark}\label{rem:moduli}
Let us recall a normal form $\Lambda(z)$ introduced in \remref{rem:HTL}.
Assume that each $\Gamma_\lambda$ is {\em non-resonant},
i.e., no two distinct eigenvalues of $\Gamma_\lambda$
dif\/fer by an integer.
Then one can show that an element $A(z) \in \g^*_k(V)$
is equivalent to $\Lambda(z)$ under the coadjoint action
if and only if there is a formal gauge transformation
$g(z) \in \Aut_{\C[[z]]}(\C[[z]] \otimes V)$
which makes $\rd -A(z)\rd z$ into $\rd - \Lambda(z)\rd z$
(see \cite[Remark~18]{0911.3863}).
In this sense the truncated formal type of $\Lambda(z)$
actually prescribe a {\em formal type}.
Hence, if each $\bO_i \subset \g^*_{k_i}(V)$
contains some normal form with non-resonant residue parts,
then the naive moduli space $\Mset(\bO_1,\dots ,\bO_n)$
gives the moduli space of meromorphic connections
on the trivial bundle $\bbP^1 \times V$
having a pole of prescribed formal type at each $t_i$.
\end{Remark}

\subsection{Star-shaped quivers of length one}\label{subsec:star}

In some special case, the naive moduli space
$\Mirr(\bO_1,\dots ,\bO_n)$
can be described as a quiver variety.
Suppose that for each $i=1,\dots ,n$,
the coadjoint orbit $\bO_i$ contains an element of the form
\[
\Xi_i(z)= \begin{pmatrix}
\xi_i(z)\,\unit_{V_i} & 0 \\
0 & \eta_i(z)\,\unit_{V'_i} \end{pmatrix}
\]
for some vector space decomposition $V=V_i \oplus V'_i$
and distinct $\xi_i, \eta_i \in \g^*_{k_i}(\C)$.
Let $d_i$ be the pole order of $\lambda_i := \xi_i - \eta_i$.
Note that $\Xi_i$ is a particular example of normal forms
introduced in \remref{rem:HTL},
and it has non-resonant residue parts
(see \remref{rem:moduli})
if and only if $d_i >1$ or $\res\limits_{z=0} (\xi_i -\eta_i) \notin \Z$.
Also, note that
$\sum\limits_{i=1}^n \tr \res\limits_{z=0} \Xi_i(z) = 0$
is a necessary condition for the non-emptiness of
$\Mset(\Xi_1,\dots ,\Xi_n)=\Mset(\bO_1,\dots ,\bO_n)$.
Indeed, if some $A(z)$ gives a point in $\Mset(\Xi_1,\dots ,\Xi_n)$,
then
\[
0=\sum_{i=1}^n \tr \res_{z=t_i} A(z) =
\sum_{i=1}^n \tr \res_{z=0} \Xi_i(z),
\]
since the function $\tr \circ \res_{z=0} \colon \g^*_{k_i}(V) \to \C$
is invariant under the coadjoint action for each $i$.

Set $I:=\{\, 0,1,\dots ,n \,\}$ and
let $\QQ=(I,\Omega,\vout,\vin)$ be
the `star-shaped quiver with $n$ legs of length one'
as drawn below
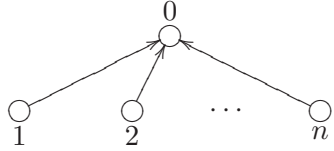
\[
\hfill
\begin{xy}
\ar@{} (25,0);(30,0)_*=0{\cdots}
\ar@{<-} (20,10)*++!D{0}*\cir<4pt>{}="A";(0,0)*++!U{1}*\cir<4pt>{}
\ar@{<-} "A";(15,0)*++!U{2}*\cir<4pt>{}
\ar@{<-} "A";(40,0)*++!U{n}*\cir<4pt>{}
\end{xy}
\hfill
\]
We set $V_0:=V$, $d_0 :=1$, which
together with the above $V_i$, $d_i$ give
an $I$-graded $\C$-vector space $\bV=\bigoplus_i V_i$
and multiplicities $\ved=(d_i) \in \Z_{>0}^I$.
Also we set
\[
\lambda_0(z):= z^{-1} \sum_{i=1}^n \res_{z=0} \eta_i(z)
\in R^1,
\]
which together with the above $\lambda_i$ gives
an element $\lambda=(\lambda_i) \in R^\ved$.
Note that
\begin{gather}
\sum_{i=1}^n \tr \res_{z=0} \Xi_i(z)  =
\sum_{i=1}^n \Bigl[ (\dim V_i) \res_{z=0} \lambda_i(z)
+ (\dim V) \res_{z=0} \eta_i(z) \Bigr] \notag\\
\hphantom{\sum_{i=1}^n \tr \res_{z=0} \Xi_i(z)}{}
=\sum_{i=1}^n (\dim V_i) \res_{z=0} \lambda_i(z)
+ (\dim V) \res_{z=0} \lambda_0(z)
 = \bv \cdot \res \lambda,\label{eq:trace}
\end{gather}
where $\bv:= \dim \bV$.
Hence $\sum\limits_{i=1}^n \tr \res\limits_{z=0} \Xi_i(z) = 0$
if and only if $\bv \cdot \res \lambda =0$,
which is a necessary condition
for the non-emptiness of $\Qset_{\QQ,\ved}(\lambda,\bv)$
(\propref{prop:dim}).

\begin{Proposition}\label{prop:moduli}
There exists a bijection from $\Qset_{\QQ,\ved}(\lambda,\bv)$ to
$\Mset(\Xi_1,\dots ,\Xi_n)$, which maps
$\Qst_{\QQ,\ved}(\lambda,\bv)$ symplectomorphically onto
$\Mirr(\Xi_1,\dots ,\Xi_n)$.
\end{Proposition}

\begin{proof}
Set $\zeta := \res\limits_{z=0} \lambda_0(z)=\sum\limits_{i=1}^n \res\limits_{z=0}\eta_i(z)$.
Then the scalar shift with $\eta_i$ induces a bijection
$\Mset(\Xi_1,\dots ,\Xi_n) \to
\Mset_\zeta(\Lambda_1,\dots ,\Lambda_n)$, where
\begin{gather}\label{eq:type}
\Lambda_i(z) := \Xi_i(z) - \eta_i(z)\,\unit_V =
\begin{pmatrix}
\lambda_i(z)\,\unit_{V_i} & 0 \\
0 & 0\,\unit_{V'_i} \end{pmatrix} \in \g^*_{d_i}(V) \subset \g^*_{k_i}(V),
\end{gather}
and it preserves the irreducibility.
As $\Lambda_i$ has the pole order $d_i$,
the $G_{k_i}(V)$-action on $\Lambda_i$ reduces to the $G_{d_i}(V)$-action
via the natural projection $G_{k_i}(V) \to G_{d_i}(V)$,
so that the orbit $G_{k_i}(V) \cdot \Lambda_i = G_{d_i}(V) \cdot \Lambda_i$
is a $G_{d_i}(V)$-coadjoint orbit.
This replacement of order has no ef\/fect on the naive moduli space.

By the def\/inition of $\QQ$,
we have $\widehat{V}_i=V_0 \otimes_\C R_1 = V$ for each $i >0$
and
\[
\bM_{\QQ,\ved}(\bV) = \bigoplus_{i=1}^n M_i, \qquad
M_i := \Hom(V,V_i \otimes_\C R_{d_i})
\oplus \Hom(V_i \otimes_\C R_{d_i}, V).
\]
Now consider the sets $Z_i \subset M_i$, $i>0$ given in \lemref{lem:orbit1}.
Since the top coef\/f\/icients of $\lambda_i \in \g^*_{d_i}(\C)$,
$i>0$ are nonzero,
\lemref{lem:stability} implies that for each $i>0$,
any point in $\mu_{\ved,i}^{-1}(\lambda_i(z)\,\unit_{V_i})$
satisf\/ies condition~\eqref{eq:stability}.
Hence
\[
\bigcap_{i=1}^n \mu_{\ved,i}^{-1}(\lambda_i(z)\,\unit_{V_i})
= \prod_{i=1}^n Z_i.
\]
Since $\dim V_i \leq \dim V$ for all $i>0$,
Lemmas~\ref{lem:orbit1} and \ref{lem:orbit2} imply that the map
\begin{gather}\label{eq:Phi}
\Phi=(\Phi_i) \colon \ \bM_{\QQ,\ved}(\bV)
\to \bigoplus_{i=1}^n \g^*_{d_i}(V), \qquad
B \mapsto  \bigl(
-B_{\slar i} (z\,\unit_{V_i \otimes R_{d_i}} - N_i)^{-1} B_{i \slar}
\bigr)
\end{gather}
induces a symplectomorphism
\[
\bigcap_{i=1}^n \mu_{\ved,i}^{-1}(\lambda_i(z)\,\unit_{V_i})
/ \prod_{i=1}^n G_{d_i}(V_i)
= \prod_{i=1}^n \left( Z_i/G_{d_i}(V_i) \right)
\xrightarrow{\simeq} \prod_{i=1}^n G_{d_i}(V) \cdot \Lambda_i,
\]
which is clearly $\GL(V)$-equivariant.
Note that
\[
\sum_{i=1}^n \res_{z=0}
\bigl(
-B_{\slar i} (z\,\unit_{V_i \otimes R_{d_i}} - N_i)^{-1} B_{i \slar}
\bigr)
= -\sum_{i=1}^n B_{\slar i} B_{i \slar}
= \res_{z=0} \mu_{\ved,0}(B).
\]
Taking the (set-theoretical) symplectic quotient
by the $\GL(V)$-action at $-\zeta\,\unit_V$,
we thus obtain a bijection from $\Qset_{\QQ,\ved}(\lambda,\bv)$ to
$\Mset_\zeta(\Lambda_1,\dots ,\Lambda_n)$.

The proof of what it maps $\Qst_{\QQ,\ved}(\lambda,\bv)$ onto
$\Mirr_\zeta(\Lambda_1,\dots ,\Lambda_n)$
is quite similar to \lemref{lem:normalirred1}.
First, assume that a point
$B \in \mu_\ved^{-1}(-\lambda\,\unit_{\bV})$ is stable.
Let $\Phi(B) = \Big (\sum\limits_{l=1}^{d_i} A_{i,l}z^{-l}\Big)$,
and assume further that a subspace $S_0 \subset V$ is invariant
under all $A_{i,l}$.
We def\/ine
\[
S_i := \sum_{l=1}^{d_i} N_i^{l-1} B_{i \slar}(S_0), \qquad i>0,
\]
and set $\bS := \bigoplus_{i \in I} S_i \subset \bV_\ved$.
Then $N_i(S_i) \subset S_i$, $B_{i \slar}(S_0) \subset S_i$ and
\[
B_{\slar i}(S_i) = \sum_l B_{\slar i}N_i^{l-1}B_{i \slar}(S_0)
= \sum_l A_{i,l}(S_0) \subset S_0
\]
imply that $\bS$ is $B$-invariant.
Since $B$ is stable,
we thus have $\bS=0$ or $\bS=\bV_\ved$,
and in particular, $S_0=0$ or $S_0=V$,
which shows that the system $\Phi(B)$ is irreducible.

Conversely, assume that the system
$\Phi(B)=(\sum_l A_{i,l}z^{-l})$ is irreducible.
Let $\bS=\bigoplus_i S_i$ be a~$B$-invariant subspace of $\bV_\ved$
satisfying $N_i(S_i) \subset S_i$ for all $i \in I$.
Then $S_0$ is invariant under all~$A_{i,l}$,
and hence $S_0=0$ or $S_0=V$.
If $S_0=0$, then for each $i>0$, we have $B_{\slar i}(S_i)=0$,
which implies $S_i=0$
since $\Ker B_{\slar i} \cap \Ker N_i =0$ by \lemref{lem:stability}
and $N_i |_{S_i}$ is nilpotent.
Dualizing the argument, we easily see that
$S_i = V_i \otimes_\C R_{d_i}$, $i>0$ if $S_0=V$.
Hence $\bS=0$ or $\bS=\bV_\ved$, which shows that
$B$ is stable.
\end{proof}

Conversely, let $\QQ=(I,\Omega,\vout,\vin)$ be as above
and suppose that an $I$-graded $\C$-vector space $\bV=\bigoplus_i V_i$
and multiplicities $\ved=(d_i)$ are given.
Suppose further that they satisfy
$\dim V_i \leq \dim V_0$ and $d_0=1$.
Set $V:=V_0$, and f\/ix a $\C$-vector space $V'_i$ of dimension
$\dim V - \dim V_i$ together with an identif\/ication
$V \simeq V_i \oplus V'_i$ for each $i>0$.
Also, for each $\lambda \in R^\ved$,
set $\zeta := \res\limits_{z=0} \lambda_0$ and let
$\Lambda_i$ be as in \eqref{eq:type}.
Then the above proof also shows that
the map $\Phi$ given in \eqref{eq:Phi} induces a bijection
$\Qset_{\QQ,\ved}(\lambda,\bv) \to \Mset_\zeta(\Lambda_1,\dots ,\Lambda_n)$
mapping $\Qst_{\QQ,\ved}(\lambda,\bv)$ symplectomorphically onto
$\Mirr_\zeta(\Lambda_1,\dots ,\Lambda_n)$.

\subsection{Middle convolution}\label{subsec:middle}

Recall the map given in \eqref{eq:Phi};
\[
\Phi=(\Phi_i) \colon \ \bM_{\QQ,\ved}(\bV)
\to \bigoplus_{i=1}^n \g^*_{d_i}(V), \qquad
B \mapsto
\bigl( -B_{\slar i} (z\,\unit_{V_i \otimes R_{d_i}} - N_i)^{-1} B_{i \slar} \bigr).
\]
Noting $\widehat{V}_0=\bigoplus_{i=1}^n V_i\otimes R_{d_i}$, we set
\[
T := \bigoplus_{i=1}^n (t_i\,\unit_{V_i \otimes R_{d_i}} + N_i)
\in \End(\widehat{V}_0).
\]
Using the natural inclusion
$\iota_i \colon V_i \otimes R_{d_i} \to \widehat{V}_0$ and
projection
$\pi_i \colon \widehat{V}_0 \to V_i \otimes R_{d_i}$,
we then have
\[
(z\,\unit_{\widehat{V}_0} - T)^{-1}
= \sum_{i=1}^n \iota_i (z-t_i-N_i)^{-1} \pi_i
= \sum_{i=1}^n \sum_{j=1}^{k_i} (z-t_i)^{-j} \iota_i N_i^{j-1} \pi_i.
\]
Thus we can write the systems $\Phi(B)$ as
\begin{gather}\label{eq:AHP}
\Phi(B)=-\sum_{i=1}^n B_{\slar i}
\bigl( (z-t_i)\,\unit_{V_i \otimes R_{d_i}} - N_i \bigr)^{-1} B_{i \slar}
=B_{0 \slar} (z\,\unit_{\widehat{V}_0} - T)^{-1}B_{\slar 0}.
\end{gather}
Such an expression of systems has been familiar since
Harnad's work~\cite{MR1309553},
and is in fact quite useful to formulate
the so-called {\em middle convolution}~\cite{0911.3863},
which was originally introduced by Katz~\cite{MR1366651}
for local systems on a punctured $\bbP^1$
and generalized by Arinkin~\cite{0808.0699}
for irregular $\mathscr{D}$-modules.

Let us def\/ine the generalized middle convolution according
to \cite{0911.3863}.
First, we introduce the following fact, which is a ref\/inement of
Woodhouse and Kawakami's observation \cite{MR2287297,kawakami}:
\begin{Proposition}[{\cite[Propositions~1~and~2]{0911.3863}}]\label{prop:datum}
Under the assumption $V \neq 0$,
for any system~$A(z)$
with poles at $t_i$, $i=1,2,\dots ,n$ and possibly a simple pole at $\infty$,
there exists a~quadruple $(W,T,X,Y)$ consisting of
\begin{itm}
\item a finite-dimensional $\C$-vector space $W$;
\item an endomorphism $T \in \End(W)$ with eigenvalues $t_i$, $i=1,2,\dots ,n$;
\item a pair of homomorphisms $(X,Y) \in \Hom(W,V) \oplus \Hom(V,W)$,
\end{itm}
such that
\begin{gather}
X(z\,\unit_W -T)^{-1}Y=A(z), \label{eq:datum-rel}\\
\Ker X_i \cap \Ker N_i =0, \qquad
\range Y_i + \range N_i =V, \label{eq:datum-stable}
\end{gather}
where $N_i$ is the nilpotent part of $T$ restricted on its
generalized $t_i$-eigenspace $W_i := \Ker (T-t_i\,\unit_W)^{\dim W}$,
and $(X_i,Y_i) \in \Hom(W_i,V)\oplus \Hom(V,W_i)$ is the block component
of $(X,Y)$ with respect to the decomposition $W=\bigoplus_i W_i$.
Moreover the choice of $(W,T,X,Y)$
is unique in the following sense: if two quadruples $(W,T,X,Y)$ and
$(W',T',X',Y')$ satisfy \eqref{eq:datum-rel} and \eqref{eq:datum-stable}, then
there exists an isomorphism $f \colon W \to W'$ such that
\[
fTf^{-1}=T', \qquad X = X' f, \qquad f Y = Y'.
\]
\end{Proposition}

The above enables us to def\/ine the middle convolution.
For a system
$A(z)=(A_i) \in \bigoplus_{i=1}^n \g^*_{k_i}(V)$,
take a quadruple $(W,T,X,Y)$ satisfying
\eqref{eq:datum-rel} and \eqref{eq:datum-stable}.
Then for given $\zeta \in \C$,
set $V^\zeta := W/\Ker (YX +\zeta\,\unit_W)$ and let
\begin{itm}
\item $X^\zeta \colon W \to V^\zeta$ be the projection;
\item $Y^\zeta \colon V^\zeta \to W$ be the injection induced from
$YX+\zeta\,\unit_W$.
\end{itm}
Now we def\/ine
\[
\mc_\zeta(A):= X^\zeta(z\,\unit_W - T)^{-1}Y^\zeta
\in \bigoplus_{i=1}^n \g^*_{k_i}(V^\zeta).
\]
By virtue of \propref{prop:datum},
the equivalence class of $\mc_\zeta(A)$ under constant gauge transformations
depends only on that of $A(z)$.
We call it the {\em middle convolution of $A(z)$ with $\zeta$}.\footnote{In \cite{0911.3863}, an explicit construction of
the quadruple $(W,T,X,Y)$ is given so that
the middle convolution $\mc_\zeta(A)$ is well-def\/ined as a system,
not as a gauge equivalence class.}

Let us come back to our situation.
The expression \eqref{eq:AHP} and
\lemref{lem:stability} (which we apply for all $i>0$) imply that
the quadruple $(\widehat{V}_0,T,B_{0 \slar},B_{\slar 0})$
satisf\/ies \eqref{eq:datum-rel} and \eqref{eq:datum-stable} for
$A(z)=\Phi(B)$.
Now assume $\lambda_0(z) \neq 0$ and
consider the middle convolution $\mc_\zeta(A)$ with
$\zeta := \res\limits_{z=0}\lambda_0$.
By the def\/inition,
the triple $(V^\zeta,B_{0 \slar}^\zeta,B_{\slar 0}^\zeta)$
satisf\/ies
\begin{gather}
B_{\slar 0}^\zeta B_{0 \slar}^\zeta =
B_{\slar 0} B_{0 \slar} + \zeta\,\unit_{\widehat{V}_0},
\label{eq:mc-rel}\\
\Ker B_{\slar 0}^\zeta =0, \qquad \range B_{0 \slar}^\zeta =V^\zeta,
\label{eq:mc-stable}
\end{gather}
i.e., it provides a full-rank decomposition of the matrix
$B_{\slar 0} B_{0 \slar} + \zeta\,\unit_{\widehat{V}_0}$.
Recall that such a~triple already appeared
in \secref{sec:reflection}; conditions
\eqref{eq:shift} and \eqref{eq:refstability}
for the $0$-th ref\/lection functor~$\frF_0$ imply that
if we take an $I$-graded $\C$-vector space
$\bV'=\bigoplus_i V'_i$ with $\dim \bV'=s_0(\bv)$ as in
\secref{subsec:proof}
and a representative $B' \in \bM_{\QQ,\ved}(\bV')$
of $\frF_0[B] \in \Qst_{\QQ,\ved}(r_0(\lambda),s_0(\bv))$,
then the triple $(V'_0,B'_{0 \slar},B'_{\slar 0})$
also satisf\/ies \eqref{eq:mc-rel} and \eqref{eq:mc-stable}
(note that $d_0=1$ and $N_0=0$).
By the uniqueness of the full-rank decomposition,
we then see that
there exists an isomorphism $f \colon V^\zeta \to V'_0$
such that $B'_{0 \slar} =f B_{0 \slar}^\zeta$,
$B'_{\slar 0} = B_{\slar 0}^\zeta f^{-1}$,
and hence
\[
\Phi(B')=
B'_{0 \slar} (z\,\unit_{\widehat{V}_0} - T)^{-1}B'_{\slar 0}
= f B_{0 \slar}^\zeta (z\,\unit_{\widehat{V}_0} - T)^{-1}
B_{\slar 0}^\zeta f^{-1} =f \mc_\zeta(A) f^{-1}.
\]
The arguments in the previous subsection
for $\bV', \lambda':=r_0(\lambda)$
show that
$\Phi \colon \bM_{\QQ,\ved}(\bV') \to \bigoplus_{i=1}^n \g^*_{d_i}(V'_0)$
induces a bijection between $\Qset_{\QQ,\ved}(r_0(\lambda),s_0(\bv))$
and $\Mset_{-\zeta}(\Lambda'_1, \dots ,\Lambda'_n)$,
where
\[
\Lambda'_i(z) =
\begin{pmatrix}
\lambda'_i(z)\,\unit_{V_i} & 0 \\
0 & 0\,\unit_{V''_i}
\end{pmatrix} \in \g^*_{d_i}(V'_0),
\qquad V'_0 \simeq V_i \oplus V''_i.
\]
We have now proved the following:

\begin{Proposition}\label{prop:mc-ref}
Let $(\QQ,\ved)$, $\lambda$, $\bv$ be as in Proposition~{\rm \ref{prop:moduli}},
and assume $\zeta := \res\limits_{z=0}\lambda_0$ is nonzero.
Under the above notation, one then has
the following commutative diagram:
\[
\xymatrix{
\Qset_{\QQ,\ved}(\lambda,\bv) \ar[r]^-{\frF_0} \ar[d]_{\Phi}
\ar@{}[rd]|{\circlearrowright} &
\Qset_{\QQ,\ved}(r_0(\lambda),s_0(\bv)) \ar[d]^{\Phi} \\
\Mset_\zeta(\Lambda_1, \dots ,\Lambda_n) \ar[r]_-{\mc_\zeta} &
\Mset_{-\zeta}(\Lambda'_1, \dots ,\Lambda'_n).
}
\]
\end{Proposition}

Next, consider the ref\/lection functors $\frF_i$ for $i>0$.
Let $[B'] = \frF_i[B]$.
Then condition \eqref{eq:refrelother} implies
\[
B'_{\slar j}(z\,\unit_{V_j \otimes R_{d_j}}-N_j)^{-1}B'_{j \slar}
= B_{\slar j}(z\,\unit_{V_j \otimes R_{d_j}}-N_j)^{-1}B_{j \slar},
\qquad j\neq 0,i,
\]
which together with~\eqref{eq:shift} shows that
the two systems $\Phi(B)$ and $\Phi'(B)$ are related via
\[
\Phi(B') =\Phi(B) - \lambda_i(z-t_i)\,\unit_V.
\]

\begin{Proposition}
Let $(\QQ,\ved), \lambda, \bv$ be as in Proposition~{\rm \ref{prop:moduli}},
and set $\zeta := \res\limits_{z=0}\lambda_0$.
For $i =1,2,\dots ,n$, one then has the following commutative diagram:
\[
\xymatrix{
\Qset_{\QQ,\ved}(\lambda,\bv) \ar[rr]^-{\frF_i} \ar[d]_{\Phi}
\ar@{}[rrd]|{\circlearrowright} &&
\Qset_{\QQ,\ved}(r_i(\lambda),s_i(\bv)) \ar[d]^{\Phi} \\
\Mset_\zeta(\Lambda_1, \dots ,\Lambda_n)
\ar[rr]_-{- \lambda_i(z-t_i)\,\unit_V} &&
\Mset_{\zeta+\res\limits_{z=0}\lambda_i}%
(\Lambda_1, \dots ,\Lambda_i -\lambda_i(z)\,\unit_V, \dots ,\Lambda_n),
}
\]
where the bottom horizontal arrow is given by the shift
$A(z) \mapsto A(z) - \lambda_i(z-t_i)\,\unit_V$.
\end{Proposition}

\begin{Remark}\label{rem:harnad}
In \cite{MR1309553}, Harnad considered two meromorphic connections
having the following symmetric description:
\[
\nabla=\rd - \bigl( S + X(z\,\unit_W -T)^{-1}Y \bigr) \rd z, \qquad
\nabla'=\rd + \bigl( T + Y(z\,\unit_V -S)^{-1}X \bigr) \rd z,
\]
where $V$, $W$ are f\/inite-dimensional $\C$-vector spaces,
$S$, $T$ are regular semisimple endomorphisms of $V$, $W$ respectively,
and $(X,Y) \in \Hom(W,V) \oplus \Hom(V,W)$
such that both $(W,T,X,Y)$ and $(V,S,Y,X)$ satisfy \eqref{eq:datum-stable}.
These have an order 2 pole at $z=\infty$
and simple poles at the eigenvalues of~$T$,~$S$ respectively.
He then proved that the isomonodromic deformations of the two systems
are equivalent. After his work, such a duality,
called the {\em Harnad duality},
was established in more general cases
by Woodhouse~\cite{MR2287297}.

Note that if $S=0$, we have $\nabla' = \rd + z^{-1} PQ\, \rd z$.
Hence on the `dual side',
the operation $\mc_\zeta$ corresponds to
just the scalar shift by $z^{-1}\zeta \rd z$.
This interpretation enables us to generalize the middle convolution further;
see~\cite{0911.3863}.
\end{Remark}

\subsection{Examples: rank two cases}\label{subsec:example}

The case $\dim V=2$ is most important because
in this case a generic element in $\g^*_{k_i}(V)$
can be transformed into an element of the form
$\Xi_i(z)=\xi_i(z) \oplus \eta_i(z)$
for some distinct $\xi_i, \eta_i \in \g^*_{k_i}(\C)$.
The dimension of $\Mirr(\Xi_1, \dots ,\Xi_n)$
can be computed as{\samepage
\[
\dim \Mirr(\Xi_1, \dots ,\Xi_n) = \dim \Qst_{\QQ,\ved}(\lambda,\bv)
= 2 -(\bv,\bv)
= 2\sum_{i=1}^n d_i -6,
\]
if it is nonempty.}

First, consider the case $\dim \Mirr(\Xi_1, \dots ,\Xi_n)=0$.
The above formula implies that the tuple $(d_1, \dots ,d_n)$
must be one of the following (up to permutation on indices):
\[
(1,1,1), \quad (2,1), \quad (3).
\]
The corresponding $(\QQ,\ved)$ have
the underlying graphs with multiplicities given by the picture below
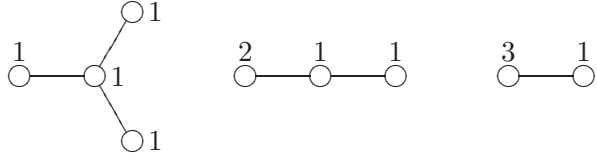
\[
\hfill
\begin{xy}
  \ar@{-} (0,0) *++!D{1} *\cir<4pt>{};
    (10,0)   *++!L{1} *\cir<4pt>{}="A"
  \ar@{-} "A";(15,8.61) *++!L{1} *\cir<4pt>{}
  \ar@{-} "A";(15,-8.61) *++!L{1} *\cir<4pt>{}
  \ar@{-} (30,0) *++!D{2} *\cir<4pt>{};
    (40,0)   *++!D{1} *\cir<4pt>{}="C"
  \ar@{-} "C";(50,0) *++!D{1} *\cir<4pt>{}
  \ar@{-} (65,0) *++!D{3} *\cir<4pt>{};
    (75,0)   *++!D{1} *\cir<4pt>{}
\end{xy}
\hfill
\]
The associated Kac--Moody algebras are respectively given by
\[
D_4, \quad C_3, \quad G_2.
\]
From Example~\ref{ex:normalization},
we see that
the ef\/fect of normalization on these quivers with multiplicities
is given as follows:
\[
D_4 \gets C_4, \quad
C_3 \to A_3, \quad
G_2 \to A_2,
\]
where the arrows represent the process of normalization.

Next consider the case $\dim \Mirr(\Xi_1, \dots ,\Xi_n)=2$.
Then the tuple $(d_1, \dots ,d_n)$
must be one of the following (up to permutation on indices):
\begin{gather}\label{eq:dim2}
(1,1,1,1), \quad (2,1,1), \quad (3,1), \quad (2,2), \quad (4).
\end{gather}
The corresponding $(\QQ,\ved)$ have
the underlying graphs with multiplicities given by the picture below
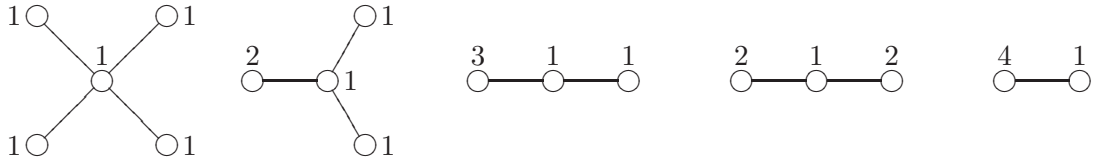
\begin{gather*}
\begin{xy}
  \ar@{-} (0,0) *++!D{1} *\cir<4pt>{}="A";
    (8.61,-8.61)   *++!L{1} *\cir<4pt>{}
  \ar@{-} "A";(-8.61,-8.61) *++!R{1} *\cir<4pt>{}
  \ar@{-} "A";(-8.61,8.61) *++!R{1} *\cir<4pt>{}
  \ar@{-} "A";(8.61,8.61) *++!L{1} *\cir<4pt>{}
  \ar@{-} (20,0) *++!D{2} *\cir<4pt>{};
    (30,0)   *++!L{1} *\cir<4pt>{}="B"
  \ar@{-} "B";(35,8.61) *++!L{1} *\cir<4pt>{}
  \ar@{-} "B";(35,-8.61) *++!L{1} *\cir<4pt>{}
  \ar@{-} (50,0) *++!D{3} *\cir<4pt>{};
    (60,0)   *++!D{1} *\cir<4pt>{}="C"
  \ar@{-} "C";(70,0) *++!D{1} *\cir<4pt>{}
  \ar@{-} (85,0) *++!D{2} *\cir<4pt>{};
    (95,0)   *++!D{1} *\cir<4pt>{}="D"
  \ar@{-} "D";(105,0) *++!D{2} *\cir<4pt>{}
  \ar@{-} (120,0) *++!D{4} *\cir<4pt>{};
    (130,0)   *++!D{1} *\cir<4pt>{}
\end{xy}
\end{gather*}
The associated Kac--Moody algebras are respectively given by
\begin{gather}\label{eq:list}
D_4^{(1)}, \quad A_5^{(2)}, \quad D_4^{(3)}, \quad C_2^{(1)}, \quad A_2^{(2)}.
\end{gather}
From Example~\ref{ex:normalization},
we see that
the ef\/fect of normalization on these quivers with multiplicities
is given as follows:
\begin{gather*}
D_4^{(1)} \gets A_7^{(2)} \gets C_4^{(1)}, \quad
A_3^{(1)} \gets A_5^{(2)} \gets C_3^{(1)}, \quad
D_4^{(3)} \to A_2^{(1)}, \quad
C_2^{(1)} \to D_3^{(2)}, \quad
A_2^{(2)} \to A_1^{(1)},
\end{gather*}
where the arrows represent the process of normalization.
Hence by performing the normalization if necessary,
we obtain the following list of (untwisted) af\/f\/ine Lie algebras:
\[
D_4^{(1)}, \quad A_3^{(1)}, \quad A_2^{(1)}, \quad C_2^{(1)}, \quad A_1^{(1)},
\]
which is well-known as the list
of Okamoto's af\/f\/ine Weyl symmetry groups
of the Painlev\'e equations of type
VI, V, \dots , II,
as mentioned in Introduction.

\begin{Remark}\label{rem:nonempty}
In all the cases appearing in \eqref{eq:dim2},
we can check that $\Mirr(\Xi_1,\dots ,\Xi_n)$ is nonempty
if and only if
$\sum\limits_{i=1}^n \tr \res\limits_{z=0} \Xi_i =0$
(recall that the `only if' part is always true).
We sketch the proof of the `if' part below.

If $(d_1,\dots ,d_n) \neq (2,2)$,
the naive moduli space $\Mirr(\Xi_1,\dots ,\Xi_n)$
is isomorphic to an ordinal quiver variety
$\Qst_\QQ(\zeta,\bv)$ for some extended Dynkin quiver $\QQ$
and $\zeta$, $\bv$ as discussed above.
The formulas \eqref{eq:product} and \eqref{eq:trace} imply
$\zeta \cdot \bv = \sum\limits_{i=1}^n \tr \res\limits_{z=0} \Xi_i$.
Furthermore, since the expected dimension of $\Qst_\QQ(\zeta,\bv)$
is two, we have $(\bv,\bv)=0$,
which implies that $\bv$ is a (positive) imaginary root
(see \cite[Proposition~5.10]{MR1104219}).
In fact, $\bv$ is the minimal positive imaginary root $\delta$
because at least one of its components is equal to one.
It is known~\cite{MR992334} that if $\zeta \cdot \delta =0$,
then $\Qst_\QQ(\zeta,\delta)$
is a deformation of a Kleinian singularity,
which is indeed nonempty\footnote{As a more direct proof,
one can check that if $\zeta \cdot \delta =0$, then
$(\zeta,\delta)$ satisf\/ies the necessary and suf\/f\/icient condition
for the non-emptiness of $\Qst_\QQ(\zeta,\bv)$ given in
\cite[Theorem~1.2]{MR1834739}.}.

Now assume $(d_1,\dots ,d_n)=(2,2)$ and
$\sum\limits_{i=1}^2 \tr \res\limits_{z=0} \Xi_i =0$.
Let $\lambda_i(z)=\lambda_{i,2}z^{-2} + \lambda_{i,1}z^{-1}$,
$\Lambda_i(z)$, $\zeta$ be as in the proof of \propref{prop:moduli},
and for instance, set
\begin{gather*}
A_1(z)  :=
\begin{pmatrix}
2 \lambda_{1,2} & -2\lambda_{1,2} \\
\lambda_{1,2} & -\lambda_{1,2}
\end{pmatrix}z^{-2} +
\begin{pmatrix}
\lambda_{1,1} + \zeta & -\lambda_{1,1} -\zeta \\
\zeta & -\zeta
\end{pmatrix}z^{-1} \\
\phantom{A_1(z) }{}\;
= \begin{pmatrix}
2 & 1 \\ 1 & 1
\end{pmatrix}
\left\{ \Lambda_1(z) +
\begin{pmatrix}
0 & 0 \\
\zeta -\lambda_{1,1} & 0
\end{pmatrix}z^{-1}
\right\}
\begin{pmatrix}
2 & 1 \\ 1 & 1
\end{pmatrix}^{-1},\\
A_2(z) := \Lambda_2(z) + \begin{pmatrix}
0 & \lambda_{1,1} +\zeta \\
-\zeta & 0
\end{pmatrix}z^{-1}.
\end{gather*}
For each $i$, using the assumption $\lambda_{i,2} \neq 0$
and the formula
\[
\begin{pmatrix}
1 & az \\ bz & 1
\end{pmatrix} \cdot \Lambda_i =
\Lambda_i(z) +
\begin{pmatrix}
0 & b\lambda_{i,2} \\ -a\lambda_{i,2} & 0
\end{pmatrix}z^{-1},
\qquad a,b \in \C,
\]
we easily see that $A_i(z)$ is contained in
the $G_2(\C^2)$-coadjoint orbit through $\Lambda_i(z)$.
Furthermore, the assumption $\sum_i \tr \res\limits_{z=0} \Xi_i =0$
implies
$\lambda_{1,1} + \lambda_{2,1} = -2\zeta$,
and hence
\[
\res_{z=0}A_1(z) + \res_{z=0}A_2(z) =
\begin{pmatrix}
\lambda_{1,1}+\lambda_{2,1} + \zeta & 0 \\
0 & -\zeta
\end{pmatrix}
= -\zeta\,\unit_{\C^2}.
\]
The assumption $\lambda_{i,2} \neq 0$ also implies that
the top coef\/f\/icients of $A_1(z)$, $A_2(z)$ have
no common eigenvector, which shows that
the system $(A_1,A_2)$ is irreducible.
Therefore the system
$(A_1+\eta_1(z)\,\unit_{\C^2},A_2+\eta_2(z)\,\unit_{\C^2})$
gives a point in $\Mirr(\Xi_1,\Xi_2)$.
\end{Remark}

\begin{Remark}\label{rem:sasano}
Our list \eqref{eq:list} of Dynkin diagrams
is obtained from Sasano's on~\cite[p.~352]{MR2493874}
by taking the transpose of the generalized Cartan matrices.
It is an interesting problem to ask
the relation between our symmetries and Sasano's.
\end{Remark}

\appendix
\section{Appendix on normalization}\label{app:misc}

In this appendix, we prove Lemmas~\ref{lem:boalch1} and \ref{lem:boalch2}.
Recall the situation discussed in \secref{subsec:shifting};
$i \in I$ is a f\/ixed pole vertex with base $j$, and
$\bO$ is the $G_{d_i}(V_j)$-coadjoint orbit through
\[
\Lambda(z)= \begin{pmatrix}
\lambda_i(z)\,\unit_{V_i} & 0 \\
0 & 0\,\unit_{V_j/V_i} \end{pmatrix},
\qquad V_j \simeq V_i \oplus V_j/V_i,
\]
where the top coef\/f\/icient $\lambda_{i,d_i}$ of $\lambda_i(z)$ is
assumed to be nonzero.
Its `normalized orbit' $\check{\bO}$ is the $B_{d_i}(V_j)$-coadjoint orbit
through the residue-free part $\Lambda^0$ of $\Lambda$.

\subsection{Proof of Lemma~\ref{lem:boalch1}}

We check that
the $B_{d_i}(V_j)$-coadjoint orbit $\check{\bO}$
is invariant under the conjugation action by $K$,
and is $K$-equivariantly
symplectomorphic to the symplectic vector space
\[
\Hom(V_j/V_i,V_i)^{\oplus (d_i-2)} \oplus \Hom(V_i,V_j/V_i)^{\oplus (d_i-2)}.
\]

Note that all the coef\/f\/icients of $\Lambda^0$ are f\/ixed by $K$,
and that the subset $B_{d_i}(V_j) \subset G_{d_i}(V_j)$ is invariant under
the conjugation by constant matrices.
Hence for any $k \in K$ and $g(z) \in B_{d_i}(V_j)$,
\[
k \big(g \cdot \Lambda^0\big) k^{-1} = \big(kgk^{-1}\big) \cdot \big(k\Lambda^0k^{-1}\big)
=\big(kgk^{-1}\big) \cdot \Lambda^0 \in \check{\bO},
\]
i.e., $\check{\bO}$ is invariant under the conjugation by $K$.
Let us calculate the stabilizer of $\Lambda^0(z)$
with respect to the coadjoint $B_{d_i}(V_j)$-action.
Suppose that
$g(z) \in B_{d_i}(V_j)$ stabilizes $\Lambda^0(z)$.
By the def\/inition, we then have
\begin{gather}\label{eq:commute}
g(z)\Lambda^0(z) = \Lambda^0(z)g(z) \mod z^{-1}\gl(V_j)[[z]].
\end{gather}
Write
\[
g(z) =
\begin{pmatrix}
G_{11}(z) & G_{12}(z) \\ G_{21}(z) & G_{22}(z)
\end{pmatrix}
\]
according to the decomposition $V_j = V_i \oplus V_j/V_i$,
and let $\lambda_i^0(z)$ be the residue-free part of $\lambda_i(z)$.
Then
\[
\big[ \Lambda^0(z),g(z) \big] =
\left[
\begin{pmatrix}
\lambda_i^0\,\unit_{V_i} & 0 \\ 0 & 0\,\unit_{V_j/V_i}
\end{pmatrix},
\begin{pmatrix}
G_{11} & G_{12} \\ G_{21} & G_{22}
\end{pmatrix}
\right] =
\begin{pmatrix}
0 & \lambda_i^0G_{12} \\ -\lambda_i^0G_{21} & 0
\end{pmatrix}.
\]
Therefore \eqref{eq:commute} is equivalent to
\[
\lambda_i^0(z)f(z) \in z^{-1}\C[[z]]
\]
for all the matrix entries $f(z)=\sum\limits_{k=1}^{d_i-1}f_kz^k$
of $G_{12}(z)$ and $G_{21}(z)$.
We can write the above condition as
\[
\begin{pmatrix}
\lambda_{i,d_i} & \lambda_{i,d_i-1} & \cdots & \lambda_{i,2} \\
0 & \lambda_{i,d_i} & \cdots & \lambda_{i,3} \\
\vdots & \ddots & \ddots & \vdots \\
0 & \cdots & 0 & \lambda_{i,d_i}
\end{pmatrix}
\begin{pmatrix}
f_{d_i-1} \\
f_{d_i-2} \\
\vdots \\
f_1
\end{pmatrix}
\in \C
\begin{pmatrix}
1 \\
0 \\
\vdots \\
0
\end{pmatrix}.
\]
Since $\lambda_{i,d_i} \neq 0$, this means
$f_k =0$ for all $k=1,2,\dots ,d_i-2$.
Hence the stabilizer is given by
\[
\lset{g(z) = \unit_{V_j} + \sum_{k=1}^{d_i-1} g_k z^k}{%
g_k \in \Lie K,\, k=1,\dots ,d_i-2, \; g_{d_i-1} \in \gl(V_j)}.
\]

The above implies that the orbit $\check{\bO}$ is naturally isomorphic to
\[
\bigl( \gl(V_j)/\Lie K \bigr)^{\oplus (d_i-2)}
\simeq \Hom(V_j/V_i,V_i)^{\oplus (d_i-2)} \oplus \Hom(V_i,V_j/V_i)^{\oplus (d_i-2)}.
\]
Let us denote an element of the vector space on the right hand side
by
\[
(a_1,\dots ,a_{d_i-2},b_1,\dots ,b_{d_i-2}), \qquad
a_k \in \Hom(V_j/V_i,V_i), \qquad b_k \in \Hom(V_i,V_j/V_i),
\]
and set $a(z) := \sum a_k z^k$, $b(z):= \sum_k b_k z^k$.
Then the isomorphism is explicitly given by
\begin{gather}\label{eq:orbit}
(a_k, b_k)_{k=1}^{d_i-2} \longmapsto
g \cdot \Lambda^0 \in \check{\bO},
\qquad
g(z):= \begin{pmatrix}
\unit_{V_i} & a(z) \\
b(z) & \unit_{V_j/V_i}
\end{pmatrix} \in B_{d_i}(V_j).
\end{gather}
It is clearly $K$-equivariant.

Let us calculate
the Kirillov--Kostant--Souriau symplectic form $\omega_{\check{\bO}}$ on $\check{\bO}$
in terms of the coordinates $(a,b)$.
Let $(\delta_l a, \delta_l b)$, $l=1,2$ be two tangent vectors
at $(a,b)$.
Then the corresponding tangent vectors
at $g \cdot \Lambda^0 \in \check{\bO}$ are given by
\[
v_l = \big[\delta_l g \cdot g^{-1}, g \Lambda^0 g^{-1}\big]
\mod z^{-1}\gl(V_j)[[z]]
\quad \in \frb_{d_i}^*(V_j),
\]
where
\[
\delta_l g := \begin{pmatrix}
0 & \delta_l a(z) \\
\delta_l b(z) & 0
\end{pmatrix} \in \frb_{d_i}(V_j), \qquad l=1,2.
\]
By the def\/inition, we have
\begin{gather}
\omega_{\check{\bO}} (v_1,v_2)  = \tr \res_{z=0} \left( g \Lambda^0 g^{-1}
[\delta_1 g \cdot g^{-1}, \delta_2 g \cdot g^{-1}] \right)
 =\tr \res_{z=0} \left( \Lambda^0
[g^{-1}\delta_1 g, g^{-1}\delta_2 g] \right) \notag\\
\phantom{\omega_{\check{\bO}} (v_1,v_2)}{}
= \tr \res_{z=0} \left( [\Lambda^0,
g^{-1}\delta_1 g] g^{-1}\delta_2 g \right). \label{eq:KKS1}
\end{gather}
Using the obvious formula
\begin{gather}\label{eq:inverse}
g(z)^{-1} = \begin{pmatrix}
\unit_{V_i} & a(z) \\
b(z) & \unit_{V_j/V_i}
\end{pmatrix}^{-1}
=\begin{pmatrix}
(\unit_{V_i}-ab)^{-1} & -a(\unit_{V_j/V_i}-ba)^{-1} \\
-b(\unit_{V_i}-ab)^{-1} & (\unit_{V_j/V_i}-ba)^{-1}
\end{pmatrix},
\end{gather}
we have
\begin{gather*}
g^{-1} \delta_1 g  = \begin{pmatrix}
(\unit_{V_i}-ab)^{-1} & -a(\unit_{V_j/V_i}-ba)^{-1} \\
-b(\unit_{V_i}-ab)^{-1} & (\unit_{V_j/V_i}-ba)^{-1}
\end{pmatrix}
\begin{pmatrix}
0 & \delta_1 a(z) \\
\delta_1 b(z) & 0
\end{pmatrix} \\
\phantom{g^{-1} \delta_1 g}{} =
\begin{pmatrix}
(\unit_{V_i}-ab)^{-1}\delta_1 b & (\unit_{V_i}-ab)^{-1}\delta_1 a \\
(\unit_{V_j/V_i}-ba)^{-1}\delta_1 b & -b(\unit_{V_i}-ab)^{-1}\delta_1 a
\end{pmatrix},
\end{gather*}
and hence
\begin{gather*}
[\Lambda^0, g^{-1}\delta_1 g] =
\begin{pmatrix}
0 & \lambda_i^0(\unit_{V_i}-ab)^{-1}\delta_1 a \\
-\lambda_i^0(\unit_{V_j/V_i}-ba)^{-1}\delta_1 b & 0
\end{pmatrix}.
\end{gather*}
Substituting it into \eqref{eq:KKS1},
we obtain
\begin{gather*}
\omega_{\check{\bO}} (v_1,v_2)
= \tr \res_{z=0} \big[ \lambda_i^0(\unit_{V_i}-ab)^{-1}\delta_1 a\,
(\unit_{V_j/V_i}-ba)^{-1}\delta_2 b \big] \\
\hphantom{\omega_{\check{\bO}} (v_1,v_2)=}{}  - \tr \res_{z=0} \big[
\lambda_i^0(\unit_{V_j/V_i}-ba)^{-1}\delta_1 b\,
(\unit_{V_i}-ab)^{-1}\delta_2 a \big],
\end{gather*}
i.e.,
\begin{gather}\label{eq:KKS2}
\omega_{\check{\bO}} = \tr \res_{z=0} \big[ \lambda_i^0(\unit_{V_i}-ab)^{-1}\rd a
\wedge (\unit_{V_j/V_i}-ba)^{-1}\rd b \big].
\end{gather}

Now we set
\begin{gather}\label{eq:basechange}
a'_k := \res_{z=0} \big[ z^k \lambda_i^0 (\unit_{V_i}-ab)^{-1}a \big],
\qquad b'_k := b_k,
\qquad  k=1,\dots ,d_i-2 .
\end{gather}
Using $(\unit_{V_i}-ab)^{-1} = \sum\limits_{l \geq 0} (ab)^l$,
we see that $a'_k$ is the sum of matrices
\[
\lambda_{i,m} (a_{p_1}b_{q_1})(a_{p_2}b_{q_2}) \cdots
(a_{p_l}b_{q_l}) a_r
\]
over all $l \geq 0$ and $m, p_1, \dots ,p_l, q_1, \dots ,q_l, r$
with $m = k+ \sum p_j + \sum q_j + r +1$.
Note that the indices for $a,b$ satisfy
\[
r \leq m - k -1 \leq d_i -k-1, \qquad
p_j, q_j \leq m-k-r-1 < d_i -k-1,
\]
and $r=d_i-k-1$ only when $m=d_i$ and $l=0$.
Thus we can write
\[
a'_k = \lambda_{i,d_i} a_{d_i -k-1} +
f_k(a_1,\dots ,a_{d_i-k-2},b_1,\dots ,b_{d_i-k-2})
\]
for some non-commutative polynomial $f_k$.
Since $\lambda_{i,d_i} \neq 0$,
the above implies that one can uniquely determine $(a_k,b_k)_{k=1}^{d_i-2}$
from $(a'_k,b'_k)_{k=1}^{d_i-2}$ in an algebraic way.
Hence $(a_k,b_k)_{k=1}^{d_i-2} \mapsto (a'_k,b'_k)_{k=1}^{d_i-2}$
is a biregular map.
By the def\/inition, it is clearly $K$-equivariant.

Let us calculate the 1-form
$\sum\limits_{k=1}^{d_i-2}\tr \rd a'_k \wedge \rd b'_k$.
First, we have
\begin{gather*}
\rd \left[ (\unit_{V_i}-ab)^{-1}a \right]
 = \rd (\unit_{V_i}-ab)^{-1} \cdot a + (\unit_{V_i}-ab)^{-1} \rd a \\
\phantom{\rd \left[ (\unit_{V_i}-ab)^{-1}a \right]}{}
= (\unit_{V_i}-ab)^{-1} \rd (ab) (\unit_{V_i}-ab)^{-1}a
+ (\unit_{V_i}-ab)^{-1} \rd a \\
\phantom{\rd \left[ (\unit_{V_i}-ab)^{-1}a \right]}{}
= (\unit_{V_i}-ab)^{-1} \rd a
\big[ b (\unit_{V_i}-ab)^{-1}a + \unit_{V_j/V_i} \big] \\
\phantom{\rd \left[ (\unit_{V_i}-ab)^{-1}a \right]=}{}
+ (\unit_{V_i}-ab)^{-1} a\rd b (\unit_{V_i}-ab)^{-1}a.
\end{gather*}
Note that the obvious equality
$b (\unit_{V_i}-ab) = (\unit_{V_j/V_i}-ba) b$ implies
\[
b (\unit_{V_i}-ab)^{-1} =
(\unit_{V_j/V_i}-ba)^{-1} b.
\]
Thus we have
\begin{align*}
\rd \left[ (\unit_{V_i}-ab)^{-1}a \right]
&= (\unit_{V_i}-ab)^{-1} \rd a
\left[ (\unit_{V_j/V_i}-ba)^{-1} ba + \unit_{V_j/V_i} \right] \\
&\quad + (\unit_{V_i}-ab)^{-1} a\rd b\, (\unit_{V_i}-ab)^{-1}a \\
&=(\unit_{V_i}-ab)^{-1} \rd a\, (\unit_{V_j/V_i}-ba)^{-1} \\
&\quad + (\unit_{V_i}-ab)^{-1} a\rd b\, (\unit_{V_i}-ab)^{-1}a,
\end{align*}
and hence
\begin{gather*}
\tr\left(\lambda_i^0 \rd \left[ (\unit_{V_i}-ab)^{-1}a \right] \wedge \rd b \right)
=\tr\big[ \lambda_i^0 (\unit_{V_i}-ab)^{-1} \rd a
\wedge (\unit_{V_j/V_i}-ba)^{-1} \rd b \big] \\
\phantom{\tr\left(\lambda_i^0 \rd \left[ (\unit_{V_i}-ab)^{-1}a \right] \wedge \rd b \right) =}{}
 + \tr\left[ \lambda_i^0 (\unit_{V_i}-ab)^{-1} a\rd b\,
\wedge (\unit_{V_i}-ab)^{-1}a\rd b \right] \\
\phantom{\tr\left(\lambda_i^0 \rd \left[ (\unit_{V_i}-ab)^{-1}a \right] \wedge \rd b \right)}{}
=\tr\big[ \lambda_i^0 (\unit_{V_i}-ab)^{-1} \rd a
\wedge (\unit_{V_j/V_i}-ba)^{-1} \rd b \big].
\end{gather*}
The above and \eqref{eq:KKS2} imply that the 1-form
$\sum\limits_{k=1}^{d_i-2} \tr \rd a'_k \wedge \rd b'_k$
coincides with $\omega_{\check{\bO}}$; indeed,
\begin{gather*}
\sum_{k=1}^{d_i-2} \tr \rd a'_k \wedge \rd b'_k
=
\sum_{k=1}^{d_i-2} \tr \res_{z=0} \big(
z^k \lambda_i^0 \rd \left[ (\unit_{V_i}-ab)^{-1}a \right]
\wedge \rd b_k \big) \\
\hphantom{\sum_{k=1}^{d_i-2} \tr \rd a'_k \wedge \rd b'_k}{}
= \res_{z=0}\tr \big(
\lambda_i^0 \rd \left[ (\unit_{V_i}-ab)^{-1}a \right]
\wedge \rd b \big) \\
\hphantom{\sum_{k=1}^{d_i-2} \tr \rd a'_k \wedge \rd b'_k}{}
= \res_{z=0}\tr\big[
\lambda_i^0 (\unit_{V_i}-ab)^{-1} \rd a
\wedge (\unit_{V_j/V_i}-ba)^{-1} \rd b \big]  =\omega_{\check{\bO}}.
\end{gather*}

Hence the map $(a'_k,b'_k)_{k=1}^{d_i-2} \mapsto g \cdot \Lambda^0$
is a $K$-equivariant symplectomorphism
\[
\Hom(V_j/V_i,V_i)^{\oplus (d_i-2)} \oplus \Hom(V_i,V_j/V_i)^{\oplus (d_i-2)}
\simeq \check{\bO}.
\]
Since this sends the origin to $\Lambda^0$,
\lemref{lem:boalch1} follows.

\subsection{Proof of Lemma~\ref{lem:boalch2}}

First, we show the following lemma:

\begin{Lemma}\label{lem:residue}
Let
\[
\check{\bO} \to \Lie K, \qquad B(z) \mapsto -\Gamma_B \in \Lie K
\]
be the $K$-moment map sending $\Lambda^0$ to zero.
Then for any $B(z) \in \check{\bO}$,
there exists $g(z) \in B_{d_i}(V_j)$ such that
\[
g(z)B(z)g(z)^{-1} = \Lambda^0(z) + z^{-1}\Gamma_B
\mod \gl(V_j)[[z]].
\]
\end{Lemma}

\begin{proof}
Let $B(z)=\sum\limits_{k=2}^{d_i} B_k z^{-k} \in \check{\bO}$,
and let $a(z),b(z),g(z)$ be as in \eqref{eq:orbit} such that
$B=g \cdot \Lambda^0$.
By the def\/inition of the $B_{d_i}(V_j)$-action, we then have
\begin{gather}\label{eq:residue1}
g(z)^{-1}B(z)g(z) = \Lambda^0(z) + z^{-1}\Gamma
\mod \gl(V_j)[[z]]
\end{gather}
for some $\Gamma \in \gl(V_j)$.
According to the decomposition $V_j = V_i \oplus V_j/V_i$,
we write it as
\[
\Gamma =
\begin{pmatrix}
\Gamma_{11} & \Gamma_{12} \\ \Gamma_{21} & \Gamma_{22}
\end{pmatrix},
\]
and set
\[
\Gamma_B :=
\begin{pmatrix}
\Gamma_{11} & 0 \\ 0 & \Gamma_{22}
\end{pmatrix},\qquad
U :=
\begin{pmatrix}
0 & \lambda_{i,d_i}^{-1}\Gamma_{12} \\ -\lambda_{i,d_i}^{-1}\Gamma_{21} & 0
\end{pmatrix},\qquad
u(z) := \unit_{V_j} + U z^{d_i-1}.
\]
Note that $\Gamma_B \in \Lie K$.
Let $\Lambda_{d_i}$ be the top coef\/f\/icient of $\Lambda^0(z)$.
Then $U$ satisf\/ies
\[
[ \Lambda_{d_i}, U ] =
\begin{pmatrix}
0 & \lambda_{i,d_i} \cdot \lambda_{i,d_i}^{-1}\Gamma_{12} \\
(-\lambda_{i,d_i}) \cdot -\lambda_{i,d_i}^{-1}\Gamma_{21} & 0
\end{pmatrix}
= \Gamma - \Gamma_B,
\]
and hence
\begin{alignat*}{2}
u(z) g(z)^{-1}B(z)g(z)u(z)
&=
u(z) (\Lambda^0(z) + z^{-1}\Gamma)u(z)^{-1}
&\quad &\mod \gl(V_j)[[z]] \\
&= \Lambda^0(z) +z^{-1}\Gamma + z^{-1}[U,\Lambda_{d_i}]
&\quad &\mod \gl(V_j)[[z]] \\
&= \Lambda^0(z) +z^{-1}\Gamma_B
&\quad &\mod \gl(V_j)[[z]].
\end{alignat*}

Now we explicitly describe $\Gamma_B$ in terms of the coordinates
$(a'_k,b'_k)_{k=1}^{d_i-2}$,
which shows that $B \mapsto -\Gamma_B$ is a $K$-moment map.
Note that the constant term of $g(z)$ is the identity,
and hence it acts trivially on $z^{-1}\gl(V_j)[[z]]/\gl(V_j)[[z]]$
by conjugation.
Therefore \eqref{eq:residue1} implies
\begin{alignat*}{2}
B(z) &= g(z)(\Lambda^0(z) +z^{-1}\Gamma)g(z)^{-1}
&\quad &\mod \gl(V_j)[[z]] \\
&= g(z)\Lambda^0(z)g(z)^{-1} + z^{-1}\Gamma
&\quad &\mod \gl(V_j)[[z]].
\end{alignat*}
Substituting \eqref{eq:inverse} into the above equality,
we have
\begin{gather}\label{eq:depend}
B(z) =
\begin{pmatrix}
\lambda_i^0(\unit_{V_i}-ab)^{-1} &
-\lambda_i^0 a(\unit_{V_j/V_i}-ba)^{-1} \\
b\lambda_i^0(\unit_{V_i}-ab)^{-1} &
-b\lambda_i^0 a(\unit_{V_j/V_i}-ba)^{-1}
\end{pmatrix}
+ \frac{\Gamma}{z}
\mod \gl(V_j)[[z]].
\end{gather}
Note that $B(z)$ and $\lambda_i^0(z)$ have no residue parts.
Looking at the block diagonal part of the above and taking the residue,
we thus obtain
\begin{gather*}
\Gamma_{11}  = -\res_{z=0}\big[
\lambda_i^0(\unit_{V_i}-ab)^{-1} \big]
 = -\sum_{l=0}^\infty \res_{z=0}\big[
\lambda_i^0(ab)^l \big]
=-\sum_{l=1}^\infty \res_{z=0}\big[
\lambda_i^0(ab)^l \big] \\
\hphantom{\Gamma_{11}}{}
=-\res_{z=0}\big[
\lambda_i^0(\unit_{V_i}-ab)^{-1}ab \big]
=-\sum_k \res_{z=0}\big[ z^k\,
\lambda_i^0(\unit_{V_i}-ab)^{-1}a \big] b_k
= -\sum_k a'_k b'_k,
\end{gather*}
and similarly,
\[
\Gamma_{22} = \res_{z=0}\big[
b\lambda_i^0 a(\unit_{V_j/V_i}-ba)^{-1}\big]
= \sum_k b'_k a'_k.
\]
Hence
\[
\Gamma_B = -\sum_k
\begin{pmatrix}
a'_k b'_k & 0 \\ 0 & -b'_k a'_k
\end{pmatrix},
\]
which gives the minus of the $K$-moment map vanishing at
$a'_k,b'_k=0$.
\end{proof}

\begin{Remark}\label{rem:depend-gauge}
The matrix $\Gamma$ in~\eqref{eq:residue1} is characterized by
$\Gamma = \res\limits_{z=0} g(z)^{-1}B(z)g(z)$,
so that it depends algebraically on~$a_k$,~$b_k$.
Hence $u(z)g(z)^{-1}$ also depends algebraically on~$a_k$,~$b_k$.
This means that one can choose $g(z)$ in the assertion of
\lemref{lem:residue} so that it depends
algebraically on $B \in \check{\bO}$.
\end{Remark}

\begin{Remark}\label{rem:depend-B}
In the above proof, let us write
\[
B(z) =
\begin{pmatrix}
B_{11}(z) & B_{12}(z) \\ B_{21}(z) & B_{22}(z)
\end{pmatrix}.
\]
Then \eqref{eq:depend} implies
\begin{alignat*}{2}
B_{11}(z) &= \lambda_i^0(\unit_{V_i}-ab)^{-1}
&\quad &\mod z^{-1}\gl(V_i)[[z]], \\
B_{12}(z) &= -\lambda_i^0 a(\unit_{V_j/V_i}-ba)^{-1}
&\quad &\mod \Hom(V_j/V_i,V_i) \otimes z^{-1}\C[[z]], \\
B_{21}(z) &= b\lambda_i^0(\unit_{V_i}-ab)^{-1}
&\quad &\mod \Hom(V_i,V_j/V_i) \otimes z^{-1}\C[[z]], \\
B_{22}(z) &= -b\lambda_i^0 a(\unit_{V_j/V_i}-ba)^{-1}
&\quad &\mod z^{-1}\gl(V_j/V_i)[[z]].
\end{alignat*}
Note that $\lambda_i^0 (\unit_{V_i}-ab)^{-1}a$ has pole order $d_i-1$ and
\[
\lambda_i^0 (\unit_{V_i}-ab)^{-1}a = \sum_{k=1}^{d_i-2} a'_k z^{-k-1}
\mod \Hom(V_j/V_i,V_i) \otimes z^{-1}\C[[z]].
\]
Set $a'(z) :=\sum\limits_{k=1}^{d_i-2} a'_k z^{-k-1}$.
Using the obvious formulas
$a(\unit_{V_j/V_i}-ba)^{-1}=(\unit_{V_i}-ab)^{-1}a$ and
$(\unit_{V_i}-ab)^{-1}=\unit_{V_i}+(\unit_{V_i}-ab)^{-1}ab$,
we can then rewrite the above four equalities as
\begin{alignat}{2}
B_{11}(z) &= \lambda_i^0\unit_{V_i} + a'b'
&\quad &\mod z^{-1}\gl(V_i)[[z]], \label{eq:depend-B-1}\\
B_{12}(z) &= -a', \label{eq:depend-B-2}\\
B_{21}(z) &= \lambda_i^0b' + b'a'b'
&\quad &\mod \Hom(V_i,V_j/V_i) \otimes z^{-1}\C[[z]], \label{eq:depend-B-3}\\
B_{22}(z) &= -b'a'
&\quad &\mod z^{-1}\gl(V_j/V_i)[[z]],\notag
\end{alignat}
which give the explicit description of $B$ in terms of the coordinates
$(a'_k,b'_k)$.
Conversely, we can describe $(a',b')$ in terms of $B$ using the above.
Indeed, \eqref{eq:depend-B-2} determines~$a'$,
and~\eqref{eq:depend-B-1} and~\eqref{eq:depend-B-3} imply
\[
B_{21}(z) = b'(z) B_{11}(z) \mod \Hom(V_i,V_j/V_i) \otimes z^{-1}\C[[z]].
\]
Writing $B_{ij}=\sum_k B_{ij,k} z^{-k}$, we then have
\[
\begin{pmatrix}
B_{21,d_i-1} & \cdots & B_{21,2}
\end{pmatrix}=
\begin{pmatrix}
b'_1 & \cdots & b'_{d_i-2}
\end{pmatrix}
\begin{pmatrix}
B_{11,d_i} & B_{11,d_i-1} & \cdots & B_{11,3} \\
0 & B_{11,d_i} & \cdots & B_{11,4} \\
\vdots & \ddots & \ddots & \vdots \\
0 & \cdots & 0 & B_{11,d_i}
\end{pmatrix}.
\]
Note that \eqref{eq:depend-B-1} also shows
$B_{11,d_i}=\lambda_{i,d_i}\unit_{V_i}$.
Hence the block matrix on the far right is invertible,
and therefore we can express $b'_k$ as
\[
b'_k = \sum_{l=2}^{d_i-1} B_{21,l} F_{lk}(B_{11,3}, \dots ,B_{11,d_i-1})
\]
with some non-commutative polynomial $F_{lk}$.
\end{Remark}

\begin{proof}[Proof of \lemref{lem:boalch2}]
We give a proof of \lemref{lem:boalch2}.
Let $\varphi \colon \check{\bO} \times M \to \g^*_{d_i}(V_j) \times M$
be the map def\/ined in its statement;
\[
\varphi(B(z), x) =(A(z),x), \qquad
A(z) := B(z)-\frac{\mu_M(x)+\zeta\,\unit_{V_j}}{z},
\]
which is clearly equivariant under the conjugation by $K$.
Now suppose that $(B(z),x) \in \check{\bO} \times M$
satisf\/ies the moment map condition
\[
\check{\mu}(B,x):=-\Gamma_B + \mu_M(x)
= -\res_{z=0}\Lambda(z)-\zeta\,\unit_{V_j}.
\]
By \lemref{lem:residue},
there exists $g(z) \in B_{d_i}(V_j)$ such that
\begin{gather}\label{eq:residue2}
B(z) = g(z)\big(\Lambda^0(z) + z^{-1}\Gamma_B\big)g(z)^{-1}
\mod \gl(V_j)[[z]].
\end{gather}
Noting that the constant term $g(0)$ of $g(z)$ is the identity,
we obtain
\begin{alignat*}{2}
A(z)
&=
g(z)\left( \Lambda^0(z) + \frac{\Gamma_B}{z}\right) g(z)^{-1}
-\frac{\mu_M(x)+\zeta\,\unit_{V_j}}{z}
&\quad &\mod \gl(V_j)[[z]] \\
&=
g(z) \Lambda^0(z) g(z)^{-1}
+\frac{\Gamma_B -\mu_M(x) -\zeta\,\unit_{V_j}}{z}
&\quad &\mod \gl(V_j)[[z]] \\
&=
g(z) \Lambda^0(z) g(z)^{-1}
+\frac{\res\limits_{z=0}\Lambda}{z}
&\quad &\mod \gl(V_j)[[z]] \\
&=
g(z) \left( \Lambda^0(z) + \frac{\res\limits_{z=0}\Lambda}{z} \right) g(z)^{-1}
&\quad &\mod \gl(V_j)[[z]] \\
&=g(z)\Lambda(z)g(z)^{-1}
&\quad &\mod \gl(V_j)[[z]],
\end{alignat*}
which implies $A(z) \in \bO$.
Since $B(z)$ has no residue,
we have $\res\limits_{z=0}A(z)=-\mu_M(x) -\zeta\,\unit_{V_j}$,
in other words, the value of the $\GL(V_j)$-moment map
\[
\mu \colon \ \bO \times M \to \gl(V_j),
\qquad (A,x) \mapsto \res_{z=0}A(z) + \mu_M(x)
\]
at $\varphi(B,x)$ is $-\zeta\,\unit_{V_j}$.
Hence $\varphi$ induces a map between the symplectic quotients
\[
\ov{\varphi}\colon \
\check{\mu}^{-1}\Big(-\res_{z=0}\Lambda -\zeta\,\unit_{V_j}\Big)/K
\longrightarrow \mu^{-1}(-\zeta\,\unit_{V_j})/\GL(V_j).
\]

We show that the above map is bijective.
Suppose that
$(B,x), (B',x') \in \check{\mu}^{-1}(-\res\limits_{z=0}\Lambda -\zeta\,\unit_{V_j})$
and $g \in \GL(V_j)$ satisfy $g \cdot \varphi(B,x) =\varphi(B',x')$.
Then $g \cdot x=x'$ and
\begin{gather*}
g \left( B(z)-\frac{\mu_M(x)+\zeta\,\unit_{V_j}}{z} \right) \! g^{-1}  =
B'(z)-\frac{\mu_M(x')+\zeta\,\unit_{V_j}}{z}
 =B'(z)-\frac{g\mu_M(x)g^{-1}\!+\zeta\,\unit_{V_j}}{z}\! \\
\phantom{g \left( B(z)-\frac{\mu_M(x)+\zeta\,\unit_{V_j}}{z} \right) g^{-1}}{}
=B'(z)-g\frac{\mu_M(x)+\zeta\,\unit_{V_j}}{z}g^{-1}.
\end{gather*}
Hence $gB(z)g^{-1}=B'(z)$.
Since $B, B' \in \check{\bO}$, their top coef\/f\/icients
are $\Lambda_{d_i}=\lambda_{i,d_i}\,\unit_{V_i} \oplus 0\,\unit_{V_j/V_i}$,
whose centralizer is $\GL(V_i) \times \GL(V_j/V_i) =K$.
By comparing the top coef\/f\/icients of $gB(z)g^{-1}, B'(z)$,
we thus obtain $g \in K$,
and hence $(B,x)$ and $(B',x')$ lie in the same $K$-orbit.
To prove the surjectivity,
suppose that $(A,x) \in \mu^{-1}(-\zeta\,\unit_{V_j})$ is given.
By using the $\GL(V_j)$-action if necessary, we may assume that
$A=g \cdot \Lambda$ for some $g(z) \in B_{d_i}(V_j)$
(if $A=g \cdot \Lambda$ for $g(z) \in G_{d_i}(V_j)$,
we replace $(A,x)$ with $g(0)^{-1} \cdot (A,x)$).
Let $B(z) \in \frb_{d_i}^*(V_j)$ be the residue-free part of $A(z)$.
Taking modulo $z^{-1}\gl(V_j)[[z]]$ of $A=g\cdot \Lambda$,
we then have $B=g \cdot \Lambda^0 \in \check{\bO}$.
Furthermore, the moment map condition for $(A,x)$
implies
\[
A(z) = B(z) + \frac{\res\limits_{z=0}A}{z}
= B(z)-\frac{\mu_M(x)+\zeta\,\unit_{V_j}}{z}.
\]
Hence $(B,x) =\varphi(A,x)$.
This shows that $\ov{\varphi}$ is surjective.

We have proved that $\ov{\varphi}$ is bijective.
Furthermore, by letting $(B,x)=(B',x')$ in the proof of
the injectivity, we see that
the stabilizer of $\varphi(B,x)$ with respect to the $\GL(V_j)$-action
is contained in that of $(B,x)$ with respect to the $K$-action.
The converse is clear from the $K$-equivariance of $\varphi$,
and hence the two stabilizers coincide.
In particular, free $K$-orbits correspond to free $\GL(V_j)$-orbits
via $\varphi$, which is the second assertion of \lemref{lem:boalch2}.

Finally, we show that $\ov{\varphi}$ preserves the symplectic structure
at points representing free orbits.
Let $(B,x)$ be a point in
the level set $\check{\mu}^{-1}(-\res\limits_{z=0}\Lambda -\zeta\,\unit_{V_j})$
whose stabilizer is trivial (so the level set
is smooth at $(B,x)$), and let $(A,x)=\varphi(B,x)$.
We take $g(z) \in B_{d_i}(V_j)$ satisfying \eqref{eq:residue2} so that
it depends smoothly on $B$,
which is possible as mentioned in \remref{rem:depend-gauge}.
Then the argument just after \eqref{eq:residue2} shows $A=g \cdot \Lambda$,
and furthermore, the smoothness of $g$ implies that
for any tangent vector $(\delta B, v)$ at $(B,x)$,
there exists $\delta g \in \frb_{d_i}(V_j)$ such that
\begin{alignat*}{2}
\delta B &= [\delta g \cdot g^{-1}, B] &\hquad &\mod z^{-1}\gl(V_j)[[z]], \\
\delta A &= [\delta g \cdot g^{-1}, A] &\hquad &\mod \gl(V_j)[[z]],
\end{alignat*}
where $(\delta A, v)=\varphi_*(\delta B,v)$
is the corresponding tangent vector at $(A,x)$.
Now let $(\delta_i B, v_i)$, $i=1,2$ be two tangent vectors
at $(B,x)$ and $\delta_i A, \delta_i g$ as above.
Let $\omega_{\bO}$ (resp.\ $\omega_M$)
be the symplectic form
on $\bO$ (resp.\ $M$).
By the def\/inition, we have
\[
\omega_\bO(\delta_1 A,\delta_2 A)
= \tr \res_{z=0}
\left( A[\delta_1 g \cdot g^{-1}, \delta_2 g \cdot g^{-1}] \right).
\]
Since $\delta_i g$ has no constant term,
we have
$[\delta_1 g \cdot g^{-1}, \delta_2 g \cdot g^{-1}] \in z^2 \gl(V_j)[[z]]$,
which implies
\[
\tr \res_{z=0}
\left( A[\delta_1 g \cdot g^{-1}, \delta_2 g \cdot g^{-1}] \right)
= \tr \res_{z=0}
\left( B[\delta_1 g \cdot g^{-1}, \delta_2 g \cdot g^{-1}] \right)
= \omega_{\check{\bO}} (\delta_1 B, \delta_2 B),
\]
and hence
\[
\omega_{\bO}(\delta_1 A,\delta_2 A) + \omega_M(v_1,v_2)
= \omega_{\check{\bO}}(\delta_1 B,\delta_2 B) + \omega_M(v_1,v_2).
\]
This shows the assertion.
\end{proof}

\subsection*{Acknowledgements}
I am grateful to Philip Boalch for stimulating conversations,
and to Professor Hiraku Nakajima for valuable comments.
This work was supported by the grants
ANR-08-BLAN-0317-01 of the Agence nationale de la recherche
and JSPS Grant-in-Aid for Scientif\/ic Research (S-19104002).

\pdfbookmark[1]{References}{ref}
\LastPageEnding

\end{document}